\documentclass[aop,MSNbibl,dvips]{arximspdf}
\usepackage{mathrsfs}
\usepackage{graphicx}
\doi{10.1214/12-AOP792} %kopijuoti is PTS
\volume{41}
\issue{4}
\pubyear{2013}
\firstpage{2880}
\lastpage{2960}

\makeatletter
\def\mid{|}

\newcommand{\R}{\mathbb{R}}
\newcommand{\K}{\mathbb{K}}
\newcommand{\Z}{\mathbb{Z}}
\newcommand{\N}{\mathbb{N}}
\newcommand{\T}{\mathbb{T}}
\newcommand{\sphere}{\mathbb{S}}
\def\build#1_#2^#3{\mathrel{\mathop{#1}\limits_{#2}^{#3}}}

\def\cal{\mathcal}
\def\ve{\varepsilon}
\def\a{{\cal A}}
\def\m{{\cal M}}
\def\g{{\cal G}}
\def\e{{\cal E}}
\def\b{{\cal B}}
\def\h{{\cal H}}
\def\W{{\cal W}}

\def\mm{\mathbf{m}}
\def\bp{\mathbf{p}}
\def\be{\mathbf{e}}
\def\bT{\mathbf{T}}

\def\ov{\overline}
\def\wh{\widehat}
\def\wt{\widetilde}

\def\bw{\mathbf{w}}
\def\la{\longrightarrow}
\def\da{\downarrow}
\def\noi{\noindent}

\newtheorem{theorem}{Theorem}[section]
\newtheorem{lemma}[theorem]{Lemma}
\newtheorem{proposition}[theorem]{Proposition}
\newtheorem{corollary}[theorem]{Corollary}
\newproclaim{definition}{Definition}[section]

\newproclaim{rema}{Remark}
\newproclaim{remas}{Remarks}
\makeatother

\begin{document}
\begin{frontmatter}

\title{Uniqueness and universality of the Brownian map}
\runtitle{Uniqueness and universality of the Brownian map}

\begin{aug}
\author{\fnms{Jean-Fran\c cois} \snm{Le Gall} \ead[label=e1]{jean-francois.legall@math.u-psud.fr}}
\runauthor{J.-F. Le Gall}
\affiliation{Universit\'e Paris-Sud and Institut universitaire de France}
\address{Math\'ematiques\\
Universit\'e Paris-Sud\\
Campus d'Orsay, 91405 ORSAY\\
France\\
\printead{e1}}
\end{aug}

% HISTORY:
\received{\smonth{5} \syear{2011}}

% ABSTRACT
%
\begin{abstract}
We consider a random planar map $M_n$ which is uniformly distributed
over the class of all rooted $q$-angulations with $n$ faces. We let
$\mm_n$
be the vertex set of $M_n$, which is equipped with the graph distance $d_{\mathrm{gr}}$.
Both when $q\geq4$ is an even integer and when $q=3$, there exists
a positive constant $c_q$ such that the rescaled
metric spaces $(\mm_n,c_qn^{-1/4}d_{\mathrm{gr}})$ converge in distribution
in the Gromov--Hausdorff sense, toward a universal limit called
the Brownian map. The particular case of triangulations
solves a question of Schramm.
\end{abstract}

% KEYWORDS


\end{frontmatter}

\tableofcontents[alignleft,level=2]

%s1 #&#
\section{Introduction}\label{sec1}

In the present work, we derive the convergence in
distribution
in the Gromov--Hausdorff sense
of several important classes of
rescaled random planar maps, toward a
universal limit called the Brownian map. This solves an open
problem that has been stated first by Oded Schramm~\cite{Sch}
in the particular case of triangulations.

Recall that a planar map is a proper embedding of
a finite connected graph in the
two-dimensional sphere, viewed up to orientation-preserving
homeomorphisms of the sphere. Loops and multiple edges are allowed in
the graph.
The faces of the map are
the connected components of the complement of edges, and the degree of
a face counts the number of edges that are incident to it, with the
convention that
if both sides of an edge are incident to the same face, this edge is
counted twice
in the degree of the face. Special cases
of planar maps are triangulations, where each face has degree~$3$,
quadrangulations,
where each face has degree $4,$ and, more generally, $q$-angulations,
where each face has degree $q$. For technical reasons, one often considers
rooted planar maps, meaning that there is a distinguished oriented edge
whose origin is called the root vertex.
Since the pioneering work of Tutte~\cite{Tu}, planar maps have been
studied thoroughly in
combinatorics, and they also arise in other areas of mathematics: See,
in particular, the book of Lando and Zvonkin~\cite{LZ} for algebraic and
geometric motivations.
Large
random planar graphs are of interest in theoretical
physics, where they serve as models of random geometry~\cite{ADJ}, in
particular, in the
theory of two-dimensional quantum gravity.

Let us introduce some notation in order to give a
precise formulation of our main result. Let $q\geq3$
be an integer. We assume that either $q=3$
or $q$ is even. The set of all rooted planar $q$-angulations
with $n$ faces is denoted by $\a^q_n$. For every integer
$n\geq1$ (if $q=3$ we must restrict our attention to even values
of $n$, since $\a^3_n$ is empty if $n$ is odd), we consider a
random planar map $\mathbf{M}_n$ that is uniformly distributed
over $\a^q_n$. We denote the vertex set of $\mathbf{M}_n$
by $\mm_n$. We equip $\mm_n$ with
the graph distance $d_{\mathrm{gr}}$, and we view $(\mm_n,d_{\mathrm{gr}})$ as
a random variable taking values in the space $\K$ of isometry classes
of compact metric spaces. We equip $\K$ with the
Gromov--Hausdorff distance $d_{\mathrm{GH}}$ (see, e.g.,~\cite{BBI}) and note that
$(\K,d_{\mathrm{GH}})$ is a Polish space.

%th1.1 #&#
\begin{theorem}
\label{mainresult}
Set
\[
c_q= \biggl(\frac{9}{q(q-2)} \biggr)^{1/4}
\]
if $q$ is even, and
\[
c_3= 6^{1/4}.
\]
There exists a random compact metric space $(\mm_\infty,D^*)$ called
the Brownian map,
which does not depend on $q$, such that
\[
\bigl(\mm_n, c_q n^{-1/4}d_{\mathrm{gr}}\bigr)
\build{\la}_{n\to\infty}^{\mathrm{(d)}} \bigl(\mm_\infty, D^*\bigr),
\]
where the convergence holds in distribution in the space $\K$.
\end{theorem}

Let us give a precise definition of the Brownian map.
We first need to introduce the random real tree called the CRT,
which can be viewed as the tree coded by a normalized Brownian
excursion, in the following sense. Let $(\be_s)_{0\leq s\leq1}$
be a normalized Brownian
excursion, that is, a positive excursion of linear Brownian motion
conditioned to have duration $1$, and set, for every $s,t\in[0,1]$,
\[
d_\be(s,t)= \be_s+ \be_t - 2
\min_{s\wedge t\leq r\leq s\vee t} \be_r.
\]
Then $d_\be$ is a (random) pseudometric on $[0,1]$, and we consider the
associated equivalence relation $\sim_\be$: for $s,t\in[0,1]$,
\[
s\sim_\be t \quad\mbox{if and only if}\quad d_\be(s,t)=0.
\]
Since $0\sim_\be1$, we may as well view $\sim_\be$ as an equivalence
relation on the
unit circle ${\sphere}^1$.
The CRT is the quotient space ${\cal T}_\be:= {\sphere}^1/ \sim_\be
$, which
is equipped with the
distance induced by $d_\be$. We write $p_\be$ for the canonical projection
from ${\sphere}^1$ onto ${\cal T}_\be$, and $\rho=p_\be(1)$. If
$u,v\in
{\sphere}^1$, we let
$[u,v]$ be the subarc of ${\sphere}^1$ going from $u$
to $v$ in clockwise order, and if $a,b\in{\cal T}_\be$, we define
$[a,b]$ as the image under the canonical projection $p_\be$ of the
smallest subarc
$[u,v]$ of ${\sphere}^1$ such that $p_\be(u)=a$ and $p_\be(v)=b$.
Roughly speaking, $[a,b]$ corresponds to the set
of vertices that one visits when going from $a$ to $b$ around the tree in
clockwise order.

We then introduce Brownian labels on the CRT. We consider
a real-valued process $Z=(Z_a)_{a\in{\cal T}_\be}$ indexed by the CRT,
such that, conditionally on ${\cal T}_\be$, $Z$ is a centered Gaussian process
with $Z_\rho=0$ and $E[(Z_a-Z_b)^2]=d_\be(a,b)$ (this presentation
is slightly informal as we are considering a random process indexed by
a random set, see Section~\ref{CRTlabels} for a more rigorous approach).
We define, for every $a,b\in{\cal T}_\be$,
\[
D^\circ(a,b)= Z_a + Z_b - 2\max \Bigl(
\min_{c\in[a,b]} Z_c, \min_{c\in
[b,a]} Z_c
\Bigr),
\]
and we put $a\simeq b$ if and only if $D^\circ(a,b)=0$. Although this
is not obvious,
it turns out that $\simeq$ is an equivalence
relation on ${\cal T}_\be$, and we let
\[
\mm_\infty:= {\cal T}_\be/ \simeq
\]
be the associated quotient space. We write $\Pi$ for the canonical
projection from
${\cal T}_\be$ onto $\mm_\infty$. We then define the distance on
$\mm_\infty$ by
setting, for every $x,y\in\mm_\infty$,
%
%e1 #&#
\begin{equation}
\label{formulaD} D^*(x,y) = \inf \Biggl\{ \sum_{i=1}^k
D^\circ(a_{i-1},a_i) \Biggr\},
\end{equation}
where the infimum is over all choices of the integer $k\geq1$ and of the
elements $a_0,a_1,\ldots,a_k$ of ${\cal T}_\be$ such that $\Pi(a_0)=x$
and $\Pi(a_k)=y$. It follows from~\cite{IM}, Theorem 3.4, that
$D^*$ is indeed a distance, and the resulting random metric space
$(\mm_\infty,D^*)$ is the Brownian map.

The present work can be viewed
as a continuation and in a sense a conclusion to our
preceding papers~\cite{IM} and~\cite{AM}. In~\cite{IM},
we proved the existence
of sequential Gromov--Hausdorff limits for rescaled uniformly distributed
rooted $2p$-angulations with $n$ faces, and
we called a Brownian map
any random compact metric space that can arise
in such limits (the name Brownian map first appeared in the
work of Marckert and Mokkadem~\cite{MM} which
was dealing with a weak form of the convergence
of rescaled quadrangulations). The main result of~\cite{IM}
used a compactness argument that required the extraction
of suitable subsequences in order to get the desired convergence.
The reason why this extraction was needed is the fact that
the limit could not be characterized completely. It was proved
in~\cite{IM} that any Brownian map can be written in the form $(\mm_\infty, D)$,
where the set $\mm_\infty$ is as described above, and
$D$ is a distance on $\mm_\infty$, for which only upper and lower
bounds were available
in~\cite{IM,AM}.
In particular, the paper~\cite{IM} provided no characterization of the
distance $D$ and it was conceivable that different sequential
limits, or different values of $q$, could lead to different
metric spaces. In the present work, we solve this uniqueness
problem by establishing the explicit formula (\ref{formulaD}), which
had been conjectured in~\cite{IM} and in
a slightly different form in~\cite{MM}. As a consequence, we obtain
the uniqueness of the Brownian map, and we get that this
random metric space is the
scaling limit of uniformly distributed $q$-angulations with
$n$ faces, for the values of $q$ discussed above. Our proofs
strongly depend on the study of geodesics in the
Brownian map that was developed in~\cite{AM}.

At this point, one should mention that the very recent paper of
Miermont~\cite{Mi-quad}
has given another proof of Theorem~\ref{mainresult} in the special
case of
quadrangulations ($q=4$). Our approach was developed
independently of~\cite{Mi-quad} and uses very different ingredients,
leading to more general results. On the other hand, the proof in \cite
{Mi-quad} gives additional
information about the properties of geodesics in $\mm_\infty$, which is
of independent interest.

Let us briefly sketch the main ingredients of our proof in the
bipartite case
where~$q$ is even. From the
main theorem of~\cite{IM}, we can find sequences $(n_k)_{k\geq1}$ of
integers converging to $\infty$ such that the random metric spaces
$(\mm_{n_k},c_qn_k^{-1/4} d_{\mathrm{gr}})$ converge in distribution to $(\mm_\infty,D)$, where
$D$ is a distance on $\mm_\infty$
such that $D\leq D^*$. Additionally, the space $(\mm_\infty,D)$ comes
with a distinguished point $x_*$, which is such that, for every $y\in
\mm_\infty$,
and every $a\in{\cal T}_\be$ such that $\Pi(a)=y$,
\[
D(y,x_*)=D^*(y,x_*)= Z_a- \min Z.
\]
The heart of the proof is now to verify
that $D=D^*$ (Theorem~\ref{MainT} below). To this end, it is enough to
prove that $D(y,y')=D^*(y,y')$ a.s.
when $y$ and $y'$ are distributed uniformly and independently on $\mm_\infty$
(the word uniformly refers to the volume measure on $\mm_\infty$,
which is the image of the normalized Lebesgue measure on $\sphere^1$ under
the projection $\Pi\circ p_\be$). By the results in~\cite{AM}, it
is known that there is an almost surely unique geodesic
path $(\Gamma(t),0\leq t\leq D(y,y'))$ from $y$ to $y'$ in the
metric space $(\mm_\infty, D)$.

Proving that $D(y,y')=D^*(y,y')$ is then essentially equivalent to
verifying that the geodesic $\Gamma$ is well approximated
(in the sense that the lengths of the two paths are not much different)
by another
continuous path going from $y$ to~$y'$, which is constructed by
concatenating pieces of geodesics toward the distinguished
point~$x_*$. To this end, we prove that, for every choice of $r\geq\ve
>0$, and
conditionally on the event $\{D(y,y')\geq r+\ve\}$, the probability
that we have either
%
%e2 #&#
\begin{eqnarray}
\label{cru-intro} D\bigl(x_*,\Gamma(r)\bigr)&=&D\bigl(x_*,\Gamma(r+\ve)\bigr)+\ve\quad \mbox{or}
\nonumber
\\[-8pt]
\\[-8pt]
\nonumber
D\bigl(x_*,\Gamma(r)\bigr)&=&D\bigl(x_*,\Gamma(r-\ve)\bigr)+\ve
\end{eqnarray}
is bounded below by $1-\ve^\beta$ when $\ve$ is small,
where $\beta>0$ is a constant. If (\ref{cru-intro}) holds,
this means that there is a geodesic from
$x_*$ to $\Gamma(r)$ that visits either $\Gamma(r-\ve)$ or $\Gamma
(r+\ve)$
and then coalesces with $\Gamma$. This is of course reminiscent of the
results of~\cite{AM} saying
that any two geodesic paths (starting from arbitrary points of $\mm_\infty$)
ending at a ``typical'' point $x$ of $\mm_\infty$ must coalesce before
hitting $x$.
The difficulty here comes from the fact that interior points of
geodesics are not typical
points of $\mm_\infty$ and so one cannot immediately rely on the
results of
\cite{AM} to establish the preceding estimate (though these results
play a crucial role in the proof).

As in many other papers investigating scaling limits for large random
planar maps,
our proofs make use of bijections between planar maps and various classes
of labeled trees. In the bipartite case, we rely on a bijection discovered
by Bouttier, Di Francesco and Guitter~\cite{BDG} between rooted and pointed
$2p$-angulations with $n$ faces and labeled $p$-trees with $n$ black
vertices (see
Section~\ref{labtrees}, in the case of triangulations we use another
bijection from~\cite{BDG}, which
is presented in Section~\ref{tricoding}). A variant of this
bijection allows us to introduce the notion of a discrete map with
geodesic boundaries (DMGB in short),
which, roughly speaking, corresponds to cutting the map along a
particular discrete
geodesic from the root vertex to the distinguished vertex. This cutting
operation
produces two distinguished geodesics, which are called the boundary
geodesics. The notion
of a DMGB turns out
to play an important role in our proofs and is also of independent interest.
The general philosophy of our approach is that a large planar map can
be obtained
by gluing together many DMGBs along their boundary geodesics.

To complete this introduction, let us mention that the idea of
studying the continuous limit of large random quadrangulations
first appeared in the pioneering paper of Chassaing and Schaeffer~\cite{CS},
which obtained detailed information about the asymptotics of
distances from the root vertex. The results of Chassaing and Schaeffer
were extended to more general classes of random maps in
several papers of Miermont and his coauthors (see, in particular,
\cite{MaMi,MieInvar}), using the bijections with trees found in~\cite{BDG}.
All these results are concerned with the profile of distances from
a particular vertex of the graph and do not provide enough information
to understand Gromov--Hausdorff limits. The understanding of these
limits would be possible if one could compute the asymptotic $k$-point
function, that is, the
asymptotic distribution of the matrix of mutual distances between
$k$ randomly chosen vertices. In the particular case of
quadrangulations, the asymptotic $2$-point function can be derived from the
results of~\cite{CS}, and the asymptotic $3$-point function has been
computed by
Bouttier and Guitter~\cite{BG0}. However, the extension of these
calculations to
higher values of $k$ seems a difficult problem.

As a final remark, Duplantier and Sheffield~\cite{DS} recently
developed a mathematical approach
to two-dimensional quantum gravity based on the Gaussian free field.
It is expected that this approach should be related to the asymptotics
of large planar maps. The very recent paper~\cite{She} contains several
conjectures in this direction. Another very appealing related question
is concerned
with canonical embeddings of the Brownian map: It is known~\cite{LGP}
that the space $(\mm_\infty,D^*)$ is a.s. homeomorphic to the
$2$-sphere ${\sphere}^2$,
and one may look for a canonical construction of a random distance $d$ on
${\sphere}^2$ such
that $(\mm_\infty,D^*)$ is a.s. isometric to $({\sphere}^2,d)$.
The random distance $d$ is expected to have nice conformal invariance
properties. Hopefully these
questions will lead to a promising new line of research in the near future.

The paper is organized as follows. Section~\ref{sec2} recalls basic facts about the
coding of $2p$-angulations by labeled trees, and known results from
\cite{IM} and~\cite{AM} about the
convergence of rescaled $2p$-angulations. Section~\ref{sec3} discusses
discrete maps with geodesic boundaries and their scaling limits. In
Section~\ref{keylemma} we prove
the traversal lemmas, which are concerned with certain properties
of geodesics in large discrete maps with geodesic
boundaries. Roughly speaking, these lemmas provide lower bounds for
the probability that a geodesic path starting from a point of one boundary
geodesic and ending at a point of the other boundary
geodesic will share a significant part of both boundary geodesics.
Section~\ref{mainEstim} proves our main estimate Lemma~\ref{mainest}, which bounds the
probability that (\ref{cru-intro}) does not hold. Section~\ref{sec6}
gives another preliminary estimate relating the distances $D$ and
$D^*$, which comes as
an easy consequence of estimates for the volume of balls proved in
\cite{AM}
(a slightly different approach to the result of Section~\ref{sec6} appears in
\cite{Mi-quad}).
Section~\ref{sec7} contains the proof of Theorem~\ref{mainresult} in the
bipartite case
where $q$ is even. The case of triangulations is treated in
Section~\ref{sec8}, and Section~\ref{sec9} discusses extensions, in particular, to the
Boltzmann distributions on bipartite planar maps considered in \cite
{MaMi}, and open problems. Finally, the \hyperref[app]{Appendix}
provides the proof of two technical lemmas.\vspace*{12pt}

\textsc{Table of notation}.
\mbox{}
%{\small

\noi$M_n$ uniform rooted and pointed $2p$-angulation with $n$ faces

\noi$\mm_n$ vertex set of $M_n$

\noi$d_{\mathrm{gr}}$ graph distance on $\mm_n$

\noi$(\tau_n,(\ell^n_v)_{v\in\tau^\circ_n})$ labeled $p$-tree
associated with $M_n$ via the BDG bijection

\noi$v^n_0,v^n_1,\ldots, v^n_{pn}$ contour sequence of $\tau^\circ_n$

\noi$d_n(i,j)=d_{\mathrm{gr}}(v^n_i,v^n_j)$

\noi$C^n$ contour function of $\tau_n^\circ$

\noi$\Lambda^n$ label function of $(\tau_n,(\ell^n_v)_{v\in\tau
^\circ_n})$

\noi$\gamma_n$ simple geodesic from the first corner of $\varnothing$
in $M_n$ or $\wt M_n$

\noi$\Delta_n=-\min\ell^n +1$

\noi$\wt M_n$ discrete map with geodesic boundaries (DMGB) associated
with $M_n$

\noi$\wt d_{\mathrm{gr}}$ graph distance in $\wt M_n$

\noi$\gamma'_n$ second distinguished boundary geodesic in $\wt M_n$

\noi$\lambda_p,\kappa_p$ scaling constants (cf. Theorem~\ref{mainIM})

\noi$r_n= \lfloor r \kappa_p^{-1} n^{1/4}\rfloor$

\noi$\sigma_n= \min\{i\geq0\dvtx \Lambda^n_i=-r_n\}$

\noi$\ov v_n = v^n_{\sigma_n}$

\noi$\psi_{n,r}(\delta)$ (half) generation of last ancestor of $\ov
v_n$ with label $> -r_n + \delta\kappa_p^{-1}n^{1/4}$

\noi$\Psi_{n,r}(\delta)$ maximal index in the contour sequence of
$\tau^\circ_n$ of this last ancestor

\noi$\be=(\be_t)_{0\leq t\leq1}$ normalized Brownian excursion

\noi${\cal T}_\be=[0,1]/ \sim_\be$ tree coded by $\be$ (CRT)

\noi$p_\be\dvtx  [0,1] \la{\cal T}_\be$ canonical projection

\noi$Z=(Z_t)_{0\leq t\leq1}$ head of Brownian snake driven by $\be$
(Brownian labels on ${\cal T}_\be$)

\noi$ \Delta= - \min Z$

\noi$s_*$ time minimizing $Z$

\noi$\mm_\infty=[0,1]/ \approx\,={\cal T}_\be/ \simeq$ Brownian map

\noi$\Pi\dvtx  {\cal T}_\be\la\mm_\infty$ canonical projection

\noi$\bp= \Pi\circ p_\be$

\noi$D$ distance on the Brownian map derived as scaling limit of graph
distances

\noi$D^\circ(s,t)= Z_s+Z_t -2\max(\min_{[s\wedge t,s\vee t]}Z_r,
\min_{[0,s\wedge t]\cup[s\vee t,1]} Z_r)$

\noi$D^\circ(a,b)=\min\{D^\circ(s,t)\dvtx p_\be(s)=a,p_\be(t)=b\}$
for $
a,b\in{\cal T}_\be$

\noi$D^*(a,b)= \inf\{ \sum_{i=1}^k D^\circ(a_{i-1},a_i)\dvtx a=a_0,a_1,\ldots, a_k=b\}$ for $a,b\in{\cal T}_\be$\vadjust{\goodbreak}

\noi$S_r=\inf\{t\in[0,1]\dvtx Z_t= -r\}$

\noi$S'_r= \sup\{t\in[0,1]\dvtx Z_t=-r\}$

\noi$\Gamma(r)=\bp(S_r)=\bp(S'_r)$ simple geodesic from $\bp(0)$ to
$\bp(s_*)$

\noi$\eta_\delta(r)=\inf\{s> S_r\dvtx \be_s=\min_{t\in[S_r,s]} \be_t\mbox{
and } Z_s=-r+\delta\}$

\noi$\eta'_\delta(r)=\sup\{s<S'_r\dvtx \be_s=\min_{t\in[s,S'_r]} \be_t\mbox
{ and } Z_s=-r+\delta\}$

%}

%s2 #&#
\section{Convergence of rescaled planar maps}\label{sec2}

%s2.1 #&#
\subsection{Labeled $p$-trees}
\label{labtrees}

A plane tree $\tau$ is a finite subset of the set
\[
{\cal U}=\bigcup_{n=0}^\infty
\N^n
\]
of all finite sequences of positive integers (including the empty
sequence $\varnothing$), which satisfies the three following conditions:
\begin{longlist}[(iii)]
\item[(i)] $\varnothing\in\tau$;
\item[(ii)] for every $v=(u_1,\ldots,u_k)\in\tau$
with $k\geq1$, the sequence $(u_1,\ldots,u_{k-1})$
also belongs to $\tau$ [$(u_1,\ldots,u_{k-1}) $ is called the
``parent'' of $v$];
\item[(iii)] for every $v=(u_1,\ldots,u_k)\in\tau$ there
exists an integer $k_v(\tau)\geq0$ such that, for every $j\in\N$,
the vertex $vj:=(u_1,\ldots,u_k,j)$ belongs to $\tau$ if
and only if $1\leq j\leq k_v(\tau)$ [the vertices of the form $vj$ with
$1\leq j\leq k_v(\tau)$
are called the children of $v$].
\end{longlist}
For every $v=(u_1,\ldots,u_k)\in\mathcal{U}$, the generation of $v$ is $|v|=k$.
The notions of an ancestor and a descendant in the tree $\tau$
are defined in an obvious way. By convention a vertex is a descendant
of itself.

Throughout this work, the integer $p\geq2$ is fixed.
A $p$-tree is a plane tree $\tau$ that satisfies the
following additional property:
For every $v\in\tau$ such that $|v|$ is odd, $k_v(\tau)=p-1$.

\begin{figure}

\includegraphics{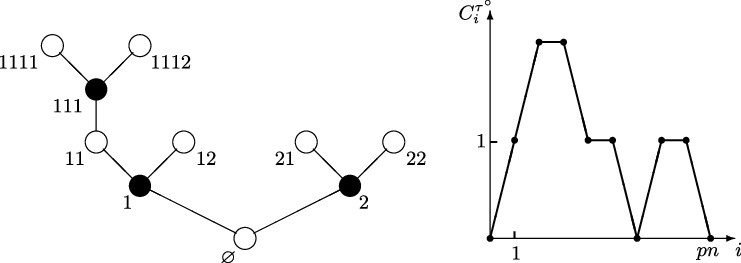}

\caption{A $3$-tree $\tau$ and the associated
contour function $C^{\tau^\circ}$ of $\tau^\circ$.}\label{fig1}
\end{figure}

If $\tau$
is a $p$-tree, vertices $v$ of $\tau$ such that $|v|$ is even are called
white vertices, and vertices $v$ of $\tau$ such that $|v|$ is odd are called
black vertices. We denote the set of all white vertices
of $\tau$ by $\tau^\circ$ and the set of all black vertices by $\tau^\bullet$.
By definition, the size $|\tau|$ of a $p$-tree $\tau$ is the number of
its black vertices.
See the left side of Figure~\ref{fig1} for an example of a $3$-tree.

%f1 #&#

A labeled $p$-tree is a pair $\theta=(\tau,(\ell_v)_{v\in\tau
^\circ})$
that consists of a $p$-tree $\tau$ and a collection of integer labels
assigned to the white vertices of
$\tau$, such that the following properties hold:
\begin{longlist}[(a)]
\item[(a)]$\ell_\varnothing=0$ and $\ell_v\in\Z$ for each $v\in
\tau^\circ$.
\item[(b)] Let $v\in\tau^\bullet$, let $v_{(0)}$ be the parent of $v$
and let
$v_{(j)}=vj$, $1\leq j\leq p-1$, be the children of $v$. Then for every
$j\in\{0,1,\ldots,p-1\}$,
$\ell_{v_{(j+1)}}\geq\ell_{v_{(j)}}-1$, where by convention
$v_{(p)}=v_{(0)}$.
\end{longlist}

Condition (b) means that if one lists the white vertices
adjacent to a given black
vertex in clockwise order, the labels of these vertices can decrease by
at most one
at each step. By definition, the size of $\theta$ is the size of $\tau$.\vadjust{\goodbreak}

Let $\tau$ be a $p$-tree
with $n$ black vertices and let $k=\#\tau-1=pn$. The
depth-first search sequence of $\tau$ is the sequence $w_0,w_1,\ldots,w_{2k}$ of vertices
of $\tau$ which is obtained by induction as follows. First
$w_0=\varnothing$, and
then for every $i\in\{0,\ldots,2k-1\}$, $w_{i+1}$ is either the first
child of
$w_i$ that has not yet appeared in the sequence $w_0,\ldots,w_i$ or
the parent
of $w_i$ if all children of $w_i$ already appear in the sequence
$w_0,\ldots,w_i$.
It is easy to verify that $w_{2k}=\varnothing$ and that all vertices of
$\tau$
appear in the sequence $w_0,w_1,\ldots,w_{2k}$ (some of them
appear more than once).

Vertices $w_i$ are white when $i$ is even
and black when $i$ is odd.
The \textit{contour sequence} of $\tau^\circ$ is by definition the sequence
$v_0,\ldots,v_k$ defined by $v_i=w_{2i}$ for every $i\in\{0,1,\ldots,k\}$.
If $v$ is a given white vertex, each index $i$ such that $v_i=v$
corresponds to a ``corner'' (angular sector) around $v$, and we abusively
speak about the corner $v_i$.

Our limit theorems for random planar maps will be derived from similar
limit theorems
for trees, which are conveniently stated in terms of the coding
functions called
the contour function and the label function. The \textit{contour function}
of $\tau^\circ$ is the discrete sequence
$C^{\tau^\circ}_0,C^{\tau^\circ}_1,\ldots,C^{\tau^\circ}_{pn}$
defined by
\[
C^{\tau^\circ}_i=\tfrac{1}{2} |v_i|\qquad\mbox{for every }0\leq i\leq pn.
\]
See Figure~\ref{fig1} for an example with $p=n=3$.
The \textit{label function} of $\theta=(\tau,(\ell_v)_{v\in\tau^\circ})$
is the discrete sequence $(\Lambda^\theta_0,\Lambda^\theta_1,\ldots,\Lambda^\theta_{pn})$ defined by
\[
\Lambda^\theta_i=\ell_{v_i}\qquad\mbox{for every }0\leq i\leq pn.
\]
From property (b) of the labels and the definition of the contour
sequence, it
is clear that $\Lambda^\theta_{i+1}\geq\Lambda^\theta_i-1$ for every
$0\leq i\leq pn-1$.
The pair $(C^{\tau^\circ},\Lambda^\theta)$ determines~$\theta$ uniquely.

We will need to consider subtrees of a $p$-tree $\tau$
branching from the ancestral line of a given white vertex.
Let $v\in\tau^\circ$, and write $v=v_j$ for some $j\in\{0,1,\ldots,pn\}$
(the choice of $j$ does not matter in what follows). The vertices
$v_i$, $j< i\leq pn$
which are not descendants of $v$
are partitioned into ``subtrees'' that can be described as follows.
First, for every white vertex $u$
that is an ancestor of $v$ distinct of $v$, we can consider the subtree
consisting of $u$ and of its descendants
that belong to the right side of the ancestral line of $v$ (or,
equivalently, that are greater
than $v$ in
lexicographical order). Second,
for every black vertex $w$ that is an ancestor of~$v$, and every child
$u$ of $w$ that is greater than $v$ in
lexicographical order, we can consider the subtree consisting of all
descendants of $u$
(including $u$ itself).
In both cases, this subtree is called a subtree branching from the
right side of the ancestral line of $v$, and
the quantity $\frac{1}{2}|u|$ is called the branching level
of the subtree.
These subtrees can be viewed as $p$-trees, modulo an obvious renaming
of the vertices that preserves the lexicographical order. In the same
way, we can partition the vertices
$v_i$, $0\leq i\leq j$ which are not descendants of $v$ into subtrees branching
from the left side of the ancestral line of~$v$.

If we start from a labeled $p$-tree $\theta=(\tau,(\ell_v)_{v\in
\tau
^\circ})$, we can assign
labels to the white vertices of each subtree in such a way that it
becomes a labeled $p$-tree: just
subtract the label $\ell_u$ of the root $u$ of the subtree from the
label of every vertex in
the subtree.

%s2.2 #&#
\subsection{The Bouttier--Di  Francesco--Guitter bijection}
\label{BDGbij}

Let $\T^p_n$ stand for the set of all labeled $p$-trees with $n$ black
vertices.
We denote the set of all
rooted and pointed $2p$-angulations with $n$ faces by $\m^p_n$. An
element of
$\m^p_n$ is thus a pair $(M,v)$ consisting of a rooted $2p$-angulation
$M\in\a^{2p}_n$ and
a distinguished vertex~$v$. By Euler's formula, the number of choices
for $v$ is
$(p-1)n+2$, independently of $M$.

We now describe the Bouttier--Di  Francesco--Guitter bijection (in short,
the BDG bijection) between $\T^p_n\times\{0,1\}$ and
$\m^{p}_n$.
This bijection can be found in Section 2 of~\cite{BDG} in the more general
setting of bipartite planar maps. Note that~\cite{BDG} deals with pointed
planar maps rather than with rooted and pointed planar maps. However,
the results described below easily follow from~\cite{BDG}
(the bijection we will use is a variant of the one presented in~\cite{IM,AM},
which was concerned with nonpointed rooted $2p$-angulations
and particular labeled $p$-trees called mobiles in~\cite{IM,AM}).

Let $\theta=(\tau,(\ell_v)_{v\in\tau^\circ})\in\T^p_n$ and let
$\ve\in
\{0,1\}$. As previously, we denote the contour
sequence of $\tau^\circ$ by $v_0,v_1,\ldots,v_{pn}$. We extend this sequence
periodically by putting $v_{pn+i}=v_i$ for every $0\leq i\leq pn$.
Suppose that the tree $\tau$ is drawn
on the sphere and add an extra vertex $\partial$.
We associate with the pair $(\theta,\ve)$ a $2p$-angulation $M$
with $n$ faces,
whose set of vertices is
\[
\mm=\tau^\circ\cup\{\partial\}
\]
and whose edges are obtained as follows: For every
$i\in\{0,1,\ldots,pn-1\}$,
\begin{itemize}
\item if $\ell_{v_i}=\min\{\ell_v\dvtx v\in\tau^\circ\}
$, draw
an edge between the corner $v_i$ and $\partial$;
\item if $\ell_{v_i}> \min\{\ell_v\dvtx v\in\tau^\circ\}$,
draw an edge between the corner $v_i$ and the corner $v_j$, where $j$
is the first index in the sequence $i+1,i+2,\ldots,i+pn-1$ such that
$\ell_{v_j}=\ell_{v_i}-1$
(we then say that $j$ is the
successor of $i$, or sometimes that $v_j$ is a successor of $v_i$).
\end{itemize}

\begin{figure}[b]

\includegraphics{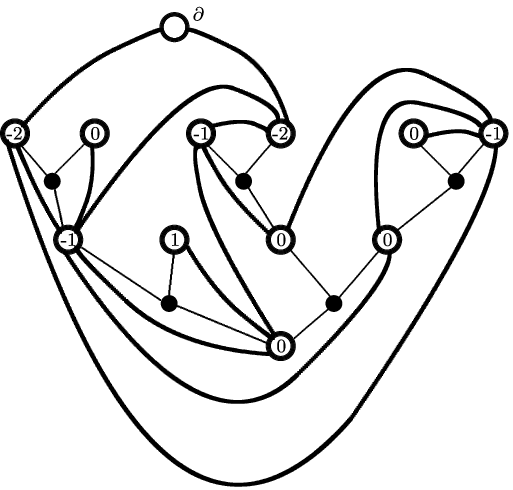}

\caption{A labeled $3$-tree $\theta$ with $5$ black vertices and the
associated $6$-angulation.}\label{fig2}
\end{figure}

Notice that condition (b)
in the definition of a $p$-tree
entails that
$\ell_{v_{i+1}}\geq\ell_{v_i}-1$ for every $i\in\{0,1,\ldots,pn-1\}$.
This ensures that whenever $\ell_{v_i}> \min\{\ell_v\dvtx v\in\tau^\circ\}
$ there is at least one
vertex among $v_{i+1},v_{i+2},\ldots,v_{i+pn-1}$ with label $\ell_{v_i}-1$.
The construction can be made in a unique way (up to
orientation-preserving homeomorphisms of the sphere)
if we impose
that edges of the map do not intersect, except possibly at their
endpoints, and
do not intersect the edges of the tree. We refer to
Section 2 of~\cite{BDG} for a more detailed description (here we will
only need the fact that edges are
generated in the way described above). The resulting planar
map $M$ is a $2p$-angulation. By definition, this $2p$-angulation is
rooted at the edge between
vertex $\varnothing$ and its successor $w=v_j$, where $j=\min\{i\in\{
1,\ldots,pn\}\dvtx \ell_i=-1\}$, and by
convention $w=\partial$ if $\min\{\ell_v\dvtx v\in\tau^\circ\}=0$. The
orientation of this edge
is specified by the variable $\ve$: if $\ve=1$, the root vertex is
$\varnothing$
and if $\ve=0$, the root vertex is~$w$. Finally, the $2p$-angulation
$M$ is pointed at the vertex $\partial$, so that we have indeed
obtained a rooted
and pointed $2p$-angulation.
Each face of $M$ contains exactly one black
vertex of
$\tau$ (see Figure~\ref{fig2}).

%f2 #&#

The preceding construction yields
a bijection from the set $\T^p_n\times\{0,1\}$ onto $\m^{p}_n$,
which is called the BDG bijection. Figure~\ref{fig2} gives an example of a
labeled $3$-tree with $5$ black vertices
(the numbers appearing
inside the circles representing white vertices are the labels assigned to
these vertices)
and shows the $6$-angulation with~$5$ faces associated with this $3$-tree
via the BDG bijection.

The
following property, which relates labels on the tree $\tau^\circ$ to distances
in the planar map $M$, plays a key role. As previously, we write
$d_{\mathrm{gr}}$ for the graph distance in the vertex
set $\mm$. Then, for every vertex $v\in\tau^\circ$, we have
%
%e3 #&#
\begin{equation}
\label{formuladgr} d_{\mathrm{gr}}(\partial, v) = \ell_v - \min
\bigl\{\ell_w\dvtx w\in\tau^\circ\bigr\} +1.
\end{equation}

If $v$ and $v'$ are two arbitrary vertices of $M$, there is no such
simple expression
for $d_{\mathrm{gr}}(v,v')$ in terms of the labels on $\tau^\circ$. However, the
following bound
is useful. Suppose that $v=v_i$ and $v'=v_j$ for some $i,j\in\{
1,\ldots,pn\}$ with $i< j$. Then,
%
%e4 #&#
\begin{equation}
\label{bounddgr} d_{\mathrm{gr}}\bigl(v,v'\bigr) \leq
\ell_{v_i} + \ell_{v_j} -2 \max \Bigl(\min_{i\leq
k\leq
j}
\ell_{v_k}, \min_{j\leq k\leq i+pn} \ell_{v_k} \Bigr)+2.
\end{equation}
See~\cite{IM}, Lemma 3.1, for a proof
in a slightly different context, which is easily adapted. This proof
makes use of simple geodesics, which are defined as follows. Let $v\in
\tau^\circ$,
and let $i\in\{0,1,\ldots,pn-1\}$ such that $v_i=v$. For every
integer $k$ such that $0\leq k\leq\ell_v - \min\{\ell_w\dvtx w\in\tau^\circ\} $, put
\[
\phi_{(i)}(k)=\min\bigl\{j\in\{i,i+1,\ldots,i+pn-1\}\dvtx \ell_{v_j}=\ell_{v}-k\bigr\},
\]
and $\omega_{(i)}(k)=v_{\phi_{(i)}(k)}$. Then, if we also set $\omega_{(i)}(d_{\mathrm{gr}}(v,\partial))=\partial$,
it easily follows from (\ref{formuladgr}) that $(\omega_{(i)}(k),0\leq
k \leq d_{\mathrm{gr}}(v,\partial))$ is a discrete
geodesic from $v$ to $\partial$ in $M$. Such a geodesic is called a (discrete)
\textit{simple geodesic}.

The bound (\ref{bounddgr}) then simply expresses the fact that the distance
between $v_i$ and $v_j$ can be bounded by the length of the path
obtained by concatenating the simple geodesics $\phi_{(i)}$ and $\phi_{(j)}$ up
to their coalescence time.

%s2.3 #&#
\subsection{The CRT}
\label{CRT}

An important role in this work is played by the random real tree called
the CRT,
which was first introduced and studied by Aldous~\cite{Al1,Al3}. For
our purposes,
the CRT is conveniently viewed as the tree coded by a normalized
Brownian excursion. Throughout this work, the notation $\be=(\be_s)_{0\leq s\leq1}$
stands for a normalized Brownian excursion (see~\cite{RY}, Chapter XII,
for basic facts
about Brownian excursion theory). Recall from Section~\ref{sec1} the definition
of the pseudometric $d_\be$ and of the associated equivalence relation
$\sim_\be$. By definition, the CRT is
the quotient space ${\cal T}_\be:= [0,1]/ \sim_\be$ and is equipped
with the
induced distance, which is still denoted by $d_\be$. It is easy to
verify that
the topology of ${\cal T}_\be$ coincides with the quotient topology.

Then $({\cal T}_\be,d_\be)$ is a random compact real tree (see
Section 2.1 of
\cite{AM} for the
definition and basic properties of compact real trees).
We write $p_\be\dvtx [0,1] \la{\cal T}_\be$ for the canonical projection.
By convention,
${\cal T}_\be$ is rooted at the point $\rho:=p_\be(0)=p_\be(1)$. The
ancestral line of
a point $a$ of the CRT is the range of the unique (up to re-parametrization)
continuous and injective path from the root to $a$. This ancestral line
is denoted by $[\![ \rho, a ]\!]$. If $a,b\in{\cal T}_\be
$, we say that
$a$ is an ancestor of $b$ (or $b$ is a descendant of $a$) if
$a\in[\![ \rho,b ]\!]$. For every
$a\in{\cal T}_\be$, we can thus define the subtree of descendants of $a$.
If $a,b\in{\cal T}_\be$, we write $a\wedge b$ for the unique vertex such
that $[\![ \rho, a ]\!]\cap[\![ \rho,b ]\!]=
[\![ \rho,a\wedge b ]\!]$.

We refer to Section 2.2 in~\cite{AM} for more information about the
coding of compact real trees by continuous functions. Many properties
related to the genealogy of ${\cal T}_\be$ can be expressed conveniently
in terms of the coding function $\be$. For instance, if $s\in[0,1]$ is
given, a
point of the form $p_\be(t)$, $t\in[0,1]$, belongs to the ancestral line
of $p_\be(s)$ if and only if
\[
\be_t=\min_{s\wedge t\leq r\leq s\vee t} \be_r.
\]
We will use such simple facts without further comment in what follows.

A leaf of ${\cal T}_\be$ is a vertex $a$ such that ${\cal T}_\be
\setminus\{a\}$
is connected. If $t\in(0,1)$, the vertex $p_\be(t)$ is a \textit{leaf} if
and only if the equivalence class of $t$ for $\sim_\be$
is a singleton. The vertex $\rho=p_\be(0)=p_\be(1)$
is also a leaf. The set of all vertices of ${\cal T}_\be$
that are not leaves is called the skeleton of ${\cal T}_\be$
and denoted by $\operatorname{Sk}({\cal T}_\be)$.

%s2.4 #&#
\subsection{Brownian labels on the CRT}
\label{CRTlabels}

Brownian labels on the CRT are another crucial ingredient of our study.
We consider a real-valued process $Z=(Z_s)_{0\leq s\leq1}$
such that, conditionally given $(\be_s)_{0\leq s\leq1}$, $Z$ is a centered
Gaussian process with covariance
\[
E[Z_sZ_t \mid\be] = \min_{s\wedge t\leq r\leq s\vee t}
\be_r.
\]

Note, in particular, that $Z_0=0$ and $E[(Z_s-Z_t)^2\mid\be]=d_\be
(s,t)$. One way of constructing the process $Z$
is via the theory of the Brownian snake~\cite{Zu}. It is easy to
verify that
$Z$ has a continuous modification, which is even H\"older continuous
with exponent $\frac{1}{4}-\ve$ for every $\ve\in(0,\frac{1}{4})$. From
now on,
we always deal with this modification. From the invariance of
the law of the Brownian excursion under time-reversal, one
immediately gets that the processes $(\be_s,Z_s)_{0\leq s\leq1}$ and
$(\be_{1-s},Z_{1-s})_{0\leq s\leq1}$ have the same distribution.

From the formula $E[(Z_s-Z_t)^2\mid\be]=d_\be(s,t)$, one obtains that
\[
Z_s=Z_t\qquad\mbox{for every }s,t\in[0,1]\mbox{ such
that }d_\be (s,t)=0,\mbox{ a.s.}
\]
Hence, we may view $Z$ as indexed by the CRT ${\cal T}_\be$, in such a
way that
$Z_s=Z_{p_\be(s)}$ for every $s\in[0,1]$.
In what follows, we write indifferently
$Z_s=Z_a$ if $s\in[0,1]$ and $a\in{\cal T}_\be$ are such that
$a=p_\be(s)$.
Using standard techniques as in the proof of the classical Kolmogorov
lemma, one checks that the mapping
${\cal T}_\be\ni a\la Z_a$ is a.s. H\"older continuous
with exponent $\frac{1}{2}-\ve$
with respect to $d_\be$, for every $\ve\in(0,\frac{1}{2})$.

It is natural (and more intuitive than the presentation we just gave) to
interpret $Z$ as a Brownian motion indexed by the CRT. Although
the latter interpretation could be justified precisely, the approach we
took is mathematically
more tractable, as it avoids constructing a random process indexed by
a random set. As we will see below, the pair $({\cal T}_\be,(Z_a)_{a\in{\cal T}_\be})$
is a continuous analog of a uniformly distributed labeled $p$-tree
with $n$ black vertices.

Throughout this work, we will use the notation
\[
\Delta=-\min_{0\leq s\leq1} Z_s.
\]
Detailed information about the distribution of $\Delta$ can be found in
\cite{Delmas}. Here we will only use the simple fact that the
topological support of the
law of $\Delta$ is the whole of $\R_+$. This can be verified by
elementary arguments.
It is known (see~\cite{LGW}, Proposition 2.5) that there is an almost surely
unique instant $s_*\in(0,1)$ such that $Z_{s_*}=-\Delta$. We will
write $a_*=p_\be(s_*)$. Note that $a_*$ is a leaf of ${\cal T}_\be$.

We say that $t\in(0,1]$ is a left-increase time of $\be$, respectively of $Z$,
if there exists $\ve\in(0,t)$ such that $\be_s\geq\be_t$, respectively
$Z_s\geq Z_t$,
for every $s\in[t-\ve,t]$. We similarly define the
notion of a right-increase time. Note that the equivalence class
of $t$ for $\sim_\be$ is a singleton if and only if $t$ is neither a
left-increase time nor a right-increase time of $\be$. The following
result is Lemma 3.2 in~\cite{LGP}.

%le2.1 #&#
\begin{lemma}
\label{increase}
With probability one, any point $t\in[0,1]$ which is a right-increase or
a left-increase time of $\be$ is neither
a right-increase nor
a left-increase time of $Z$.
\end{lemma}

We set for every $r\geq0$,
\[
S_r=\inf\bigl\{s\in[0,1]\dvtx Z_s=-r\bigr\}
\]
with the usual convention $\inf\varnothing=\infty$. Note that
$S_r<\infty$ if and only if $r\leq\Delta$. If $r\in(0,\Delta]$, then
by definition $S_r$ is a left-increase time of $Z$, and Lemma~\ref{increase} implies that the equivalence class of
$S_r$ for $\sim_\be$ is a singleton, so that $p_\be(S_r)$
is a leaf of ${\cal T}_\be$ (the latter property is also true for $r=0$).

The following lemma shows that, in some sense,
labels do not vary too much between $S_{r}$ and $S_{r+\ve}$
when $\ve$ is small.

%le2.2 #&#
\begin{lemma}
\label{auxil}
There exists a constant $\beta_0\in(0,1)$ such that the following holds.
Let $\mu,A,\kappa$ be three reals with $0<\mu<A$ and $\kappa\in(0,1)$.
There exists a constant $C_{A,\mu,\kappa}$ such that, for every
$r\in[\mu,A]$ and $\ve\in(0,\mu/2)$,
\[
P \Bigl[\{ S_{r}\leq1-\kappa\}\cap \Bigl\{\sup_{s\in[S_{r-\ve},S_r]}
Z_s\geq-r+\sqrt{\ve} \Bigr\} \Bigr] \leq C_{A,\mu,\kappa}
\ve^{\beta_0}.
\]
\end{lemma}

Our proof of Lemma~\ref{auxil} depends on certain fine properties of
the Brownian snake,
which are also used in the proof of another more difficult lemma (Lemma~\ref{snakelemma} below). For this reason, we postpone the proof of both
results to
the \hyperref[app]{Appendix}.

For every $s,t\in[0,1]$ such that $s\leq t$, we set
\[
D^\circ(s,t)=D^\circ(t,s) = Z_s+
Z_t -2\max \Bigl(\min_{r\in[s,t]} Z_r,
\min_{r\in[t,1]\cup[0,s]} Z_r \Bigr).
\]
We then set, for every
$a,b\in{\cal T}_\be$,
\[
D^\circ(a,b)=\min\bigl\{ D^\circ(s,t)\dvtx s,t\in[0,1],
p_\be(s)=a, p_\be (t)=b\bigr\}.
\]
This is equivalent to the definition given in the introduction.
Suppose that $D^\circ(a,b)=0$ for some $a,b\in{\cal T}_\be$ with
$a\not=b$. Then
we can find $s,t\in[0,1]$ such that $p_\be(s)=a$, $p_\be(t)=b$ and
$D^\circ(s,t)=0$.
Clearly, $s$ and $t$ must be (right or left) increase times of $Z$ and
Lemma~\ref{increase} implies that
both $a$ and $b$ are leaves of ${\cal T}_\be$.

As a function on ${\cal T}_\be\times{\cal T}_\be$, $D^\circ$ does
not satisfy the
triangle inequality, but we can set, for every $a,b\in{\cal T}_\be$,
\[
D^*(a,b)=\inf \Biggl\{ \sum_{i=1}^k
D^\circ(a_{i-1},a_i) \Biggr\},
\]
where the infimum is over all choices of the integer $k\geq1$ and
of $a_0,\ldots,a_k\in{\cal T}_\be$ such that $a_0=a$ and $a_k=b$.
Then $D^*$ is a pseudometric on ${\cal T}_\be$, and obviously $D^*\leq
D^\circ
$. It will sometimes be convenient
to view $D^*$ as a function on $[0,1]^2$, by setting
\[
D^*(s,t)=D^*\bigl(p_\be(s),p_\be(t)\bigr)
\]
for every $s,t\in[0,1]$.

As a consequence of
Theorem 3.4 in~\cite{IM}, the property $D^*(a,b)=0$ holds if and only
if $D^\circ(a,b)=0$,
for every $a,b\in{\cal T}_\be$, a.s. (to be precise, the results of
\cite{IM}
are formulated in terms of
a pair $(\ov\be, \ov Z)$ which corresponds to re-rooting the CRT at
the vertex $p_\be(s_*)$ with a minimal
label---see Section 2.4 in~\cite{IM}---however, the preceding
formulation easily follows from
the results stated in~\cite{IM}).

%s2.5 #&#
\subsection{Convergence toward the Brownian map}
\label{convBrma}

For every integer $n\geq1$, let $M_n$ be a random rooted and pointed
$2p$-angulation, which is uniformly distributed
over the set $\m_n^{p}$. We can write $M_n$ as the image under the BDG
bijection
of a pair $(\theta_n,\ve_n)$, where
$\theta_n=(\tau_n,(\ell^n_v)_{v\in\tau^\circ_n})$ is a random
labeled $p$-tree
and $\ve_n$ is a random variable with values in $\{0,1\}$. Clearly,
$\theta_n$ is uniformly distributed
over the set $\T^p_n$ (and $\ve_n$ is uniformly distributed over $\{
0,1\}$).
We write $v^n_0,v^n_1,\ldots,v^n_{pn}$ for the contour sequence of
$\tau^\circ_n$.
We denote the contour function of
$\tau_n^\circ$ by $C^n=(C^n_i)_{0\leq i\leq pn}$ and
the label function of $\theta_n$ by
$\Lambda^n=(\Lambda^n_i)_{0\leq i\leq pn}$.
We extend the definition of both $C^n$ and $\Lambda^n$
to the real interval $[0,pn]$ by linear interpolation.\vadjust{\goodbreak}

Let $\mm_n$ stand for the vertex set of $M_n$. Thanks to
the BDG bijection, we have the
identification
\[
\mm_n=\tau^\circ_n\cup\{\partial\},
\]
where $\partial$ denotes the distinguished vertex of $M_n$. We also
observe that the notation $\mm_n$ is consistent
with Section~\ref{sec1}, since the random rooted $2p$-angulation
$\mathbf{M}_n$ obtained from $M_n$ by ``forgetting'' the distinguished
vertex of $M_n$ is uniformly distributed over $\a^{2p}_n$. Therefore,
when proving Theorem~\ref{mainresult}, we may assume that
the random metric space $(\mm_n,d_{\mathrm{gr}})$ is constructed from $M_n$
as explained above.

If $i,j\in\{0,1,\ldots,pn\}$,\vspace*{1pt}
we set $d_n(i,j)=d_{\mathrm{gr}}(v^n_i,v^n_j)$. We have then
$|\Lambda^n_i-\Lambda^n_j|\leq d_n(i,j)$ by
(\ref{formuladgr}) and the triangle inequality. As in~\cite{IM}, Section 3,
we extend the definition of $d_n(s,t)$ to noninteger values
of $(s,t)\in[0,pn]^2$ by setting
\begin{eqnarray*}
{d}_n(s,t)&=&\bigl(s-\lfloor s\rfloor\bigr) \bigl(t-\lfloor t\rfloor
\bigr){d}_n\bigl(\lceil s\rceil,\lceil t\rceil\bigr) + \bigl(s-
\lfloor s\rfloor\bigr) \bigl(\lceil t\rceil-t\bigr) {d}_n\bigl(\lceil
s\rceil,\lfloor t\rfloor\bigr)
\\
&&{}+ \bigl(\lceil s\rceil-s\bigr) \bigl(t-\lfloor t\rfloor\bigr){d}_n
\bigl(\lfloor s\rfloor,\lceil t\rceil\bigr) + \bigl(\lceil s\rceil-s\bigr) \bigl(
\lceil t\rceil-t\bigr){d}_n\bigl(\lfloor s\rfloor,\lfloor t\rfloor
\bigr),
\end{eqnarray*}
where $\lfloor t\rfloor= \max\{k\in\Z\dvtx k\leq t\}$ and $\lceil
t\rceil
=\min\{k\in\Z\dvtx k>t\}$.

The following theorem shows that the contour
and label processes and the distance process associated
with $M_n$ have a joint scaling limit, at least along a suitable
sequence of integers converging to $\infty$. This result is closely
related to~\cite{IM}, Theorem 3.4.
To simplify notation, we set
\[
\lambda_p=\frac{1}{2}\sqrt{\frac{p}{p-1}},\qquad
\kappa_p= \biggl(\frac{9}{4p(p-1)} \biggr)^{1/4}.
\]

%th2.3 #&#
\begin{theorem}
\label{mainIM}
From every sequence of integers converging to $\infty$, we can extract
a subsequence $(n_k)_{k\geq1}$
along which the following convergence in distribution of continuous
processes holds:
%
%e5 #&#
\begin{eqnarray}
\label{basicconv} && \bigl(\lambda_p n^{-1/2}
C^n_{pnt}, \kappa_p n^{-1/4}
\Lambda^n_{pnt }, \kappa_p n^{-1/4}
d_n( pns, pnt) \bigr)_{0\leq s\leq1, 0\leq
t\leq1}
\nonumber
\\[-8pt]
\\[-8pt]
\nonumber
&&\qquad \build{\la}_{n\to\infty}^{(d)} \bigl(\be_t,Z_t,D(s,t) \bigr)_{0\leq s\leq1,0\leq t\leq
1},
\end{eqnarray}
where
the pair $(\be,Z)$ is as in Section~\ref{CRTlabels}, and
$(D(s,t))_{0\leq s\leq1,0\leq t\leq1}$ is a continuous random process
such that the function $(s,t)\la D(s,t)$ defines a
pseudometric on $[0,1]^2$, and the following properties hold:
\begin{longlist}[(a)]
\item[(a)] $D(s,s_*)=Z_s+\Delta=D^\circ(s,s_*)$ for every
$s\in[0,1] $;
\item[(b)] $D(s,t)\leq D^*(s,t)\leq D^\circ(s,t)$ for every $s,t\in[0,1]$.
\end{longlist}
For every $s,t\in[0,1]$, we put $s\approx t$ if $D(s,t)=0$. Then, a.s.
for every
$s,t\in[0,1]$,
the property $s\approx t$ holds if and only if $D^*(s,t)=0$ or,
equivalently, $D^\circ(p_\be(s),p_\be(t))=0$.

Finally, set $\mm_\infty= [0,1]/ \approx$ and equip $\mm_\infty$
with the distance induced by~$D$, which is still denoted by $D$. Then,
along the same sequence where the convergence~(\ref{basicconv}) holds,
the random compact metric spaces
\[
\bigl(\mm_n,\kappa_p n^{-1/4}d_{\mathrm{gr}}
\bigr)
\]
converge in distribution to $(\mm_\infty,D)$ in the sense of the
Gromov--Hausdorff convergence.
\end{theorem}

\begin{remas*}(a) The bound $D(s,t)\leq D^\circ(s,t)$ is
an analog of the bound~(\ref{bounddgr}). Since $D$ satisfies the triangle
inequality, this bound immediately gives\break $D(s, t)\leq D^*(s,t)$
[and $D^*(s, t)\leq D^\circ(s,t)$ is true by definition as we already
noticed].\vspace*{-6pt}
\begin{longlist}[(b)]
\item[(b)] The convergence of the first two components in (\ref{basicconv})
does not require the use of a subsequence; see~\cite{MaMi}.

\item[(c)]
The identity $D(s,s_*)=Z_s+\Delta$ is a continuous analog of
formula (\ref{formuladgr}).

\item[(d)]
It is not hard to prove that equivalence classes for $\approx$ can contain
at most~$3$ points (see the discussion in~\cite{IM}, Section 3).
Moreover, if $s$ and $t$ are distinct points of $[0,1)$ such that
$s\approx t$, then we have either $p_\be(s)=p_\be(t)$ or $D^\circ
(s,t)=0$, but
these two properties cannot hold simultaneously by Lemma~\ref{increase}.
\end{longlist}
\end{remas*}

\begin{pf*}{Proof of Theorem~\ref{mainIM}} Although this theorem is very close to the results of \cite
{IM}, it cannot be
deduced immediately from that paper, because~\cite{IM} deals with rooted
$2p$-angulations, where the associated tree is constructed by using distances
from the root vertex, whereas in our setting of rooted and pointed
$2p$-angulations
the associated tree is obtained by considering the distances from the
distinguished
vertex. Still, the arguments in Section~3 of~\cite{IM} can be adapted to
the present setting. The convergence of the first two components in
(\ref{basicconv})
is deduced from~\cite{MaMi}, Theorem 8 (we should note that
\cite{MaMi} deals with the so-called height process, which is a variant
of the contour process, and
the corresponding variant of the label process,
but it is easy to verify that limit theorems for the height process can be
translated in terms of the contour process; see, for example, Section 1.6 in
\cite{trees}). From this convergence, the tightness of
the laws of the processes $ (n^{-1/4} d_n( pns, pnt)  )_{0\leq s\leq1, 0\leq
t\leq1}$ is derived exactly as in~\cite{IM}, Proposition 3.2, or in
\cite{Buzios}, Section 6, in the particular
case $p=2$. It follows that the
convergence (\ref{basicconv}) holds along a suitable subsequence and, via
the Skorokhod representation theorem, we may even assume that this
convergence holds a.s. The other assertions of the theorem are then obtained
in a straightforward way (see
Section 3 of~\cite{IM} or Section 6 of~\cite{Buzios}), with the
exception of the fact that
$D(s,t)=0$ implies $D^*(s,t)=0$.
To verify the latter fact, one can reproduce the rather delicate
arguments of~\cite{IM}, Section 4,
in the present setting. Alternatively, one can use the estimates for
the volume of\vadjust{\goodbreak}
balls proved in~\cite{AM}, Section 6, and follow the ideas that will be
developed below in Section
6 to get a sharper comparison estimate between $D$ and $D^*$.
We leave the details to the reader.
\end{pf*}

We will write $\bp$ for the canonical projection from
$[0,1]$ onto $\mm_\infty=[0,1]/\approx$. As a consequence
of the bound $D\leq D^\circ$, this projection is continuous when
$[0,1]$ is equipped with the usual Euclidean distance. The volume measure
$\operatorname{Vol}$ on $\mm_\infty$ is the image of the Lebesgue measure
on $[0,1]$ under the projection~$\bp$.

From the characterization of the equivalence
relation $\approx$, we see that $\mm_\infty$
can be viewed as well as a quotient space of ${\cal T}_\be$, for the
equivalence
relation $\simeq$ defined by $a\simeq b$ if and only if $D^\circ(a,b)=0$
(this is consistent with the presentation we gave in Section~\ref{sec1}).
We then write
$\Pi$ for the canonical projection from ${\cal T}_\be$ onto $\mm_\infty$
in such a way that $\bp=\Pi\circ p_\be$. Noting that
the topology on ${\cal T}_\be$ is the quotient topology and that
$\bp$ is continuous, it follows that $\Pi$ is also continuous.
We set $x_*=\bp(s_*)= \Pi(a_*)$. Note that property~(a)
in the theorem identifies all distances from $x_*$ in $\mm_\infty$
in terms of the label process $Z$.

We can define $D^*(x,y)$ for every
$x,y\in\mm_\infty$, so that
$D^*(\Pi(a),\Pi(b))=D^*(a,b)$ for every
$a,b\in{\cal T}_\be$. Then $D^*$ is also a random distance on $\mm_\infty$.
Most of what follows is devoted to
proving that $D(x,y)=D^*(x,y)$ for every $x,y\in\mm_\infty$.
If this equality holds, the limiting space in Theorem~\ref{mainIM}
coincides with $(\mm_\infty, D^*)$ and in particular does not
depend on the choice of the sequence $(n_k)_{k\geq1}$.
The statement of Theorem~\ref{mainresult} (in the bipartite case
when $q=2p$ is even) follows.

Notice that we already know by
property (b) of the theorem that $D\leq D^*$ and that
an easy compactness argument shows that the topologies
induced, respectively, by $D$ and by $D^*$ on $\mm_\infty$ coincide, as it
was already noted in~\cite{IM}. Furthermore, it
is immediate from properties (a) and (b) in the theorem that
%
%e6 #&#
\begin{equation}
\label{distance-root} D^*(x_*,x)= D(x_*,x)\qquad\mbox{for every }x\in
\mm_\infty.
\end{equation}

%s2.6 #&#
\subsection{Geodesics in the Brownian map}
\label{GeoBm}

If $x,y$ are points in a metric space $(E,d)$, a (continuous) geodesic
from $x$ to $y$ is
a path $(\omega(t),0\leq t\leq d(x,y))$ such that
$\omega(0)=x$, $\omega(d(x,y))=y$ and
$d(\omega(t),\omega(t'))=t'-t$ for every $0\leq t\leq t'\leq d(x,y)$.
The metric space $(E,d)$ is called geodesic if for any two
points $x,y\in E$ there is (at least) one geodesic
from $x$ to $y$.

From general results about Gromov--Hausdorff limits of
geodesic spaces~\cite{BBI}, Theorem 7.5.1, we get that
$(\mm_\infty, D)$ is almost surely a geodesic space.
Detailed information about the geodesics in $\mm_\infty$
has been obtained in~\cite{AM}, and we summarize
the results that will be needed below.

Let $s\in[0,1]$. For every $r\in[0, D(s,s_*)]$, we set
\[
\varphi_s(r)=\cases{ %
\inf\bigl\{t
\in[s,1]\dvtx Z_t=Z_s-r\bigr\}, &\quad $\mbox{if }\min\bigl
\{Z_t\dvtx t\in[s,1]\bigr\} \leq Z_s-r,$
\vspace*{2pt}\cr
\inf\bigl\{t\in[0,s]\dvtx Z_t=Z_s-r\bigr\},&\quad
$\mbox{otherwise.}$ }
\]
Since $D(s,s_*)= Z_s+\Delta$, the preceding definition makes sense.
For every $r\in[0, D(s,s_*)]$, set
\[
\Gamma_s(r)= \bp\bigl(\varphi_s(r)\bigr).
\]
By construction, $D^\circ(\varphi_s(r),\varphi_s(r'))=r'-r$ for every
$0\leq r\leq r'\leq D(s,s_*)$. On the other hand, by property (a)
of Theorem~\ref{mainIM}, we have also
\begin{eqnarray*}
r'-r&=&D^\circ\bigl(\varphi_s(r),
\varphi_s\bigl(r'\bigr)\bigr) \geq D\bigl(
\varphi_s(r),\varphi_s\bigl(r'\bigr)
\bigr)
\\[2pt]
&\geq& D\bigl(s_*,\varphi_s\bigl(r'\bigr)\bigr)-D
\bigl(s_*,\varphi_s(r)\bigr) =r'-r.
\end{eqnarray*}
It follows that $\Gamma_s$ is a geodesic in $(\mm_\infty, D)$.
Using property (b) of Theorem~\ref{mainIM}, we have then
$D^*(\Gamma_s(r),\Gamma_s(r'))=r'-r$ for every
$0\leq r\leq r'\leq D(s,s_*)$, and thus $\Gamma_s$ is also a geodesic
in $(\mm_\infty,D^*)$.

The geodesics of the form $\Gamma_s$
are called \textit{simple geodesics}. They are indeed the continuous
analogs of the discrete simple geodesics discussed at
the end of Section~\ref{BDGbij}.

The following theorem reformulates the main results of~\cite{AM} in
our setting.

%th2.4 #&#
\begin{theorem}
\label{geodesicroot}
All geodesics in $(\mm_\infty,D)$ from an arbitrary vertex of $\mm_\infty$
to $x_*$ are simple geodesics, and therefore also geodesics in $(\mm_\infty,D^*)$.
\end{theorem}

\begin{pf} For the same reason that was discussed in the proof of
Theorem~\ref{mainIM}, this result is not a mere restatement of
Theorem 7.4 and Theorem~7.6 in~\cite{AM}. However, it can be deduced
from these results along the following lines. Showing that all geodesics
from an arbitrary vertex of $\mm_\infty$
to $x_*$ are simple geodesics is easily seen to be equivalent to verifying
that a geodesic ending at $x_*$ cannot visit the skeleton
$\operatorname{Sk}({\cal T}_\be)$, except possibly at its starting point.
However, points of the skeleton are exactly those from which there are
(at least) two distinct simple geodesics. Hence, supposing that there
exists a geodesic ending at $x_*$ that visits the skeleton at a strictly
positive time, one could construct two geodesics $\omega$ and $\omega'$
starting from the same point and both ending at $x_*$ such that
$\omega(t)=\omega'(t)$ for every $t\in[0,\ve]$, for some $\ve>0$.
By the
invariance of the Brownian map under uniform re-rooting (Theorem 8.1 of
\cite{AM})
and the main results of~\cite{AM}, this does not occur.
\end{pf}

If $x\in\mm_\infty$
is such that $\bp^{-1}(x)$ is a singleton, Theorem~\ref{geodesicroot}
shows that there is a unique geodesic
from $x$ to $x_*$.
The particular case $x=\bp(0)$ plays an important role in the
remaining part of this work. In this case $\bp^{-1}(x)=\{0,1\}$, a.s.,
but it is trivial
that $\Gamma_0=\Gamma_1$, so that there is a.s. a unique
geodesic from $\bp(0)$ to $x_*$. To simplify notation, we will write
$\Gamma=\Gamma_0$ for this unique geodesic. We note that we have
$\varphi_0(r)=S_r$,
for every $r\in[0,\Delta]$, where $S_r$ was
introduced in Section~\ref{CRTlabels}, and, thus,
$\Gamma(r)=\bp(S_r)$.

%s3 #&#
\section{Maps with geodesic boundaries}\label{sec3}

%s3.1 #&#
\subsection{Discrete maps with geodesic boundaries}
\label{Dimageobo}

We will now describe a variant of the BDG bijection that produces a
$2p$-angulation with a boundary. We start from
a labeled $p$-tree $\theta=(\tau,(\ell_v)_{v\in\tau^\circ})$
with $n$ black vertices, and we set
\[
\delta= -\min\bigl\{\ell_v\dvtx v\in\tau^\circ\bigr\}+1.
\]
We use again the notation $v_0,v_1,\ldots,v_{pn}$ for the
contour sequence of $\tau^\circ$. We write~$M$
for the rooted and pointed $2p$-angulation associated with $\theta$ via
the BDG bijection (we should have fixed $\ve\in\{0,1\}$ to
determine the orientation of the root edge, but the choice of $\ve$
is irrelevant in what follows), and $d_{\mathrm{gr}}$ for the graph distance
on the vertex set $\mm$.

We then add $\delta-1$ vertices $\wt v_1,\wt v_2,\ldots,\wt v_{\delta
-1}$ to the tree $\tau$
in the following way. If $k=k_\varnothing(\tau)$ is the number of children
of $\varnothing$ in $\tau$, we put $\wt v_1=(k+1)$, $\tilde v_2=(k+1,1)$,
$\tilde v_3=(k+1,1,1)$ and so on until $\wt v_{\delta
-1}=(k+1,1,1,\ldots,1)$.
For notational convenience, we also set $\wt v_0=\varnothing$ and $\wt
v_\delta=\partial$.
Then $\wt\tau:=\tau\cup\{\wt v_1,\ldots,\wt v_{\delta-1}\}$ is
again a
plane tree (but no longer a $p$-tree). By convention, we put
\[
\wt\tau^\circ=\tau^\circ\cup\{\wt v_1,\ldots,\wt
v_{\delta-1}\}.
\]
We thus view $\wt v_1,\ldots,\wt v_{\delta-1}$ as white vertices with labels
$\ell_{\wt v_i}=-i$ for $i=1,\ldots,{\delta-1}$.

Now recall the construction of edges in the BDG bijection: For
every $i\in\{0,1,\ldots, pn-1\}$ with $\ell_i>-\delta+1$, the corner
$v_i$ is connected by an edge to the corner $v_j$, where
$j\in\{i,i+1,\ldots,i+pn-1\}$ is the successor of $i$. Note that every
corner of
$\tau$ corresponds to one corner of $\wt\tau$ (the vertex
$\varnothing$
has one more corner in $\wt\tau$,
except in
the particular case $\delta=1$).
To construct the planar map with a boundary, we follow rules similar to those
of the
BDG bijection. We start by drawing an edge between $v_i$ and $\partial$,
for all $i\in\{0,1,\ldots,pn-1\}$ such that $\ell_{v_i}=\min\ell$. Then,
let $i\in\{0,1,\ldots,pn-1\}$ such that $\ell_{v_i}>\min\ell$.
If the successor $j$ of $i$ is in $\{i+1,i+2,\ldots,pn\}$, we draw an
edge between $v_i$ and~$v_j$,
as we did before. However,
if the successor of
$i$ is in $\{pn+1,\ldots,i+pn-1\}$, we instead draw an edge between
$v_i$ and $\wt v_{-\ell_i+1}$ (since each new vertex $\wt v_j$ is
assigned the
label $-j$, $v_i$ is again connected by an edge to the next vertex
of $\wt\tau$ with a smaller label). Finally, for every $i\in\{
0,1,\ldots,\delta-1\}$,
we also draw an edge between $\wt v_i$ and $\wt v_{i+1}$ (in
particular, we draw an
edge between $\wt v_{\delta-1}$ and $\partial$).

%f3 #&#
\begin{figure}

\includegraphics{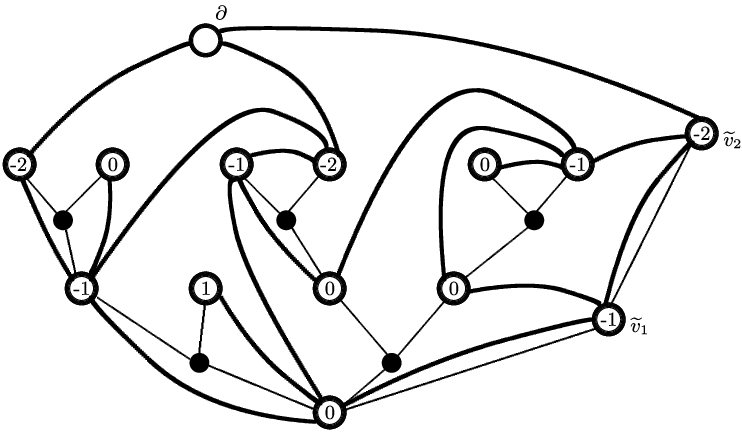}

\caption{The DMGB associated with the $6$-angulation of Figure \protect\ref{fig2}. In
this case, $\delta=3$ and the
two extra vertices $\wt v_1$ and $\wt v_2$ appear on the right of the
figure. The map is bounded by
the two boundary geodesics connecting the root of the tree to the
vertex $\partial$.}\label{fig3}
\end{figure}

The preceding construction gives a planar map $\wt M$ with vertex set
$\wt\mm=\wt\tau^\circ\cup\{\partial\}$ (see Figure~\ref{fig3} for an example).
The planar map $\wt M$ is in general not a $2p$-angulation. Leaving
aside the special case
$\delta=1$, where $\wt M=M$, the map $\wt M$ can be viewed
as a $2p$-angulation with a boundary. Indeed, it is not hard to verify
that every face
of $\wt M$ has degree $2p$ (and corresponds to one face in the planar
map $M$),
with the exception of one face, which has degree $2 \delta$ and is
bounded by the two
geodesics from $\varnothing$ to $\partial$ that are defined as follows:
$\gamma(0)=\wt\gamma(0)=\varnothing$, $\gamma(\delta)=\wt\gamma
(\delta
)=\partial$,
and
for every $i\in\{1,\ldots, \delta-1\}$,
\begin{eqnarray*}
\gamma(i)&=&v_{\phi(i)}\qquad \mbox{where }\phi(i)=\min\{j\geq0\dvtx \ell_{v_j}=-i\},
\\
\wt\gamma(i)&=&\wt v_i.
\end{eqnarray*}
Let $\wt d_{\mathrm{gr}}$ be the graph
distance on the vertex set $\wt\mm$.
The following properties are easily checked:
\begin{longlist}[(iii)]
\item[(i)] $\gamma$ and $\wt\gamma$ are two geodesics from
$\varnothing
$ to
$\partial$ in $\wt M$, that intersect only
at their initial and final points;
\item[(ii)]
$d_{\mathrm{gr}}(v,v')\leq\wt d_{\mathrm{gr}}(v,v')$ for every $v,v'\in\tau^\circ$;
\item[(iii)] $\wt d_{\mathrm{gr}}(v,\partial)=\ell_v+\delta$ for every $v\in
\wt
\tau^\circ$, and,
in particular, $\wt d_{\mathrm{gr}}(v,\partial)=d_{\mathrm{gr}}(v,\partial)$ for every
$v\in\tau^\circ$;
\item[(iv)] $\wt d_{\mathrm{gr}}(\varnothing,v)=d_{\mathrm{gr}}(\varnothing,v)$ for every
$v\in\tau^\circ$.
\end{longlist}
Informally, $M$ can be recovered from $\wt M$ by gluing the two
geodesics $\gamma$
and $\wt\gamma$ onto each other (and, in particular, identifying
$v_{\phi(i)}$ with $\wt v_i$ for every \mbox{$i=1,\ldots,\delta-1$}). This
explains why distances from $\varnothing$
or from $\partial$ are the same in $M$ and in $\wt M$, whereas other distances
may be different.
Note that the geodesic $\gamma$
coincides with the discrete simple geodesic\vspace*{1pt} $\omega_{(0)}$ introduced
at the end of Section~\ref{BDGbij}.

We will say that $\wt M$ is the discrete map with geodesic boundaries
(in short, the DMGB) associated with $M$. Notice that the boundary of
$\wt M$
is only piecewise geodesic since it consists of the union of two
geodesics from
$\varnothing$ to $\partial$. We sometimes say that $\gamma$,
respectively $\gamma'$, is the left boundary geodesic, respectively the right
boundary geodesic,
of $\wt M$.

The definition of discrete simple geodesics can be extended to $\wt M$ in
the following way. Recall the notation at the end of Section~\ref{BDGbij},
and let $i\in\{0,\ldots,pn-1\}$. If the minimal label on $\tau^\circ$
is attained at $v_j$ for some $j\in\{i,\ldots,pn\}$,\vspace*{1pt} we just put
$\wt\omega_{(i)}=\omega_{(i)}$, which is also a geodesic from
$v_i$ to $\partial$ in $\wt M$. On the other hand, if the
preceding property does not hold, there is a unique integer
$k\in\{1,\ldots,d_{\mathrm{gr}}(v_i,\partial)-1\}$ such that $\phi_{(i)}(k-1)\leq pn$
and $\phi_{(i)}(k) > pn$. Then the edge of $M$ between
$\omega_{(i)}(k-1)$ and $\omega_{(i)}(k)$ does not exist in $\wt M$,
but instead there is an edge of $\wt M$ between $\omega_{(i)}(k-1)$
and~$\wt v_{k'}$, where $k'=k-\ell_{v_{i}}$. So we can put
$\wt\omega_{(i)}(j)=\omega_{(i)}(j)$ if $j\leq k-1$ and
$\wt\omega_{(i)}(j)=\wt\gamma(j-\ell_{v_{i}})$ if $k\leq j\leq
d_{\mathrm{gr}}(v_i,\partial)$,
and $\wt\omega_{(i)}$ is again a geodesic from $v_i$ to $\partial$ in~$\wt M$.

%s3.2 #&#
\subsection{Scaling limits}
\label{scalim}

We now apply the construction of the preceding subsection
to a random $2p$-angulation $M_n$ that is uniformly distributed
over the set $\m_n^{p}$. We let
$\theta_n=(\tau_n,(\ell^n_v)_{v\in\tau_n})$ be the labeled $p$-tree
associated with $M_n$, and
we write $v^n_0,\ldots,v^n_{pn}$ for the contour sequence of $\tau^\circ_n$.
As previously, we also write
$(C^n_i)_{0\leq i\leq pn}$ for the contour function of $\tau^\circ_n$
and $(\Lambda^n_i)_{0\leq i\leq pn}$ for the label
function of $\theta_n$.

The DMGB associated with $M_n$ is denoted by $\wt M_n$. We also let
$\mm_n$ and $\wt\mm_n$ denote, respectively, the vertex set of $M_n$
and the vertex set of $\wt M_n$.

Recall the definition of the function $d_n$ before Theorem \ref
{mainIM}. For every $i,j\in\{0,1,\ldots,pn\}$, we also set
\[
\wt d_n(i,j)=\wt d_{\mathrm{gr}}\bigl(v^n_i,v^n_j
\bigr).
\]
A simple adaptation of the proof of (\ref{bounddgr})
gives the bound
\[
d_n(i,j)\leq\wt d_n(i,j)\leq d_n^\bullet(i,j),
\]
where, for every $i,j\in\{0,1,\ldots,pn\}$,
\[
d_n^\bullet(i,j)=\Lambda^n_i +
\Lambda^n_j - 2\min_{i\wedge j\leq
k\leq
i\vee j} \Lambda^n_k+2.
\]
Similarly as in the case of $d_n$, we extend the definition of
$\wt d_n$ to
$[0,pn]\times[0,pn]$ by linear interpolation.
The next proposition reinforces the joint convergence (\ref{basicconv})
in Theorem~\ref{mainIM}
by considering also the distance $\wt d_n$ jointly with the contour
and label processes and the distance $d_n$.

%pr3.1 #&#
\begin{proposition}
\label{convergenceDMGB}
From every sequence of integers converging to $\infty$, we can extract
a subsequence $(n_k)_{k\geq1}$
along which the following convergence in distribution of continuous processes
indexed by $s,t\in[0,1]$ holds:
\begin{eqnarray}
\label{basicconv2}\qquad && \bigl(\lambda_pn^{-1/2}
C^n_{pnt}, \kappa_pn^{-1/4}
\Lambda^n_{pnt}, \kappa_pn^{-1/4}d_n(pns,pnt),
\kappa_pn^{-1/4} \wt d_n(pns, pnt) \bigr)
\nonumber
\\[-8pt]
\\[-8pt]
\nonumber
 && \qquad \build{\la}_{n\to\infty}^{(d)} \bigl(
\be_t,Z_t,D(s,t), \wt D(s,t) \bigr),
\nonumber
\end{eqnarray}
where
$(\be_t)_{0\leq t\leq1}$,
$(Z_t)_{0\leq t\leq1}$ and $(D(s,t))_{0\leq s\leq1,0\leq t\leq1}$ are
as in Theorem~\ref{mainIM}, and
$(\wt D(s,t))_{0\leq s\leq1,0\leq t\leq1}$ is a continuous random process
such that $D\leq\wt D$ and the function $(s,t)\la\wt D(s,t)$ defines a
pseudometric on $[0,1]^2$. We put $s\equiv t$ if and only if $\wt
D(s,t)=0$. The property
$s\equiv t$ holds if and only if at least one of the following
two conditions holds:
\begin{longlist}[(a)]
\item[(a)] $s\sim_{\be} t$;
\item[(b)] $Z_s=Z_t= \min_{r\in[s\wedge t,s\vee t]} Z_r$.
\end{longlist}
Finally, along the same sequence where the convergence (\ref{basicconv2})
holds, we have the joint convergence in distribution of random metric
spaces in the Gromov--Hausdorff sense:
\[
\bigl(\bigl(\mm_n,\kappa_p n^{-1/4}d_{\mathrm{gr}}
\bigr), \bigl(\wt\mm_n,\kappa_p n^{-1/4}\wt
d_{\mathrm{gr}}\bigr) \bigr) \build{\la}_{n\to\infty}^{(d)} \bigl((
\mm_\infty, D),(\wt\mm_\infty, \wt D) \bigr),
\]
where $(\mm_\infty, D)$ is as in Theorem~\ref{mainIM},
$\wt\mm_\infty=[0,1]/\equiv$, and $\wt D$ is the induced distance
on $\wt\mm_\infty$.
\end{proposition}

\begin{pf} From the bound $\wt d_n(i,j)\leq d_n^\bullet(i,j)$, we can
use the same arguments
as in the proof of~\cite{IM}, Proposition 3.2, to verify that the
sequence of
laws of the processes $(n^{-1/4}\wt d_n(pns,pnt))_{0\leq s,
t\leq1}$ is tight in the space of all probability measures
on $C([0,1]^2,\R)$. To be specific, we write for every $s,t,s',t'\in[0,1]$,
\begin{eqnarray*}
&&\bigl|n^{-1/4}\wt d_n(pns,pnt)- n^{-1/4}\wt
d_n\bigl(pns',pnt'\bigr)\bigr|
\\
&&\qquad \leq n^{-1/4}\bigl(\wt d_n\bigl(pns,pns'
\bigr)+ \wt d_n\bigl(pnt,pnt'\bigr)\bigr)
\\
&&\qquad \leq n^{-1/4}\bigl(d^\bullet_n
\bigl(pns,pns'\bigr)+ d^\bullet_n
\bigl(pnt,pnt'\bigr)\bigr).
\end{eqnarray*}
By (\ref{basicconv}), the processes $(\kappa_p n^{-1/4}d^\bullet_n(pns,pnt))_{0\leq
s,t\leq1}$
converge in distribution to the process
\[
\Bigl( Z_s + Z_t - 2\min_{s\wedge t\leq r\leq s\vee t}
Z_r \Bigr)_{0\leq
s,t\leq1}.
\]
It then follows that, for every fixed $\delta>0$, the quantity
\[
P \Bigl(\sup_{|s-s'|\leq\varepsilon, |t-t'|\leq\varepsilon} \bigl|n^{-1/4}\wt d_n(pns,pnt)-
n^{-1/4}\wt d_n\bigl(pns',pnt'
\bigr)\bigr| > \delta \Bigr)
\]
can be made arbitrarily small, uniformly in $n$, by choosing
$\varepsilon>0$ small enough. This yields the
desired tightness property.

Using also the convergence (\ref{basicconv}), we see that we
can extract a subsequence along which the convergence (\ref
{basicconv2}) holds,
and obviously the processes $\be, Z$ and $D$ satisfy the same
properties as in Theorem~\ref{mainIM}.
From now on we restrict our attention to values of $n$
in this subsequence. Using the Skorokhod representation theorem,
we may assume throughout the proof that the convergence
(\ref{basicconv2}) holds a.s.

From the analogous properties for $\wt d_n$, it is immediate that
$\wt D$ is symmetric and satisfies the triangle inequality. Note that the
bound $d_n\leq\wt d_n$ implies that $D\leq\wt D$.

Let us now verify that $\wt D(s,t)=0$ if and only if
(at least) one of the two conditions (a) and (b) holds. First, if (a)
holds, the
same argument as in the proof of Proposition 3.3(iii) in~\cite{IM}
shows that $\wt D(s,t)=0$. Then, by passing to the limit $n\to\infty$
in the bound
\[
n^{-1/4} \wt d_n\bigl(\lfloor pns\rfloor,\lfloor pnt\rfloor
\bigr) \leq n^{-1/4} d_n^\bullet\bigl(\lfloor pns
\rfloor,\lfloor pnt\rfloor\bigr),
\]
we easily get that, a.s. for every $s,t\in[0,1]$,
\[
\wt D(s,t)\leq Z_s + Z_t - 2 \min_{s\wedge t\leq r\leq s\vee t}
Z_r.
\]
If (b) holds, the right-hand side vanishes, which immediately gives
\mbox{$\wt D(s,t)=0$}.

Conversely, suppose that $\wt D(s,t)=0$, and without loss of generality
assume that $s<t$. Since $D\leq\wt D$, we have also
$D(s,t)=0$ and, by Theorem~\ref{mainIM}, we know that either (a) holds
(in which case we are done) or
\[
Z_s=Z_t= \max \Bigl(\min_{r\in[s,t]}
Z_r, \min_{r\in[t,1]\cup[0,s]} Z_r \Bigr).
\]
If
\[
Z_s=Z_t=\min_{r\in[s,t]} Z_r,
\]
then (b) holds. So we concentrate on the case where
%
%e7 #&#
\begin{equation}
\label{geobotech1} Z_s=Z_t=\min_{r\in[t,1]\cup[0,s]}
Z_r.
\end{equation}
Assuming that this equality holds and that $s\sim_\be t$
does not hold, we will arrive at a contradiction, which will complete the
proof of our characterization of the equivalence relation $\equiv$.
We may assume that $s>0$ and $t<1$ [the case $s=0$, $t=1$ is excluded, and
then we note that $\min_{[0,\ve]}Z<0$
and $\min_{[1-\ve,1]}Z<0$, for every $\ve\in(0,1)$, a.s. by Lemma
\ref
{increase}]. Then we can find
positive integers $i_n,j_n\in\{0,1,\ldots,pn\}$,
with $i_n\leq j_n$, such that $s=\lim(pn)^{-1}i_n$
and $t=\lim(pn)^{-1}j_n$, and we have
%
%e8 #&#
\begin{equation}
\label{geobotech2} \lim_{n\to\infty} \kappa_p n^{-1/4}
\wt d_{\mathrm{gr}}\bigl(v^n_{i_n},v^n_{j_n}
\bigr) = \wt D(s,t)=0.
\end{equation}
From (\ref{geobotech1}) and the fact that the minimum
of $Z$ is attained at a unique time, we know that, for $n$ large,
the minimum of $\ell^n$ will be attained (only) in $\{i_n,\ldots,j_n\}$.
Let $k_n\in\{i_n,\ldots,j_n\}$ be the largest integer
such that $\ell^n_{v^n_{k_n}} = \min\ell^n$, and write
$[\![ \varnothing, v^n_{k_n} ]\!]$
for the ancestral line of $v^n_{k_n}$ in $\tau_n$.
By
construction, if an edge of $\wt\mm_n$ connects a point of
$\{v^n_{k_n},v^n_{k_n+1},\ldots, v^n_{pn})$ to a point of
$\{v^n_0,v^n_1,\ldots, v^n_{k_n}\}$, then (at least) one of these
two points must belong to $[\![ \varnothing, v^n_{k_n} ]\!]$.
Therefore, if $\omega_n$ is a geodesic path from $v^n_{i_n}$
to $v^n_{j_n}$ in $\wt\mm_n$, it must either visit $\partial$ or intersect
$[\![ \varnothing, v^n_{k_n} ]\!]$ at (at least) one
point, which
may be written in the form $v^n_{\ell_n}$ with
$\ell_n\in\{0,1,\ldots,pn\}$. The case when $\omega_n$ visits
$\partial$
does not occur when $n$ is large, since this would imply that
$Z_s=Z_t=\min Z$,
which is absurd. In the other case, we can find a subsequence
of the sequence
$(pn)^{-1}\ell_n$ that converges to $r\in[0,1]$, and
automatically $p_\be(r)$ belongs to the ancestral line of the
vertex $a_*=p_\be(s_*)$ minimizing $Z$. Furthermore, it is also clear
from (\ref{geobotech2}) that $\wt D(s,r)=\wt D(t,r)=0$ and, therefore,
$D(s,r)=D(t,r)=0$. By Theorem~\ref{mainIM}, we must have
$D^\circ(p_\be(s),p_\be(r))=0$. However,
$p_\be(s)$ is a leaf of ${\cal T}_\be$ (by Lemma~\ref{increase}),
whereas $p_\be(r)$ is a point of $\operatorname{Sk}({\cal T}_\be)$. This
contradicts our
previous observation that, if $a,b\in{\cal T}_\be$ with $a\not=b$,
$D^\circ(a,b)=0$ may hold only if
$a$ and $b$ are both leaves of~${\cal T}_\be$. This contradiction
completes the proof of the characterization of the property
$\wt D(s,t)=0$.

We still have to prove the last convergence of the proposition. The
almost sure
convergence of the random compact metric spaces $(\mm_n,\kappa_p
n^{-1/4}d_{\mathrm{gr}})$
toward $(\mm_\infty,D)$ is easily derived from the (almost sure) convergence
(\ref{basicconv2}) as in the first part of the proof of Theorem 3.4 in
\cite{IM}.
A similar argument will give the almost sure
convergence of $(\wt\mm_n,\kappa_p n^{-1/4}\wt d_{\mathrm{gr}})$
toward $(\wt\mm_\infty,\wt D)$. Let us provide details for the
sake of completeness. We first observe that we may discard the
extra vertices that we added to $\mm_n$ in order to define $\wt\mm_n$.
Indeed, it is immediate that the Hausdorff distance between
$\mm_n$ [viewed as a compact subset of the metric space
$(\wt\mm_n,\wt d_{\mathrm{gr}})$] and $\wt\mm_n$ is bounded by $1$,
and so the Gromov--Hausdorff convergence of $(\wt\mm_n,\kappa_p
n^{-1/4}\wt d_{\mathrm{gr}})$
will follow from that of $(\mm_n,\kappa_p n^{-1/4}\wt d_{\mathrm{gr}})$. For
the same reason,
we may replace $\mm_n$ by $\mm_n\setminus\{\partial\}$. We then
construct a
correspondence between $\mm_n\setminus\{\partial\}$ and $\wt\mm_\infty
$ by saying that, for
every $i\in\{0,1,\ldots,pn\}$ and $s\in[0,1]$, the vertex $v^n_i$ is in
correspondence
with the equivalence class of $s$ in $\wt\mm_\infty=[0,1]/\equiv$ if
$|i-pns|\leq1$. Thanks to the convergence (\ref{basicconv2}), we can
easily verify that
the distortion of this convergence, when $\mm_n\setminus\{\partial\}$
is equipped
with the distance $\kappa_p n^{-1/4}\wt d_{\mathrm{gr}}$ and $\wt\mm_\infty$
with $\wt D$, tends to $0$ a.s. as $n\to\infty$. This completes the
proof.
\end{pf}

Let us state some important properties of the space $(\wt\mm_\infty,\wt D)$.
In the following proposition, as well as in the remaining part of this
section, we consider the processes $(\be,Z,D,\wt D)$ and the associated
random metric spaces $(\mm_\infty,D)$ and $(\wt\mm_\infty,\wt D)$ that
arise from the
convergences of the preceding proposition via the choice of
a suitable subsequence.
We write $\wt\bp$ for the canonical projection from
$[0,1]$ onto $\wt\mm_\infty= [0,1]/\equiv$. Recall the notation
\[
\Delta=- \min\bigl\{Z_s\dvtx s\in[0,1]\bigr\}.
\]

%pr3.2 #&#
\begin{proposition}
\label{propertiesGeoBo}
\textup{(i)} For every $s,t\in[0,1]$,
\[
\wt D(s,t)\leq Z_s + Z_t - 2\min_{s\wedge t\leq r\leq s\vee t}
Z_r.\vspace*{-6pt}
\]
\begin{longlist}[(iii)]
\item[(ii)] For every $s\in[0,1]$, $\wt D(0,s)= D(0,s)$.
\item[(iii)] For every $s\in[0,1]$, $\wt
D(s,s_*)=D(s,s_*)=Z_s+\Delta$.
\item[(iv)]For every $t\in[0,\Delta]$, put
\begin{eqnarray*}
\wt\Gamma(t)&=&\wt\bp\bigl(\inf\bigl\{s\in[0,1]\dvtx Z_s=-t\bigr\}
\bigr),
\\
\wt\Gamma'(t)&=&\wt\bp\bigl(\sup\bigl\{s\in[0,1]\dvtx Z_s=-t\bigr\}\bigr).
\end{eqnarray*}
Then $\wt\Gamma$ and $\wt\Gamma'$ are two geodesic paths from $\wt
\bp
(0)$ to $\wt\bp(s_*)$
in $(\wt\mm_\infty,\wt D)$, which intersect only at their initial and
final points.
\end{longlist}
\end{proposition}

\begin{pf} Property (i) was already derived in the preceding proof.
Properties~(ii)
and (iii) follow from the analogous properties of a DMGB stated at the
end of Section~\ref{Dimageobo} by a straightforward passage to
the limit. Let us verify (iv). First it is immediate that $\wt\Gamma
(0)=\wt\Gamma'(0)=\wt\bp(0)=\wt\bp(1)$,
and $\wt\Gamma(\Delta)=\wt\Gamma'(\Delta)=\wt\bp(s_*)$. Then,
from (ii)
or (iii), we have
\[
\wt D\bigl(\wt\bp(0),\wt\bp(s_*)\bigr)=\wt D(0,s_*)=D(0,s_*)=\Delta.
\]
On the other hand, for every $0\leq t\leq t'\leq\Delta$, (i) gives
\[
\wt D\bigl(\wt\Gamma(t),\wt\Gamma\bigl(t'\bigr)\bigr)\leq
t'-t.
\]
Thanks to the triangle inequality, this implies that
$\wt D(\wt\Gamma(t),\wt\Gamma(t'))= t'-t$ for every $0\leq t\leq
t'\leq
\Delta$. The fact
that $\wt\Gamma'$ is a geodesic path is proved in a similar way. Finally,
the property $\wt\Gamma(t)\not=\wt\Gamma'(t)$ for $t\in(0,\Delta)$
follows from
the characterization of the equivalence relation $\equiv$ in Proposition
\ref{convergenceDMGB}, using also Lemma~\ref{increase}.
\end{pf}

We will now explain how the space $(\wt\mm_\infty,\wt D)$ can be constructed
from $(\mm_\infty,D)$ by ``cutting'' the surface $(\mm_\infty,D)$ along
the geodesic
$\Gamma$, which produces the two geodesics $\wt\Gamma$ and $\wt
\Gamma'$.
Such surgery is common in the study of the geometry of surfaces, but
since we are
working in a singular setting we will proceed with some care.

We set
\[
{\cal R}_\Gamma= \bigl\{\Gamma(t)\dvtx 0<t<\Delta\bigr\} \subset
\mm_\infty,\vadjust{\goodbreak}
\]
and write $\ov{\cal R}_\Gamma={\cal R}_\Gamma\cup\{\bp(0),\bp
(s_*)\}$ for the
closure of ${\cal R}_\Gamma$.
We consider a set $\mm^\bullet_\infty$ which is obtained from
$\mm_\infty$ by duplicating every point of ${\cal R}_\Gamma$. Formally,
\[
\mm^\bullet_\infty = (\mm_\infty\setminus{\cal
R}_\Gamma) \cup\bigl\{(x,0)\dvtx x\in{\cal R}_\Gamma\bigr\} \cup
\bigl\{(x,1)\dvtx x\in{\cal R}_\Gamma\bigr\}.
\]

We then define a topology on $\mm^\bullet_\infty$ by the following
prescriptions:
\begin{itemize}
\item If $x\in\mm_\infty\setminus\ov{\cal R}_\Gamma
$, a subset
of $\mm^\bullet_\infty$ is a neighborhood of $x$ in $\mm^\bullet_\infty$
if and only if it contains a neighborhood of $x$ in $\mm_\infty$.
\item A subset $V$
of $\mm^\bullet_\infty$ is a neighborhood of $\bp(0)$, respectively of
$\bp(s_*)$,
in $\mm^\bullet_\infty$ if and only if there exists a neighborhood
$U$ of $\bp(0)$, respectively of $\bp(s_*)$, in $\mm_\infty$, and $\ve>0$
such that
\[
V\supset \bigl( (U\setminus{\cal R}_\gamma)\cup\bigl\{\bigl(\Gamma
(t),1\bigr)\dvtx 0\leq t\leq\ve \bigr\} \cup\bigl\{\bigl(\Gamma(t),0\bigr)\dvtx 0\leq t\leq
\ve\bigr\} \bigr),
\]
respectively,
\[
V\supset \bigl( (U\setminus{\cal R}_\gamma)\cup\bigl\{\bigl(
\Gamma(t),1\bigr)\dvtx \Delta-\ve \leq t\leq\Delta\bigr\} \cup\bigl\{\bigl(
\Gamma(t),0\bigr)\dvtx \Delta-\ve\leq t\leq\Delta\bigr\} \bigr).
\]
\item If $x\in{\cal R}_\Gamma$, a subset $V$
of $\mm^\bullet_\infty$ is a neighborhood of $(x,0)$, respectively of $(x,1)$,
in $\mm^\bullet_\infty$
if and only if there exists a neighborhood $U$ of $x$ in $\mm_\infty$
such that
$V$ contains $U\cap\bp([0,s_*])$, respectively $U\cap\bp([s_*,1])$.
\end{itemize}

We write $\pi$ for the obvious projection from $\mm^\bullet_\infty$
onto $\mm_\infty$, and note that $\pi$ is continuous. We define a metric
$D^\bullet$ on $\mm^\bullet_\infty$ by setting, for every $x,y\in
\mm^\bullet_\infty$,
\[
D^\bullet(x,y) = \inf\bigl\{ L(\pi\circ g)\dvtx g\in C\bigl(
\mm^\bullet_\infty,x\to y\bigr)\bigr\},
\]
where $ C(\mm_\infty^\bullet,x\to y)$ stands for the set of all
continuous paths
$g\dvtx [0,1]\to\mm^\bullet_\infty$ such that $g(0)=x$ and $g(1)=y$, and
$L(\pi\circ g)$
denotes the length of the path $\pi\circ g$ in $(\mm_\infty,D)$.
Informally, the
paths of the form $\pi\circ g$ are those paths from $\pi(x)$ to $\pi(y)$
in $\mm_\infty$ that do not cross the geodesic $\Gamma$.

%pr3.3 #&#
\begin{proposition}
\label{identBMGB}
The metric spaces $(\wt\mm_\infty,\wt D)$ and $(\mm_\infty^\bullet,
D^\bullet)$
are almost surely isometric.
\end{proposition}

\begin{pf} This proposition is not needed in the derivation of
our main result, and so we only sketch the proof. We first observe that
there is
an obvious bijection~$h$ from $\mm_\infty^\bullet$ onto $\wt\mm_\infty$
such that,
for every $t\in(0,\Delta)$,
\begin{eqnarray*}
h\bigl(\bigl(\Gamma(t),0\bigr)\bigr)&=& \wt\Gamma(t),
\\
h\bigl(\bigl(\Gamma(t),1\bigr)\bigr)&=& \wt\Gamma'(t).
\end{eqnarray*}
Indeed, every $x\in\mm_\infty\setminus{\mathcal R}_\Gamma$ clearly
corresponds to
exactly one point $y$ of $\wt\mm_\infty$ and we take $h(x)=y$.

We then need to verify that $h$ is an isometry. Since $(\wt\mm_\infty,\wt D)$
is a geodesic space (as a Gromov--Hausdorff limit of rescaled graphs),
we know
that, for every $z_1,z_2\in\mm_\infty^\bullet$,
\[
\wt D\bigl(h(z_1),h(z_2)\bigr)=\inf\bigl\{ \wt L(\wt f)\dvtx f\in C\bigl(\wt\mm_\infty, h(z_1)\to h(z_2)
\bigr)\bigr\},
\]
where $C(\wt\mm_\infty, h(z_1)\to h(z_2))$ is the set of all continuous
paths $\wt f\dvtx [0,1] \la\wt\mm_\infty$
such that $\wt f(0)=h(z_1)$ and $\wt f(1)=h(z_2)$, and $ \wt L(\wt f)$
denotes the
length of $\wt f$ in~$ \wt\mm_\infty$. It is easy to verify that
$\wt
f\in C(\wt\mm_\infty, h(z_1)\to h(z_2))$ if and only if it can be
written in the form $\wt f = h\circ g$,
where $g\in C(\mm_\infty^\bullet,z_1\to z_2)$. Moreover, we have then
%
%e9 #&#
\begin{equation}
\label{identBMGBtech} \wt L(\wt f)= L(\pi\circ g).
\end{equation}
Once (\ref{identBMGBtech}) has been established, it readily follows
from the preceding
formulas for $\wt D$ and $D^\bullet$ that we have $\wt D(h(z_1),h(z_2))
= D^\bullet(z_1,z_2)$ for every $z_1,z_2\in\mm_\infty^\bullet$, so
that $h$
is an isometry. We leave the details of the proof of (\ref
{identBMGBtech}) to the reader.
\end{pf}

%s3.3 #&#
\subsection{A technical lemma}
\label{techlemm}

We will now use the results of the preceding subsection to derive a
technical lemma that will play an important role later in this work.
Recall the notation introduced at the beginning of Section~\ref{scalim}.
In particular, the random $2p$-angulation $M_n$ is uniformly distributed
over the set $\m_n^{p}$, and the DMGB associated with $M_n$ is denoted
by $\wt M_n$.
We put $\Delta_n= d_{\mathrm{gr}} (\varnothing,\partial)$ (where $d_{\mathrm{gr}}$
refers to the graph distance
in $M_n$) and
following the end of Section~\ref{Dimageobo}, we introduce the two
distinguished
geodesics from $\varnothing$ to $\partial$ in $\wt M_n$, which are
denoted by
$\gamma_n$ and $\gamma'_n$.

%le3.4 #&#
\begin{lemma}
\label{traverse}
We can find two positive constants $\ve$ and $\eta$ such that, for every
sufficiently large integer $n$,
\[
P \Bigl[10 n^{1/4}\leq\Delta_n \leq11n^{1/4},
\mathop{\mathop{\min_{n^{1/4}\leq i\leq9n^{1/4}
}}}_{n^{1/4}
\leq j\leq9n^{1/4}} \wt d_{\mathrm{gr}}
\bigl(\gamma_n(i),\gamma'_n(j)\bigr)\geq\ve
n^{1/4} \Bigr] \geq\eta.
\]
\end{lemma}

\begin{rema*}Lemma~\ref{traverse} is related to the fact that the
(continuous) geodesics
$\wt\Gamma$ and $\wt\Gamma'$ in Proposition~\ref{propertiesGeoBo} do
not intersect
except at their initial and final points. In the discrete setting,
``interior points'' of the geodesics
$\gamma_n$ and $\gamma'_n$ stay at a distance of order $n^{1/4}$.
\end{rema*}

\begin{pf} Set $\ve_k=2^{-k}$ and $\eta_k=2^{-k}$ for every
integer $k\geq0$. We argue by contradiction and assume that
for every $k\geq0$ we can find an integer $n_k\geq k$
such that
%
%e10 #&#
\begin{equation}
\label{traverstech1}P \Bigl[10n_k^{1/4}\leq
\Delta_{n_k} \leq11n_k^{1/4}, \mathop{\mathop{
\min_{n_k^{1/4}\leq i\leq9n_k^{1/4} }}_{n_k^{1/4}
\leq j\leq9n_k^{1/4}}} \wt d_{\mathrm{gr}}\bigl(
\gamma_{n_k}(i),\gamma'_{n_k}(j)\bigr)\geq
\ve_k n_k^{1/4} \Bigr] < \eta_k.\hspace*{-30pt}
\end{equation}
From Proposition~\ref{convergenceDMGB} and replacing the sequence
$(n_k)_{k\geq0}$ by a
subsequence, we may assume that the convergence (\ref{basicconv2})
holds along the
sequence $(n_k)_{k\geq0}$. Notice that the bound (\ref{traverstech1})
remains valid after
this replacement. By using the Skorokhod\vadjust{\goodbreak} representation theorem, we may
even assume
that (\ref{basicconv2}) holds almost surely. From now on
until the end of the proof, we consider only values of $n$ belonging to
the sequence $(n_k)_{k\geq0}$,
even if this is not indicated in the notation.

We then consider the random closed subsets $K$ and $K'$ of $[0,1]$,
which are
defined by
\begin{eqnarray*}
K&=&\Bigl\{t\in[0,1]\dvtx Z_t=\min_{0\leq s\leq t} Z_s
\in[\kappa_p,9\kappa_p]\Bigr\},
\\[-2pt]
K'&=&\Bigl\{t\in[0,1]\dvtx Z_t=\min_{t\leq s\leq1}
Z_s \in[\kappa_p,9\kappa_p]\Bigr\}.
\end{eqnarray*}
We recall that by definition, for every $n$, and $0\leq i\leq\Delta_n-1$,
\[
\gamma_n(i)=v^n_{\phi_n(i)}\qquad\mbox{where }
\phi_n(i)=\min\bigl\{j\geq0\dvtx \ell^n_{v^n_j}=-i
\bigr\}.
\]
No similar formula holds for $\gamma'_n$, but we can write, for every
$1\leq i\leq\Delta_n$,
\[
\wt d_{\mathrm{gr}}\bigl(\gamma'_n(i),v^n_{\phi'_n(i)}
\bigr)=1\qquad\mbox{where }\phi'_n(i)=\max\bigl\{j\leq pn\dvtx
\ell^n_{v^n_j}=-i+1\bigr\}.
\]
The latter equality easily follows from the construction of edges in
the DMGB. We thus have, for every
$i\in\{0,\ldots,\Delta_n-1\}$ and $j\in\{1,\ldots,\Delta_n\}$,
%
%e11 #&#
\begin{equation}
\label{traverstech2} \wt d_{\mathrm{gr}}\bigl(\gamma_n(i),
\gamma'_n(j)\bigr)\leq\wt d_{\mathrm{gr}}
\bigl(v^n_{\phi
_n(i)},v^n_{\phi'_n(j)}\bigr)+1 =
\wt d_n\bigl(\phi_n(i),\phi'_n(j)
\bigr) +1.
\end{equation}
Recall that $\ell^n_{v^n_j}=\Lambda^n_j$ and that the sequence of processes
$(\kappa_pn^{-1/4}\Lambda^n_{pnt})_{0\leq t\leq1}$ converges almost surely
to $(Z_t)_{0\leq t\leq1}$ by (\ref{basicconv2}). Elementary arguments
using the latter
convergence and the definition of the functions $\phi_n$ and $\phi'_n$
then show that, on the event $\{\Delta>10\kappa_p\}$,
%
%e12 #&#
\begin{equation}
\label{traverstech3} \sup_{n^{1/4}\leq i\leq9n^{1/4}} d\biggl(\frac{\phi_n(i)}{pn},K\biggr)
\build {\la }_{n\to\infty}^{\mathrm{a.s.}} 0,\qquad \sup_{n^{1/4}\leq j\leq9n^{1/4}} d\biggl(
\frac{\phi'_n(j)}{pn},K'\biggr) \build {\la }_{n\to\infty}^{\mathrm{a.s.}}
0,\hspace*{-30pt}
\end{equation}
where $d$ refers to the usual Euclidean distance on $[0,1]$. Notice
that, on the event $\{\Delta> 10\kappa_p\}$, we have
$\Delta_n \geq10n^{1/4}$ for $n$ large enough, a.s., and, in
particular, $\phi_n(i)$ and $\phi'_n(j)$ are well defined
for $n^{1/4}\leq i\leq9n^{1/4}$ and $n^{1/4}\leq j\leq9n^{1/4}$.

From (\ref{traverstech2}), (\ref{traverstech3}) and the convergence
(\ref{basicconv2}), we now get,
on the event $\{10\kappa_p<\Delta<11\kappa_p\}$,
\[
\liminf_{k\to\infty} \Bigl(\kappa_pn_k^{-1/4}
\mathop{\mathop{\min_{n_k^{1/4}\leq i\leq
9n_k^{1/4}}}}_{n_k^{1/4}\leq
j\leq9n_k^{1/4}} \wt d_{\mathrm{gr}}
\bigl(\gamma_{n_k}(i),\gamma'_{n_k}(j)\bigr)
\Bigr) \geq\inf_{t\in K,t'\in K'} \wt D\bigl(t,t'\bigr),
\]
almost surely. In particular, for every $\ve>0$,
%
%e13 #&#
\begin{eqnarray}
\label{traverstech5}\qquad &&\liminf_{k\to\infty}P \Bigl[10n_k^{1/4}
\leq\Delta_{n_k} \leq 11n_k^{1/4}, \mathop{
\mathop{\min_{n_k^{1/4}\leq j\leq9n_k^{1/4}}}_
{n_k^{1/4}\leq i\leq9n_k^{1/4} }} \wt d_{\mathrm{gr}}\bigl(
\gamma_{n_k}(i),\gamma'_{n_k}(j)\bigr)\geq\ve
n_k^{1/4} \Bigr]
\nonumber
\\[-8pt]
\\[-8pt]
\nonumber
&&\qquad \geq P \Bigl[10\kappa_p<\Delta<11\kappa_p,
\inf_{t\in K,t'\in K'} \wt D\bigl(t,t'\bigr) > \kappa_p^{-1}
\ve \Bigr].
\end{eqnarray}
The characterization of the equivalence relation $\equiv$ in
Proposition~\ref{convergenceDMGB}
shows that $\wt D(t,t')>0$ for every $t\in K$ and $t'\in K'$, a.s. on
the event $\{\Delta>10\kappa_p\}$.
By compactness, we have thus
\[
\inf_{t\in K,t'\in K'} \wt D\bigl(t,t'\bigr)>0
\]
a.s. on that event. In particular, we can fix $\ve>0$ so that the
right-hand side of (\ref{traverstech5}) is (strictly) positive.
This gives a contradiction with (\ref{traverstech1}), and this
contradiction completes the proof.
\end{pf}

In the next section we will use a minor extension of Lemma \ref
{techlemm}, concerning the
case when our random $p$-tree has a random number of black vertices.
Let $\mu>0$ and, for every
(sufficiently large) integer $n$, consider a random labeled $p$-tree
$\wh\theta_n$
whose size belongs to $[\mu^4n,2\mu^4n]$, and such that the
conditional distribution
of $\wh\theta_n$ given its size is uniform. With $\wh\theta_n$ we
associate a DMGB as
explained in Section~\ref{Dimageobo}, and we let $\wh\gamma_n$ and
$\wh\gamma'_n$
be, respectively, the left and right boundary geodesics in this map. We
also denote by
$\wh\Delta_n$ the common length of these geodesics. We can apply
Lemma~\ref{techlemm}
to $\wh\theta_n$ after conditioning on its size, and we get, for all
sufficiently large $n$,
%
%e14 #&#
\begin{eqnarray}
\label{traversext}&& P \Bigl[\bigl\{10\mu n^{1/4}\leq\wh
\Delta_n \leq15\mu n^{1/4}\bigr\}
\nonumber
\\[-8pt]
\\[-8pt]
\nonumber
&&\qquad{}\cap \Bigl\{ \mathop{
\mathop{\min_{2^{1/4}\mu n^{1/4}\leq i\leq9\mu n^{1/4}
}}_{2^{1/4}\mu n^{1/4}\leq j\leq9\mu
n^{1/4}}} \wt d_{\mathrm{gr}}\bigl(\wh
\gamma_n(i),\wh\gamma'_n(j)\bigr)\geq\ve \mu
n^{1/4} \Bigr\} \Bigr] \geq\eta.
\end{eqnarray}

%s4 #&#
\section{The traversal lemmas}
\label{keylemma}

We use the notation
of Sections~\ref{scalim} and~\ref{techlemm}.

%le4.1 #&#
\begin{lemma}
\label{keylem}
We can find a constant $\alpha_0>0$ such
that the following holds. For every $\alpha\in(0,\alpha_0)$, for
every choice
of the constants $\beta_1$ and $\beta_2$ such that $15\alpha<\beta_1<\beta_2$,
and for
every sufficiently large integer $n$, the probability of the event
%
%{\fontsize{10.3}{12.3}{\selectfont
\begin{eqnarray*}
&&\bigl\{\beta_1n^{1/4}<\Delta_n<
\beta_2 n^{1/4}\bigr\}\\
&&\qquad{}\cap \biggl\{ \wt d_{\mathrm{gr}}
\bigl(\gamma_n\bigl(\bigl\lfloor\alpha n^{1/4}\bigr\rfloor
\bigr),\gamma'_n\bigl(\bigl\lfloor \alpha
n^{1/4}\bigr\rfloor\bigr)\bigr)\\
&&\hspace*{18pt}\qquad= \wt d_{\mathrm{gr}}\biggl(\gamma_n\biggl(\biggl\lfloor
\frac{\alpha}{3} n^{1/4}\biggr\rfloor \biggr),\gamma'_n
\biggl(\biggl\lfloor\frac{\alpha}{3} n^{1/4}\biggr\rfloor\biggr)\biggr)
+2\biggl(\bigl\lfloor\alpha n^{1/4}\bigr\rfloor- \biggl\lfloor
\frac{\alpha}{3} n^{1/4}\biggr\rfloor \biggr) \biggr\}
\end{eqnarray*}
%}}
%
is bounded below by a positive constant independent of $n$.
\end{lemma}

We note that, provided that $\Delta_n\geq\lfloor\alpha
n^{1/4}\rfloor
$, one has
\[
\wt d_{\mathrm{gr}}\biggl(\gamma_n\bigl(\bigl\lfloor\alpha
n^{1/4}\bigr\rfloor\bigr),\gamma_n\biggl(\biggl\lfloor
\frac
{\alpha}{3} n^{1/4}\biggr\rfloor\biggr)\biggr) =\bigl\lfloor\alpha
n^{1/4}\bigr\rfloor- \biggl\lfloor\frac{\alpha}{3} n^{1/4}
\biggr\rfloor
\]
and similarly if $\gamma_n$ is replaced by $\gamma'_n$.
Therefore, the event considered in the lemma holds if and only if
$\beta_1n^{1/4}<\Delta_n<\beta_2 n^{1/4}$
and there exists a geodesic from $\gamma_n(\lfloor\alpha
n^{1/4}\rfloor
)$ to $\gamma'_n(\lfloor\alpha n^{1/4}\rfloor)$
in $\wt M_n$, which visits both points $\gamma_n(\lfloor\frac{\alpha
}{3} n^{1/4}\rfloor)$
and $\gamma'_n(\lfloor\frac{\alpha}{3} n^{1/4}\rfloor)$ in this order.

To simplify notation, we will write
$u_n=v^n_{\lfloor pn/2\rfloor}$\vspace*{1pt} in the remaining part of this section.
An important role in the proof below
will be played by subtrees branching
from the right side of the ancestral line of $u_n$ (see the end of
Section~\ref{labtrees}).

\begin{pf*}{Proof of Lemma~\ref{keylem}} Recall the notation $C^n$ for
the contour function and $\Lambda^n$
for the label function of the labeled $p$-tree tree $\theta_n=(\tau_n,(\ell^n_v)_{v\in\tau^\circ_n})$. We know from (\ref{basicconv})
that the pair of processes $(\lambda_pn^{-1/2} C^n_{pnt},\kappa_pn^{-1/4}\Lambda^n_{pnt})_{0\leq t\leq1}$ converges in distribution toward
$(\be_t,Z_t)_{0\leq t\leq1}$ [this convergence does not require the
use of
a subsequence; see remark (a) after Theorem~\ref{mainIM}]. We will use
this convergence in distribution to get that, with a probability bounded
from below
when $n$ is large, the pair $(C^n,\Lambda^n)$ satisfies certain
properties, which can then be expressed
in terms of properties of the subtrees branching from the right side of
the ancestral line of $u_n$.

We let $\ve>0$ be the constant appearing in Lemma~\ref{traverse}, and
we put
\[
A=\biggl\lfloor\frac{2}{\ve}\biggr\rfloor+1.
\]
We determine $\alpha_0$ by the condition $p\alpha_0^4A=\frac{1}{8}$.
Then we fix
$\alpha\in(0,\alpha_0)$ and $\beta_1,\beta_2$ such that $15\alpha
<\beta_1<\beta_2$.

Let $F$ be the event where the following properties hold:
\begin{enumerate}[(a)]
\item[(a)] We have $\be_{1/2}>\lambda_p$ and $Z_{1/2}<-\kappa_p\alpha$.
Moreover, for any vertex $a$ of ${\cal T}_\be$ that is an ancestor of
$p_\be
(1/2)$ in ${\cal T}_\be$ and is such
that $\frac{\lambda_p}{2}\leq d_\be(\rho,a)\leq\lambda_p$, we have
$3\kappa_p\alpha< Z_a <4\kappa_p \alpha$.
\item[(b)] $ \Delta=-\min\{Z_s\dvtx 0\leq s\leq1\} \in(\kappa_p\beta_1,\kappa_p\beta_2)$.
\item[(c)] For every $s$ such that either $0\leq s\leq\sup\{r\leq
\frac
{1}{2}\dvtx \be_r\leq\frac{\lambda_p}{2}\}$
or $\inf\{r\geq\frac{1}{2}\dvtx \be_r\leq\frac{\lambda_p}{2}\}\leq
s\leq
1$, we have $Z_s>- \kappa_p\alpha/6$.
\item[(d)] There exist $A$ subintervals $[s_1,t_1],\ldots,[s_A,t_A]$ of
$[0,1]$, with $ \frac{1}{2}<s_1<t_1<\cdots<s_A<t_A<1$,
such that, for every $i\in\{1,\ldots,A\}$:
\begin{enumerate}[(d1)]
\item[(d1)] $\be_{s_i}=\be_{t_i}=\min\{\be_s\dvtx \frac{1}{2}\leq
s\leq t_i\}
$, and $\be_{s_i}\in(\frac{\lambda_p}{2},\lambda_p)$.
\item[(d2)] $\alpha^4<t_i-s_i<2\alpha^4$.
\item[(d3)] $-15\kappa_p \alpha< \min_{s\in[s_i,t_i]}Z_s -Z_{s_i}
<-10\kappa_p \alpha$.
\end{enumerate}
\end{enumerate}
Let us comment on condition (d). By (d1), the intervals
$[s_1,t_1],\ldots,[s_A,t_A]$ correspond to excursions
of the process $(\be_{(1/2)+s})_{0\leq s\leq1/2}$ above its minimum
process. In particular, for every
$i\in\{1,\ldots,A\}$, $p_\be(s_i)=p_\be(t_i)$ belongs to the ancestral
line of $p_\be(1/2)$. In terms of the tree
${\cal T}_\be$ coded by $\be$, each interval $[s_i,t_i] $ can be interpreted
as a subtree branching from
the ancestral line of $p_\be(1/2)$ at level $\be_{s_i}$. Condition (d2)
then gives bounds for the ``mass'' of
this subtree, and condition~(d3) provides bounds for the minimal
relative label on this subtree.

Simple arguments show that $P(F)>0$. The conditions that do not involve
the label process $Z$ are easily seen
to hold with positive probability [note that our choice of $\alpha$
such that $p\alpha^4A<\frac{1}{8}$
makes it possible to fulfill condition (d2)]. The fact that the other
conditions then also hold with positive probability
requires a little more work, but we leave the details to the reader.

For every integer $n\geq1$, we then let $F_n$ be the
event where the following properties hold:
\begin{enumerate}[(a$'$)]
\item[(a$'$)] We have $C^n_{\lfloor pn/2\rfloor}=\frac
{1}{2}|u_n|>\sqrt {n}$ and $\ell^n_{u_n}<-\alpha n^{1/4}$. Moreover, if $v$
is a vertex of $\tau_n^\circ$ that is an ancestor of $u_n$ and such
that $\frac{1}{2}\sqrt{n}\leq\frac{1}{2}|v|\leq\sqrt{n}$, we have
$3\alpha n^{1/4} < \ell^n_v <4\alpha n^{1/4}$.
\item[(b$'$)] $ \Delta_n=1-\min\{\ell^n_v\dvtx v\in\tau_n^\circ\} \in
(\beta_1n^{1/4},\beta_2n^{1/4})$.
\item[(c$'$)] For every vertex $v$ of $\tau_n^\circ$ that belongs to a
subtree branching from the left side
or from the right side of the ancestral line of $u_n$ at level
(strictly) less than $\frac{1}{2}\sqrt{n}$,
we have $\ell^n_v\geq- \frac{\alpha}{6}n^{1/4}$.
\item[(d$'$)] There exist at least $A$ subtrees $\tau_{n,1},\ldots,\tau_{n,A}$ branching from the
right side of the ancestral line of $u_n$, such that, for every $i\in\{
1,\ldots,A\}$: %
\begin{enumerate}[(d1$'$)]
\item[(d1$'$)] The branching level of $\tau_{n,i}$ belongs to $[\frac
{1}{2}\sqrt{n},\sqrt{n}]$.
\item[(d2$'$)] $\alpha^4n<|\tau_{n,i}|<2\alpha^4n$.
\item[(d3$'$)] The minimal difference between the label of a vertex of
$\tau_{n,i}$ and the label
of its root belongs to $[-15\alpha n^{1/4},-10\alpha n^{1/4}]$.
\end{enumerate}
\end{enumerate}
In condition (d2$'$), we recall that the size $|\tau_{n,i}|$ is the
number of black vertices of $\tau_{n,i}$. See Figure~\ref{fig4} for a rough
illustration of conditions (a$'$), (c$'$), (d$'$).

%f4 #&#
\begin{figure}

\includegraphics{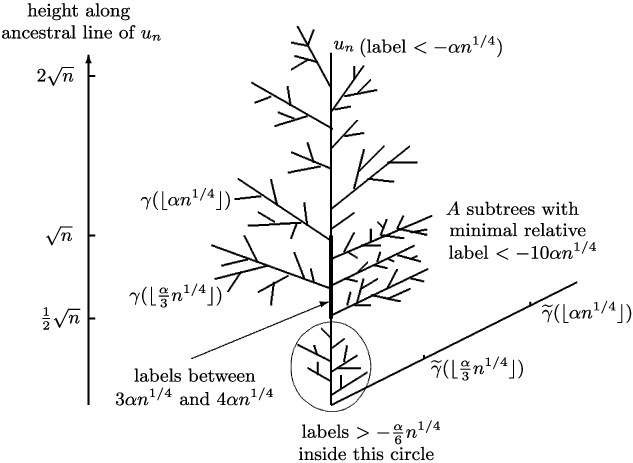}

\caption{Illustration of the proof. The geodesic from $\wt\gamma
(\lfloor
\alpha n^{1/4}\rfloor)$ to $\gamma(\lfloor\alpha n^{1/4}\rfloor)$
cannot visit the part of the ancestral line of $u_n$ between height
$\frac{1}{2}\sqrt{n}$ and height $\sqrt{n}$. If it does not visit
the part of the ancestral line between $0$ and $\frac{1}{2}\sqrt{n}$,
it has to cross the $A$ subtrees.}\label{fig4}
\end{figure}

The convergence in distribution of $(\lambda_pn^{-1/2}
C^n_{pnt},\kappa_pn^{-1/4}\Lambda^n_{pnt})_{0\leq t\leq1}$ toward
$(\be_t,Z_t)_{0\leq t\leq1}$ now implies that
\[
\liminf_{n\to\infty} P(F_n)\geq P(F).
\]
To see this, first note that we can replace the convergence in
distribution by an almost sure convergence, thanks to
the Skorokhod representation theorem. Then on the event $F$, the almost
sure (uniform) convergence
of $(\lambda_pn^{-1/2} C^n_{ pnt})_{0\leq t\leq1}$ toward $(\be_t)_{0\leq t\leq1})$ will imply the existence
of subintervals $[m_1,n_1],\ldots,[m_A,n_A]$ of $[pn/2,pn]$ such that
properties analogous to (d1),(d2)
hold for these subintervals and for the contour function $C^n$. From
the relation between the contour
function and the tree $\tau_n$, we then get, still on the event $F$ and
for large enough $n$, the existence
of subtrees satisfying the properties in (d$'$). The remaining part of
the argument is
straightforward.\looseness=-1

Fix $\nu>0$ such that $P(F)>\nu$. We can then find $n_0$ such that
$P(F_n)\geq\nu$ for
every $n\geq n_0$. Let us fix $n\geq n_0$ and argue under the
conditional probability $P(\cdot\mid F_n)$.
We can determine the choice of the subtrees $\tau_{n,1},\ldots,\tau_{n,A}$ by saying that we choose the first
$A$ subtrees branching from the right side of the ancestral line of
$u_n$ and satisfying the
conditions (d1$'$), (d2$'$), (d3$'$), and order them in lexicographical order.
As mentioned in Section~\ref{labtrees}, we can
view each $\tau_{n,i}$ as a (random) $p$-tree, via an obvious renaming
of the vertices, and we can
equip the vertices of this $p$-tree with labels obtained by taking the
difference of the original labels
(in $\theta_n$) with the label of the root of $\tau_{n,i}$. In this way
we obtain a random labeled
$p$-tree, which we denote by $\theta_{n,i}$, for every $i\in\{
1,\ldots,A\}$. Let $k_1,\ldots,k_A$
be integers with $\alpha^4n<k_i<2\alpha^4n$ for $i\in\{1,\ldots,A\}$.
We claim that
under the measure $P(\cdot\mid F_n)$ and conditionally on the event $\{
|\tau_{n,1}|=k_1,\ldots,
|\tau_{n,A}|=k_A\}$, the random labeled $p$-trees $\theta_{n,1},\ldots,\theta_{n,A}$
are independent, and the conditional distribution of each $\theta_{n,i}$ is uniform
over labeled $p$-trees with $k_i$ black vertices subject to the
constraint that the minimal
label belongs to $[-15\alpha n^{1/4},-10\alpha n^{1/4}]$. This follows
from the fact that
the tree $\theta_n$ is uniformly distributed, and the conditions
(a$'$), (b$'$), (c$'$) do not depend
on the properties of the trees $\theta_{n,i}$ [in the case of (b$'$), we
note that, because of (a$'$) and
the assumption $\alpha<\beta_1/15$, the minimal label in $\tau_n$
will certainly
not be attained at a vertex of one of the subtrees $\tau_{n,i}$].

We write $\wt M_{n,i}$ for the DMGB associated with $\theta_{n,i}$, and
we let $\gamma_{n,i}$ and $\gamma'_{n,i}$
be, respectively, the left and right boundary geodesic in $\wt M_{n,i}$.
Let $G_n$ be the intersection of $F_n$ with
the event where
%
%e15 #&#
\begin{equation}
\label{traverskeylem} \mathop{\mathop{\min_{2^{1/4}\alpha n^{1/4}\leq j\leq
9\alpha n^{1/4} }}_{2^{1/4}\alpha n^{1/4}
\leq k\leq9\alpha n^{1/4}}}
\wt d_{\mathrm{gr}}\bigl(\gamma_{n,i}(j),\gamma'_{n,i}(k)
\bigr)\geq\ve\alpha n^{1/4}
\end{equation}
for every $i\in\{1,\ldots,A\}$ [in (\ref{traverskeylem}), $\wt d_{\mathrm{gr}}$
obviously stands for the graph distance in $\wt M_{n,i}$]. From (\ref
{traversext}) and the preceding considerations,
we can find $n_1\geq n_0$ such that, for every $n\geq n_1$,
$P(G_n)\geq\eta^A P(F_n)\geq\eta^A \nu$. To complete the
proof of Lemma~\ref{keylem}, it now suffices to verify that the event
considered in this lemma contains $G_n$.

So suppose that $G_n$ holds. We already know that $ \Delta_n \in
(\beta_1n^{1/4},\beta_2n^{1/4})$ by (b$'$). Next
consider a geodesic path $\omega_n$ from $\gamma'_n(\lfloor\alpha
n^{1/4}\rfloor)$ to $\gamma_n(\lfloor\alpha n^{1/4}\rfloor)$
in $\wt M_n$. Recall that the label of both $\gamma_n(\lfloor\alpha
n^{1/4}\rfloor)$ and $\gamma'_n(\lfloor\alpha n^{1/4}\rfloor)$
is equal to $-\lfloor\alpha n^{1/4}\rfloor$. From the trivial bound
\[
\wt d_{\mathrm{gr}}\bigl(\gamma'_n\bigl(\bigl\lfloor
\alpha n^{1/4}\bigr\rfloor\bigr),\gamma'_n
\bigl(\bigl\lfloor \alpha n^{1/4}\bigr\rfloor\bigr)\bigr) \leq2 \bigl
\lfloor\alpha n^{1/4}\bigr\rfloor
\]
and the fact that labels correspond to distances from $\partial$ in
$\wt M_n$ (modulo a shift
by a fixed quantity), we immediately see that the path $\omega_n$
cannot visit a vertex whose label
is positive or strictly smaller than $-2\lfloor\alpha n^{1/4}\rfloor$.
To simplify notation, write $w_{n,i}$ for the
(white) vertex at generation $2i$ on the ancestral line of $u_n$, for
every $i\in\{0,1,\ldots,\frac{|u_n|}{2}\}$. It follows
from (a$'$) and the preceding considerations that $\omega_n$ does not
visit the set
$\{w_{n,i}\dvtx \frac{1}{2}\sqrt{n}\leq i\leq\sqrt{n}\}$.

%
%%
%
%
%
%
%
%(\lfloor
%
%
%(\lfloor
%
%
%(\lfloor
%
%(\lfloor
%
%$<-\alpha
%n^{1/4}$})}}
%
%{\alpha}{6}n^{1/4}$}}
%
%n^{1/4}$}}
%
%and $4\alpha n^{1/4}$}}
%
%
%
%

We claim that $\omega_n$ must visit the set $H_n:=\{w_{n,i}\dvtx 0\leq
i<\frac{1}{2}\sqrt{n}\}$. If the claim holds, the proof is
easily completed. Indeed, suppose that $\omega_n$ visits the vertex
\mbox{$w\in H_n$}. Then we can construct a geodesic
path $\wh\omega_n$, respectively $\wh\omega'_n$, that connects $w$ to
$\gamma_n(\lfloor\alpha n^{1/4}\rfloor)$, respectively
to $\gamma'_n(\lfloor\alpha n^{1/4}\rfloor)$, and visits the point
$\gamma_n(\lfloor\frac{\alpha}{3} n^{1/4}\rfloor)$,
respectively the point $\gamma'_n(\lfloor\frac{\alpha}{3} n^{1/4}\rfloor)$.
To construct~$\wh\omega_n$, pick any $k\in\{0,\ldots,\lfloor
pn/2\rfloor\}$ such that
$v^n_k=w$ and consider the simple geodesic $\wt\omega_{(k)}$ as
defined in Section~\ref{Dimageobo}. By condition (c$'$), this simple
geodesic will coalesce with
$\gamma_n$
at a point of the form $\gamma_n(j)$ with $j< \lfloor\frac{\alpha}{3}
n^{1/4}\rfloor$.
Therefore, we can just let $\wh\omega_n$ coincide with $\omega_{(k)}$
up to its hitting time of $\gamma_n(\lfloor\alpha n^{1/4}\rfloor)$.
Similarly, to construct $\wh\omega'_n$, we pick any $k'\in\{\lfloor
pn/2 \rfloor,\ldots,pn\}$ such that
\mbox{$v^n_{k'}=w$}. By condition (c$'$) again, the simple geodesic $\wt\omega_{(k')}$ will coalesce with
$\gamma'_n$ at a point of the form $\gamma'_n(j)$ with $j< \lfloor
\frac
{\alpha}{3} n^{1/4}\rfloor$, and we let
$\wh\omega'_n$ coincide with $\wt\omega_{(k')}$ up to its hitting time
of $\gamma'_n(\lfloor\alpha n^{1/4}\rfloor)$.
The concatenation of $\wh\omega_n$ and
$\wh\omega'_n$ gives a geodesic path from $\gamma'_n(\lfloor\alpha
n^{1/4}\rfloor)$ to $\gamma_n(\lfloor\alpha n^{1/4}\rfloor)$
that visits both $\gamma_n(\lfloor\frac{\alpha}{3} n^{1/4}\rfloor)$
and $\gamma'_n(\lfloor\frac{\alpha}{3} n^{1/4}\rfloor)$,
as desired.

It remains to verify the claim. We argue by contradiction and suppose
that $\omega_n$
does not visit $H_n$. Recall that $\omega_n$ does not visit the set $\{
w_{n,i}\dvtx \frac{1}{2}\sqrt{n}\leq i\leq\sqrt{n}\} $ either,
and notice that the condition $\ell^n_{u_n}<-\alpha n^{1/4}$ ensures
that $\gamma_n(\lfloor\alpha n^{1/4}\rfloor)$ belongs
to the left side of the ancestral line of $u_n$. Also recall that
labels along
$\omega_n$ must remain in the range\vadjust{\goodbreak} $[-2\lfloor\alpha n^{1/4}\rfloor,0]$. From these observations, the properties
(a$'$) and (d3$'$)
and the construction of edges in~$\wt M_n$, it follows that $\omega_n$
must visit each of the trees
$\tau_{n,A},\tau_{n,A-1},\ldots,\tau_{n,1}$ in this order before it
can hit the left side of the ancestral line of $u_n$.
Furthermore, the path $\omega_n$ hits $\tau_{n,A}$ for the first time
at a vertex $v$
such that the following property holds: There is no occurrence of the
label $\ell^n_v-1$ among vertices that appear
in the part of the contour sequence of $\tau^\circ_n$ corresponding to
$\tau_{n,A}$
after the last occurrence of $v$. Indeed, this property is needed for
$v$ to be connected to a vertex on the right side of $\tau_{n,A}$.
If we now view $v$ as a vertex of the DMGB $\wt M_{n,A}$, this means
that $v$ is connected by
an edge to the vertex $\gamma'_{n,A}(k)$, where $-k+1$ is the
difference between $\ell^n_v$
and the label of the root of $\tau_{n,A}$. Since $-2\lfloor\alpha
n^{1/4}\rfloor\leq\ell^n_v\leq0$, property
(a$'$) implies that $3\alpha n^{1/4}-p\leq k\leq6 \alpha n^{1/4}+p+1$
(notice that although the root of
$\tau_{n,A}$ may not belong to the ancestral line of $u_n$, its label
differs by at
most $p$ from the label of a vertex in this ancestral line, whose
generation belongs to
$[\frac{1}{2}\sqrt{n},\sqrt{n}]$). We may assume that $p+1\leq
\alpha
n^{1/4}$ and so we get that
$2\alpha n^{1/4}\leq k\leq7 \alpha n^{1/4}$. Similarly, the last
vertex of $\omega_n$ belonging to $\tau_{n,A}$
before $\omega_n$ first hits $\tau_{n,A-1}$ can be written as $\gamma_{n,A}(j)$
for some $j$ such that $2\alpha n^{1/4}\leq j \leq7 \alpha n^{1/4}$.
We can now use (\ref{traverskeylem}) to obtain that the
time spent by $\omega_n$ between its first visit of $\tau_{n,A}$ and
its first visit of $\tau_{n,A-1}$ is at least
$\ve\alpha n^{1/4} $. The same lower bound holds for the time between
the first hitting time of $\tau_{n,i}$
and the first hitting time of $\tau_{n,i-1}$, for every
$i=A,A-1,\ldots,2$. We conclude that the length
of $\omega_n$ is bounded below by
\[
A\times\ve\alpha n^{1/4} > 2\alpha n^{1/4}\geq \wt
d_{\mathrm{gr}}\bigl(\gamma_n\bigl(\bigl\lfloor\alpha
n^{1/4}\bigr\rfloor\bigr),\gamma'_n\bigl(\bigl
\lfloor \alpha n^{1/4}\bigr\rfloor\bigr)\bigr),
\]
by our choice of $A$. This contradiction completes the proof.
\end{pf*}

In our applications, we will need a version of Lemma~\ref{keylem} where
the roles of the root vertex
and of the distinguished vertex $\partial$ of $\mm_n$ are
interchanged. To this end, we will
rely on a symmetry property of rooted and pointed $2p$-angulations that
we state in the next lemma.

%le4.2 #&#
\begin{lemma}
\label{coupling}
We can construct a random $2p$-angulation $M_n$ with two oriented edges
$e_n$ and $e^*_n$
and two distinguished vertices $\partial_n$ and $\partial^*_n$, in such
a way that:
\begin{enumerate}[(ii)]
\item[(i)] Both $(M_n,e_n,\partial_n)$ and $(M_n,e^*_n,\partial^*_n)$
are uniformly distributed
over roo\-ted and pointed $2p$-angulations with $n$ faces.
\item[(ii)] If $\rho_n$, respectively $\rho^*_n$, is the origin of $e_n$,
respectively of $e^*_n$, we have
\[
P\bigl(d_{\mathrm{gr}}\bigl(\partial_n,\rho^*_n\bigr)
> p\bigr)\leq\frac{2}{(p-1)n},\qquad P\bigl(d_{\mathrm{gr}}\bigl(
\partial^*_n,\rho_n\bigr) > p\bigr)\leq\frac{2}{(p-1)n}.
\]
\end{enumerate}
\end{lemma}

\begin{pf} We start from a uniformly distributed rooted and pointed
$2p$-angulation $M_n$ with $n$ faces,
given (as previously) as the image under the BDG bijection of
a uniformly distributed\vadjust{\goodbreak} labeled $p$-tree $(\tau_n,\ell^n_v)_{v\in
\tau
_n^\circ}$ with $n$ black vertices
and an independent Bernoulli variable $\ve$ with parameter $1/2$. We
let $e_n$
be the root edge of $M_n$, and as previously we write $\partial$
for the distinguished vertex of $M_n$. We can easily equip $M_n$ with
another distinguished oriented edge
$e^*_n$ by using the following device: We choose an independent random
variable $U_n$
uniformly distributed over $\{0,1,\ldots,pn-1\}$ and we let $e^*_n$ be
the edge
generated by the $U_n$th step of the BDG construction, oriented
uniformly at
random independently of $U_n$. In this way, and ``forgetting'' the distinguished
vertex $\partial$, we get a triplet $(M_n,e_n,e^*_n)$ which is
uniformly distributed over $2p$-angulations with $n$ faces and two
oriented edges.
Note that we distinguish the first oriented edge and the second one, so that
$(M_n,e_n,e^*_n)\not= (M_n,e^*_n,e_n)$ unless $e_n=e^*_n$. However, it
is easy to see (from the fact that a $2p$-angulation with $n$ faces has
always $pn$ edges) that the triplets $(M_n,e_n,e^*_n)$ and $ (M_n,e^*_n,e_n)$
have the same distribution.

Let $\rho_n$ and $\rho^*_n$ be the respective origins of $e_n$ and $e^*_n$.
Although $\rho^*_n$ is not uniformly distributed over the vertex set
$\mm_n$ of $M_n$, we can
construct a uniformly distributed random vertex $\partial_n$ that will be
close to $\rho^*_n$ with high probability. To do so, recall the
notation $C^n$ for the contour function, and $v^n_0,v^n_1,\ldots,v^n_{pn}$
for the contour sequence of $\tau^\circ_n$. If
$C^n_{U_n+1}\geq C^n_{U_n}$, let $\wh\partial_n$ be equal
to $v^n_{U_n+1}$. On the other hand, if
$C^n_{U_n}> C^n_{U_n+1}$, we let $\wh\partial_n$ be chosen uniformly
at random
among the $p-1$ children of the black vertex which is the parent
of $v^n_{U_n}$ in $\tau_n$. Then it is not hard to see that $\wh
\partial_n$ and
$\rho^*_n$ are on the boundary of the same face of $M_n$, so that
$d_{\mathrm{gr}}(\wh\partial_n,\rho^*_n)\leq p$. Furthermore, a~moment's
thought shows
that $\wh\partial_n$ is uniformly distributed over $\mm_n\setminus
\{
\varnothing,\partial\}$.
So, independently of the other random quantities, we may define
\[
\partial_n=\cases{ %
 \wh
\partial_n, &\quad $\mbox{with probability } \displaystyle\frac{(p-1)n}{(p-1)n+2},$
\vspace*{2pt}\cr
\varnothing,&\quad $\mbox{with probability } \displaystyle\frac{1}{(p-1)n+2},$
\vspace*{2pt}\cr
 \partial,&\quad $\mbox{with probability } \displaystyle\frac{1}{(p-1)n+2},$}
\]
so that the triplet $(M_n,e_n,\partial_n)$ is uniformly distributed over
rooted and pointed $2p$-angulations with $n$ faces, and the first bound
in (ii) holds by the preceding considerations.

To complete the proof, we select independently of $e_n^*$ and of the
other random quantities a random vertex
$\wt\partial_n$ uniformly distributed over $\mm_n$. Applying the
(inverse) BDG bijection to the (uniformly distributed)
rooted and pointed $2p$-angulation $(M_n,e^*_n,\wt\partial_n)$, we
can associate
with the edge $e_n$ a random variable $U^*_n$ uniformly distributed
over $\{0,1,\ldots,pn-1\}$ (just as $U_n$ was associated with $e^*_n$).
By duplicating the preceding argument, we then construct from $U^*_n$
a random vertex $\partial^*_n$ such that $(M_n,e^*_n,\partial^*_n)$
is uniformly distributed and the second bound in (ii) holds.
\end{pf}

If $\omega=(\omega(0),\omega(1),\ldots,\omega(k))$ is a path in
$M_n$, the length of $\omega$ is $|\omega|=k$, and the reversed path
$(\omega(|\omega|),\omega(|\omega|-1),\ldots,\omega(0))$ is
denoted by
$\ov\omega$.
Recall the constant $\alpha_0$ introduced in Lemma~\ref{keylem}.

%le4.3 #&#
\begin{lemma}
\label{keylem2}
Let $\alpha\in(0,\alpha_0)$ and $\beta_1,\beta_2$ such that
$15\alpha
<\beta_1<\beta_2$.
For every integer $n\geq1$, and every $\delta>0$,
consider the event
\begin{eqnarray*}
\e_{n,\delta}&=&\bigl\{\beta_1n^{1/4}<
\Delta_n<\beta_2 n^{1/4}\bigr\}\\
&&{} \cap \biggl\{ \wt
d_{\mathrm{gr}}\bigl(\ov\gamma_n\bigl(\bigl\lfloor\alpha
n^{1/4}\bigr\rfloor\bigr),\ov\gamma'_n\bigl(
\bigl\lfloor\alpha n^{1/4}\bigr\rfloor\bigr)\bigr)
,\\
&&\hspace*{18pt} \geq\wt d_{\mathrm{gr}}\biggl(\ov\gamma_n\biggl(\biggl\lfloor
\frac{\alpha}{3} n^{1/4}\biggr\rfloor\biggr)\ov\gamma'_n
\biggl(\biggl\lfloor\frac{\alpha}{3} n^{1/4}\biggr\rfloor\biggr)\biggr)
+\biggl(\frac{4\alpha}{3} -\delta\biggr)n^{1/4} \biggr\}.
\end{eqnarray*}
There exist
a real sequence $(\delta_n)_{n\geq0}$ decreasing to $0$, and a
constant $a_1>0$ such that
\[
P(\e_{n,\delta_n})\geq a_1
\]
for every sufficiently large integer $n$.

Furthermore, if the event $\e_{n,\delta_n}$ holds, we have also
%
%e16 #&#
\begin{eqnarray}
\label{keylemext2}&& \wt d_{\mathrm{gr}}\bigl(\ov\gamma'_n
\bigl(j'\bigr),\ov\gamma_n(j)\bigr)
\nonumber
\\[-8pt]
\\[-8pt]
\nonumber
&&\qquad\geq\wt
d_{\mathrm{gr}}\biggl(\ov\gamma'_n
\bigl(j'\bigr), \ov\gamma_n\biggl(\biggl\lfloor
\frac{\alpha}{3} n^{1/4}\biggr\rfloor\biggr)\biggr) + j- \biggl\lfloor
\frac{\alpha}{3} n^{1/4}\biggr\rfloor -\delta_n
n^{1/4}-2
\end{eqnarray}
for every $j\in\{\lfloor\frac{ \alpha}{3} n^{1/4}\rfloor,\ldots,
\lfloor
\alpha n^{1/4}\rfloor\}$ and
$j'\in\{0,1,\ldots, \lfloor\alpha n^{1/4}\rfloor\}$. The same bound
holds if the roles of
$\ov\gamma_n$ and $\ov\gamma'_n$ are interchanged.
\end{lemma}

\begin{pf} We start by proving the existence of a constant $a_2>0$ such
that, for every $\delta>0$,
%
%e17 #&#
\begin{equation}
\label{keylem200} \liminf_{n\to\infty} P(\e_{n,\delta}) \geq
a_2.
\end{equation}
The first part of the lemma follows: Just construct by induction a
monotone increasing sequence $(n_k)_{k\geq0}$ such that $P(\e_{n,2^{-k}})\geq a_2/2$
for every $n\geq n_k$, and put $\delta_n= 2^{-k}$ for $n_k\leq n< n_{k+1}$.

We consider a random $2p$-angulation $M_n$ with two oriented edges
$e_n$ and $e^*_n$
and two distinguished vertices $\partial_n$ and $\partial^*_n$, such
that properties
(i) and (ii) of Lemma~\ref{coupling} hold. We also write $\rho_n$ and
$\rho^*_n$
for the respective origins of $e_n$ and $e^*_n$. With the uniformly
distributed rooted and pointed
$2p$-angulation $(M_n,e_n,\partial_n)$ we associate the DMGB $\wt M_n$,
and similarly
with $(M_n,e^*_n,\partial^*_n)$ we associate the DMGB $\wt M_n^*$. We
write $\gamma_n$
and $\gamma'_n$, respectively $\gamma^*_n$ and ${\gamma^{*}_n}'$,
for the left and right boundary geodesics in $\wt M_n$, respectively in $\wt
M^*_n$. The common length
of $\gamma_n$ and $\gamma'_n$, respectively
of $\gamma^*_n$ and ${\gamma^*_n}'$, is denoted by $\Delta_n$, respectively by
$\Delta^*_n$. Notice that
$\gamma_n$, respectively\vadjust{\goodbreak} $\gamma^*_n$, can also be viewed as a geodesic in $M_n$,
which starts from one of the two vertices incident to $e_n$, respectively to
$e^*_n$, and ends
at $\partial_n$, respectively at $\partial^*_n$.

Thanks to property (ii) of Lemma~\ref{coupling}, we can use $\gamma_n^*$, or rather the
time-reversed path $\ov\gamma_n^*$, to construct an ``approximate''
geodesic from $\rho_n$ to
$\partial_n$ in~$M_n$. To do so, we concatenate a geodesic path from
$\rho_n$ to $\partial^*_n$ with the path $\ov\gamma_n^*$, and then
with a geodesic path from the initial point of $\gamma^*_n$ (which is
either $\rho^*_n$
or a neighbor of $\rho^*_n$) to $\partial_n$.
Let $\gamma_n^{**}$ be the path resulting from this concatenation.
From property~(ii) of Lemma~\ref{coupling}, the length
of $\gamma_n^{**} $ is bounded above by $d_{\mathrm{gr}}(\rho_n,\partial_n)+4(p+1)$, with
probability at least $1-2((p-1)n)^{-1}$. From Proposition 1.1 in \cite
{AM} (and the fact that this
result also holds for approximate geodesics as discussed in the
introduction of~\cite{AM}),
we know that the paths $\gamma_n$ and $\gamma_n^{**}$
must be close to each other with high probability when $n$ is large.
More precisely, we get
for every $\ve>0$,
%
%e18 #&#
\begin{equation}
\label{keylemtech1} \lim_{n\to\infty} P \Bigl[\max_{i\geq0}
d_{\mathrm{gr}}\bigl(\gamma_n(i\wedge \Delta_n),\ov
\gamma_n^* \bigl(i\wedge\Delta_n^*\bigr)\bigr) >\ve
n^{1/4} \Bigr] =0.
\end{equation}

In what follows, we suppose that the event considered in Lemma~\ref{keylem}
holds for the DMGB $\wt M_n$. We fix $\delta>0$ and using the property
(\ref{keylemtech1}),
we will show that we have on this event
%
%e19 #&#
\begin{eqnarray}
\label{keylemtech2}&& \wt d_{\mathrm{gr}}\bigl(\ov\gamma^*_n\bigl(
\bigl\lfloor\alpha n^{1/4}\bigr\rfloor\bigr),\ov\gamma^{*\prime}_n\bigl(\bigl\lfloor\alpha n^{1/4}
\bigr\rfloor\bigr)\bigr)
\nonumber
\\[-8pt]
\\[-8pt]
\nonumber
&&\qquad\geq\wt d_{\mathrm{gr}}\biggl(\ov\gamma^*_n
\biggl(\biggl\lfloor\frac{\alpha}{3} n^{1/4}\biggr\rfloor\biggr),
\ov\gamma^{*\prime}_n\biggl(\biggl\lfloor
\frac{\alpha}{3} n^{1/4}\biggr\rfloor \biggr)\biggr)
 +\biggl(\frac{4\alpha}{3} -\delta\biggr)n^{1/4},
\end{eqnarray}
except possibly on a set of probability tending to $0$ when $n\to
\infty
$. Our claim (\ref{keylem200})
will follow since the DMGBs $\wt M_n$ and $\wt M_n^*$ have the same
distribution [and $P(|\Delta_n-\Delta_n^*|\geq\ve n^{1/4})$
tends to $0$ as $n\to\infty$, for every $\ve>0$].

Without loss of generality, we can assume that $0<\delta<\alpha/2$.
We write $B_{M_n}(v,r)$, respectively $B_{\wt M_n}(v,r)$, $B_{\wt M_n^*}(v,r)$
for the open ball of radius $r$ centered at $v$ in $M_n$, respectively in $\wt
M_n$, in $\wt M_n^*$.
We first note that a geodesic path in $\wt M_n$ from
$\gamma'_n(\lfloor\alpha n^{1/4}\rfloor)$ to $\gamma_n(\lfloor
\alpha
n^{1/4}\rfloor)$
cannot visit $B_{\wt M_n}(\partial_n,\delta n^{1/4})$,
because otherwise its length would be at least $2(\Delta_n -\delta
n^{1/4}-\lfloor\alpha n^{1/4}\rfloor)
\geq2(\beta_1-\delta-\alpha)n^{1/4}$, which is clearly impossible. For
the same reason, a geodesic
path in $\wt M_n^*$ from
$\ov\gamma^{*\prime}_n(\lfloor\alpha n^{1/4}\rfloor)$ to $\ov
\gamma^*_n(\lfloor\alpha n^{1/4}\rfloor)$ does not
visit $B_{\wt M_n^*}(\gamma_n^*(0),\delta n^{1/4})$, except perhaps on
a set of
probability tending to $0$ as $n$ tends to infinity, which we may discard.

To simplify notation, set ${\cal R}(\gamma_n)=\{\gamma_n(i)\dvtx 0\leq
i\leq\Delta_n\}$. Let $i_0\in\{1,\ldots,\break\Delta_n-1\}$.
A path $\omega=(\omega(i),0\leq i\leq|\omega|)$ in $M_n$ is said to be
an \textit{admissible loop}
from $\gamma_n(i_0)$ in $M_n$ if $\omega(0)=\omega(|\omega|)=
\gamma_n(i_0)$, if
$\omega$ does not visit $\gamma_n(0)$ or $\gamma_n(\Delta_n)$, and
if there exist integers $k_\omega$ and $\ell_\omega$ such that
$0\leq
k_\omega<\ell_\omega\leq|\omega|$,
and the
following holds:
\begin{longlist}[(iii)]
\item[(i)] $\omega(i)\in{\cal R}(\gamma_n)$ if and only if $0\leq
i\leq
k_\omega$ or $\ell_\omega\leq i\leq|\omega|$;
\item[(ii)] $\omega(k_\omega)$ is connected to $\omega(k_\omega
+1)$ by
an edge starting from a corner
belonging to the left side of $\gamma_n$ [here and later, $\gamma_n$ is
oriented from
$\gamma_n(0)$ to $\gamma_n(\Delta_n)$];
\item[(iii)] $\omega(\ell_\omega)$ is connected to $\omega(\ell_\omega
-1)$ by an edge starting from a corner
belonging to the right side of $\gamma_n$.
\end{longlist}
In the same way, we define a $*$-\textit{admissible loop} in $M_n$ by
replacing $\gamma_n$
with $\ov\gamma^*_n$ and $\Delta_n$ with $\Delta_n^*$ everywhere in
the previous
definition. Note that an admissible loop (resp., a $*$-admissible loop)
winds exactly once in clockwise order around
the point $\gamma_n(\Delta_n)=\partial_n$ [resp., around $\ov\gamma_n^*(\Delta_n^*)=\gamma_n^*(0)$]
in the twice punctured sphere ${\sphere}^2\setminus\{\gamma_n(0),
\gamma_n(\Delta_n)\}$
(resp., in ${\sphere}^2\setminus\{\gamma^*_n(0), \gamma^*_n(\Delta^*_n)\}$).

The following properties are easily checked from the relation between $M_n$
and $\wt M_n$ or $\wt M^*_n$: %
\begin{longlist}[(a)]
\item[(a)] If $\omega$ is an admissible loop from $\gamma_n(i_0)$,
then we can find a path $\wt\omega$ in $\wt M_n$ such that $\wt
\omega
(0)=\gamma'_n(i_0)$,
$\wt\omega(|\wt\omega|)=\gamma_n(i_0)$ and $|\wt\omega|=|\omega|$.
\item[(b)] Similarly, if $\omega^*$ is a $*$-admissible loop from
$\ov
\gamma_n^*(i_0)$,
then we can find a path $\wt\omega^*$ in $\wt M_n^*$ such that $\wt
\omega(0)=\ov\gamma^{*\prime}_n(i_0)$,
$\wt\omega(|\wt\omega|)=\ov\gamma_n^*(i_0)$ and $|\wt\omega^*|=|\omega^*|$.
\item[(c)] Let $\wt\omega$ be a path in $\wt M_n$ that does not visit
$\gamma_n(0)$ or $\gamma_n(\Delta_n)$, such that
$\wt\omega(0)=\gamma_n'(i_0)$ and $\wt\omega(|\wt\omega|)=\gamma_n(i_0)$, and such that
$\wt\omega$ visits $\gamma'_n(i)$ and $\gamma_n(i)$ in this order,
for some
$i\in\{1,\ldots,\Delta_n-1\}$. Then we can
find an admissible loop $\omega$ from $\gamma_n(i_0)$ such that
$|\omega
|\leq|\wt\omega|$,
and such that $\gamma_n(i)\in\{\omega(j)\dvtx 0\leq j\leq k_\omega\}
\cap\{
\omega(j)\dvtx \ell_\omega\leq j\leq|\omega|\}$. If,
for some $r,r'>0$, $\wt\omega$ does not visit $B_{\wt M_n}(\gamma_n(0),r)
\cup B_{\wt M_n}(\gamma_n(\Delta_n),r')$, then $\omega$
can be constructed so that it does not visit $B_{M_n}(\gamma_n(0),r)
\cup B_{M_n}(\gamma_n(\Delta_n),r')$.
\item[(d)] Similarly, if $\wt\omega^*$ is any path in $\wt M_n^*$ that
does not
visit $\ov\gamma_n(0)$ or $\ov\gamma_n(\Delta_n^*)$ and is such that
$\wt\omega(0)=\ov\gamma^{*\prime}_n(i_0)$ and $\wt\omega^*(|\wt\omega^*|)={\ov\gamma_n^*}(i_0)$,
then we can find a $*$-admissible loop $\omega^*$ from $\ov\gamma_n^*(i_0)$ such that
$|\omega^*|\leq|\wt\omega^*|$. If, for some $r,r'>0$, $\wt\omega^*$
does not visit $B_{\wt M_n^*}(\ov\gamma_n^*(0),r)
\cup B_{\wt M_n^*}(\ov\gamma_n(\Delta_n^*),r')$, then $\omega^*$
can be constructed so that it does not visit $B_{M_n}(\ov\gamma_n^*(0),r)
\cup B_{M_n}(\ov\gamma_n(\Delta_n^*),r')$.
\end{longlist}

Let $\wt\omega_n^*$ be a geodesic in $\wt M_n^*$ from
$\ov\gamma^{*\prime}_n(\lfloor\alpha n^{1/4}\rfloor)$ to $\ov
\gamma^*_n(\lfloor\alpha n^{1/4}\rfloor)$.
As mentioned above, we may assume that $\wt\omega_n^*$
does not
visit $B_{\wt M_n^*}(\gamma_n^*(0),\delta n^{1/4})$. Also, if
$\wt\omega_n^*$ visits $B_{\wt M_n^*}(\partial_n^*,\frac{\delta
}{4} n^{1/4})$,
then it readily follows that (\ref{keylemtech2}) holds. So we may
restrict our attention to the
case when $\wt\omega_n^*$ does not visit this ball.

By property (d) above, we can construct a $*$-admissible loop $\omega^*$ from
$\ov\gamma_n^*(\lfloor\alpha n^{1/4}\rfloor)$, such that $|\omega^*|
\leq\wt d_{\mathrm{gr}}(\ov\gamma^*_n(\lfloor\alpha n^{1/4}\rfloor),\ov\gamma^{*\prime}_n(\lfloor\alpha n^{1/4}\rfloor))$
and such that $\omega^*$
does not visit $B_{ M_n}(\gamma_n^*(0),\delta n^{1/4})$
or $B_{M_n}(\partial_n^*,\frac{\delta}{4} n^{1/4})$.
Now pick a geodesic $g_n$ from $\gamma_n(\lfloor\alpha n^{1/4}\rfloor)$
to $\ov\gamma^*_n(\lfloor\alpha n^{1/4}\rfloor)$ and note that,
thanks to
(\ref{keylemtech1}), we have $|g_n|\leq\frac{\delta}{16} n^{1/4}$,
except on
a set of probability tending to $0$ which we may discard. By concatenating
$g_n$, $\omega^*$
and the time-reversed path $\ov g_n$, we get a path $\omega$
whose length is bounded above by $|\omega^*| +\frac{\delta
}{8}n^{1/4}$, and
which starts and ends at $\gamma_n(\lfloor\alpha n^{1/4}\rfloor)$.
Moreover, we can also use property (ii) of Lemma~\ref{coupling} to get
that, except on an event of probability tending to $0$
as $n\to\infty$, $\omega$ does not visit
the balls $B_{M_n}(\partial_n,\frac{\delta}{2} n^{1/4})$ and
$B_{M_n}(\gamma_n(0),\frac{\delta}{8} n^{1/4})$.
Of course,
$\omega$ need not be an admissible loop. However,
since $\omega^*$ is a $*$-admissible loop and therefore winds exactly
once in clockwise order
around $\ov\gamma_n^*(\Delta_n^*)=\gamma_n^*(0)$
in the twice punctured sphere ${\sphere}^2\setminus\{\gamma^*_n(0),
\gamma^*_n(\Delta^*_n)\}$,
it follows that $\omega$ also winds exactly once around $\gamma_n(\Delta_n)=\partial_n$
in clockwise order in ${\sphere}^2\setminus\{\gamma_n(0), \gamma_n(\Delta_n)\}$.
Hence, a simple topological
argument shows that there must exist a subinterval $\{k_n,\ldots,\ell_n\}$
of $\{0,1,\ldots,|\omega|\}$, with $k_n<\ell_n$,
such that $\omega(k_n)\in{\cal R}(\gamma_n)$,
$\omega(\ell_n)\in{\cal R}(\gamma_n)$, $\omega(i)\notin{\cal
R}(\gamma_n)$
if $k_n<i<\ell_n$, $\omega(k_n)$ is connected to $\omega(k_n+1)$
by an edge that starts from the left side of $\gamma_n$
and $\omega(\ell_n)$ is connected to $\omega(\ell_n-1)$
by an edge that starts from the right side of $\gamma_n$. It follows
that we can find an admissible loop $\omega'$ from $\gamma(\lfloor
\alpha n^{1/4}\rfloor)$
such that $|\omega'|\leq|\omega| \leq|\omega^*| +\frac{\delta
}{8}n^{1/4}$.
By property (a), we have then
%
%e20 #&#
\begin{eqnarray}
\label{keylemtech3} && \wt d_{\mathrm{gr}}\bigl(\gamma_n\bigl(\bigl
\lfloor\alpha n^{1/4}\bigr\rfloor\bigr),\gamma'_n
\bigl(\bigl\lfloor \alpha n^{1/4}\bigr\rfloor\bigr)\bigr)
\nonumber
\\[-8pt]
\\[-8pt]
\nonumber
&&\qquad \leq\bigl|\omega'\bigr| \leq \wt d_{\mathrm{gr}}\bigl(\ov
\gamma^*_n\bigl(\bigl\lfloor\alpha n^{1/4}\bigr\rfloor\bigr),
\ov \gamma^{*\prime}_n\bigl(\bigl\lfloor\alpha
n^{1/4}\bigr\rfloor\bigr)\bigr) +\frac{\delta}{8}n^{1/4}.
\nonumber
\end{eqnarray}

Now recall that we are arguing on the event of Lemma~\ref{keylem}.
Hence, we know that there exists
a geodesic path $\wt\omega_n$ in $\wt M_n$ from $\gamma'_n(\lfloor
\alpha n^{1/4}\rfloor)$ to $\gamma_n(\lfloor\alpha n^{1/4}\rfloor)$,
that visits $\gamma'_n(\lfloor\frac{\alpha}{3}n^{1/4}\rfloor)$ and
$\gamma_n(\lfloor\frac{\alpha}{3}n^{1/4}\rfloor)$
in this order. We already noticed that $\wt\omega_n$ does not visit the
ball $B_{\wt M_n}(\partial_n,\delta n^{1/4})$. If $\wt\omega_n$ visits
the ball
$B_{\wt M_n}(\gamma_n(0), \frac{\delta}{8}n^{1/4})$, then
\[
\wt d_{\mathrm{gr}}\bigl(\gamma_n\bigl(\bigl\lfloor\alpha
n^{1/4}\bigr\rfloor\bigr),\gamma'_n\bigl(\bigl
\lfloor \alpha n^{1/4}\bigr\rfloor\bigr)\bigr)\geq2\bigl\lfloor\alpha
n^{1/4}\bigr\rfloor- \frac{\delta}{4} n^{1/4}
\]
and it follows from (\ref{keylemtech3}) that
\[
\wt d_{\mathrm{gr}}\bigl(\ov\gamma^*_n\bigl(\bigl\lfloor\alpha
n^{1/4}\bigr\rfloor\bigr),\ov \gamma^{*\prime}_n
\bigl(\bigl\lfloor\alpha n^{1/4}\bigr\rfloor\bigr)\bigr) \geq2\bigl
\lfloor\alpha n^{1/4}\bigr\rfloor- \frac{\delta}{2} n^{1/4},
\]
from which (\ref{keylemtech2}) is immediate. So we may assume that
$\wt
\omega_n$ does not
visit $B_{\wt M_n}(\gamma_n(0), \frac{\delta}{8}n^{1/4})$.
It follows from property (c) that there exists
an admissible loop $\omega_n$ from $\gamma_n(\lfloor\alpha
n^{1/4}\rfloor)$, which visits
$\gamma_n(\lfloor\frac{\alpha}{3} n^{1/4}\rfloor)$ both between times~$0$ and $k_{\omega_n}$
and between times $\ell_{\omega_n}$ and $|\omega_n|$, and has length
$|\omega_n|\leq|\wt\omega_n|=\wt d_{\mathrm{gr}}(\gamma_n(\lfloor\alpha
n^{1/4}\rfloor),  \gamma'_n(\lfloor\alpha n^{1/4}\rfloor))$.
Moreover, $\omega_n$ does not visit $B_{M_n}(\partial_n,\delta n^{1/4})$
or $B_{M_n}(\gamma_n(0), \frac{\delta}{8}n^{1/4})$.
Let $p_n\in\{0,1,\ldots,k_{\omega_n}\}$ and
$q_n\in\{\ell_{\omega_n},\ldots,|\omega_n|\}$ such that $\omega_n(p_n)=\omega(q_n)= \gamma_n(\lfloor\frac{\alpha}{3}
n^{1/4}\rfloor)$.
Notice that necessarily
\[
q_n-p_n\leq|\omega_n| - 2\biggl(\bigl\lfloor
\alpha n^{1/4}\bigr\rfloor- \biggl\lfloor \frac
{\alpha}{3}
n^{1/4}\biggr\rfloor\biggr).
\]
Let $h_n$ be a geodesic path in $M_n$ from $\ov\gamma_n^*(\lfloor
\frac
{\alpha}{3} n^{1/4}\rfloor)$ to $ \gamma_n(\lfloor\frac{\alpha}{3}
n^{1/4}\rfloor)$, and let $\omega'_n$ be the path obtained by concatenating
$h_n$, $(\omega_n(p_n+i),0\leq i\leq q_n-p_n)$ and $\ov h_n$ in this
order. Notice that by\vadjust{\goodbreak}
(\ref{keylemtech1}), we have $|\omega'_n|\leq q_n-p_n+ \frac{\delta}{4}
n^{1/4}$ outside a set
of probability tending to $0$ as $n\to\infty$, which we may discard. By
the same topological argument as previously, we
can find a $*$-admissible loop $\omega^*_n$
from $\ov\gamma_n^*(\lfloor\frac{\alpha}{3} n^{1/4}\rfloor)$
such that $|\omega^*_n|\leq|\omega'_n|$. By property (b), we have now
%
%e21 #&#
\begin{eqnarray}
\label{keylemtech4}
&&\wt d_{\mathrm{gr}}\biggl(\ov\gamma^*_n
\biggl(\biggl\lfloor\frac{\alpha}{3} n^{1/4}\biggr\rfloor \biggr),\ov\gamma^{*\prime}_n\biggl(\biggl\lfloor\frac{\alpha}{3}
n^{1/4}\biggr\rfloor\biggr)\biggr)\nonumber \\
&&\qquad\leq\bigl|\omega^*_n\bigr|\leq\bigl|
\omega'_n\bigr|\leq q_n-p_n+
\frac{\delta}{4} n^{1/4}
\\
&&\qquad\leq\wt d_{\mathrm{gr}}\bigl(\gamma_n\bigl(\bigl\lfloor\alpha
n^{1/4}\bigr\rfloor\bigr),\gamma'_n\bigl(\bigl
\lfloor\alpha n^{1/4}\bigr\rfloor\bigr)\bigr) - 2\biggl(\bigl\lfloor
\alpha n^{1/4}\bigr\rfloor- \biggl\lfloor\frac{\alpha}{3}
n^{1/4}\biggr\rfloor\biggr)+ \frac{\delta}{4} n^{1/4}.
\nonumber
\end{eqnarray}
We now notice that (\ref{keylemtech2}) follows from (\ref{keylemtech3})
and (\ref{keylemtech4}),
which completes the proof of~(\ref{keylem200}) and of the first part of
the lemma.

To prove the second part of the lemma, we first note that if $\e_{n,\delta_n}$ holds,
we have also, for every integer $i,j,i',j'$ such that $\lfloor\frac{
\alpha}{3} n^{1/4}\rfloor\leq i\leq j\leq\lfloor\alpha
n^{1/4}\rfloor$
and $\lfloor\frac{ \alpha}{3} n^{1/4}\rfloor\leq i'\leq j'\leq
\lfloor
\alpha n^{1/4}\rfloor$,
%
%e22 #&#
\begin{equation}
\label{keylemext} \qquad\wt d_{\mathrm{gr}}\bigl(\ov\gamma_n(j),\ov
\gamma'_n\bigl(j'\bigr)\bigr) \geq\wt
d_{\mathrm{gr}}\bigl(\ov\gamma_n(i),\ov\gamma'_n
\bigl(i'\bigr)\bigr) +j+j'-i-i' -
\delta_n n^{1/4}-2.
\end{equation}
This immediately follows from the triangle inequality, which gives
\[
\wt d_{\mathrm{gr}}\bigl(\ov\gamma_n(j),\ov\gamma'_n
\bigl(j'\bigr)\bigr)+ 2\bigl\lfloor\alpha n^{1/4}\bigr
\rfloor- j- j' \geq\wt d_{\mathrm{gr}}\bigl(\ov
\gamma_n\bigl(\bigl\lfloor\alpha n^{1/4}\bigr\rfloor\bigr),
\ov \gamma'_n\bigl(\bigl\lfloor\alpha n^{1/4}
\bigr\rfloor\bigr)\bigr)
\]
and
\[
\wt d_{\mathrm{gr}}\bigl(\ov\gamma_n(i),\ov\gamma'_n
\bigl(i'\bigr)\bigr)\leq\wt d_{\mathrm{gr}}\biggl(\ov
\gamma_n\biggl(\biggl\lfloor\frac{\alpha}{3} n^{1/4}\biggr
\rfloor\biggr),\ov\gamma'_n\biggl(\biggl\lfloor
\frac
{\alpha}{3} n^{1/4}\biggr\rfloor\biggr)\biggr) + i+i'
- 2\biggl\lfloor\frac{\alpha}{3} n^{1/4}\biggr\rfloor.
\]
Then, if $\lfloor\frac{\alpha}{3} n^{1/4}\rfloor\leq j'\leq\lfloor
\alpha n^{1/4}\rfloor$, the
bound (\ref{keylemext2}) follows from the special case $i'=j'$,
$i=\lfloor\frac{\alpha}{3} n^{1/4}\rfloor$
in (\ref{keylemext}). If $0\leq j'< \lfloor\frac{\alpha}{3}
n^{1/4}\rfloor$, consider a geodesic from
$\ov\gamma'_n(j')$ to $\ov\gamma_n(j)$, and a geodesic from $\ov
\gamma'(\lfloor\frac{\alpha}{3} n^{1/4}\rfloor)$
to $\ov\gamma(\lfloor\frac{\alpha}{3} n^{1/4}\rfloor)$,
and observe that, by a topological argument,
these two geodesics must intersect, say, at a vertex~$v$. Because $v$
belongs to a geodesic from
$\ov\gamma'(\lfloor\frac{\alpha}{3} n^{1/4}\rfloor)$
to $\ov\gamma(\lfloor\frac{\alpha}{3} n^{1/4}\rfloor)$,
the case $j'=\lfloor\frac{\alpha}{3} n^{1/4}\rfloor$ in (\ref
{keylemext2}) easily gives
\[
\wt d_{\mathrm{gr}}\bigl(v,\ov\gamma_n(j)\bigr)\geq \wt
d_{\mathrm{gr}}\biggl(v, \ov\gamma_n\biggl(\biggl\lfloor
\frac{\alpha}{3} n^{1/4}\biggr\rfloor \biggr)\biggr) + j- \biggl\lfloor
\frac{\alpha}{3} n^{1/4}\biggr\rfloor-\delta_n
n^{1/4}-2.
\]
Then, by adding $\wt d_{\mathrm{gr}}(\ov\gamma'_n(j'),v)$ to both sides of this
inequality, we arrive at the desired
bound (\ref{keylemext2}) also in the case $0\leq j'< \lfloor\frac
{\alpha}{3} n^{1/4}\rfloor$. The case when
the roles of
$\ov\gamma_n$ and $\ov\gamma'_n$ are interchanged is treated similarly.
\end{pf}

%s5 #&#
\section{The main estimate}
\label{mainEstim}

In this section and in the next two ones, we consider the setting of
Theorem~\ref{mainIM}, and we assume that the
convergence (\ref{basicconv}) holds almost surely along a suitable
sequence $(n_k)_{k\geq0}$. We use the\vadjust{\goodbreak}
notation introduced at the beginning of Section~\ref{convBrma}.
Recall from Section~\ref{GeoBm} the notation $\Gamma=(\Gamma
(r),0\leq r\leq\Delta)$ for the
geodesic from $\bp(0)$ to $\bp(s_*)$ in $\mm_\infty$. For every
$r\in
[0,\Delta]$, we have $\Gamma(r)=\bp(S_r)$, where
\[
S_r:=\inf\{s\geq0\dvtx Z_s=-r\}.
\]
Let $r>0$ and argue under the conditional probability measure $P(\cdot
\mid\Delta\geq r)$.
The main result of this section (Lemma~\ref{mainest}) shows that, with
a probability close to $1$ when $\ve>0$ is small,
for every point $z$ of $\mm_\infty$ ``sufficiently far''
from $\Gamma(r)$, either there is a geodesic from $z$ to $\Gamma(r)$
that visits $\Gamma(r-\ve)$
or there is a geodesic from $z$ to $\Gamma(r)$ that visits $\Gamma
(r+\ve)$.

We fix $\mu\in(0,1/2)$, $A>\mu$ and $\kappa\in(0,1/4)$.
We assume that $\mu\leq r\leq A$. The forthcoming estimates will depend
on $\mu, A$ and $\kappa$, but not
on the choice of $r$ in the interval $[\mu,A]$.

We start by introducing some notation. For every $\delta\in(0,r)$, we set
\[
\eta_\delta(r):=\inf\Bigl\{s\geq S_r\dvtx \be_s=\min_{t\in[S_r,s]} \be_t\mbox{ and }
Z_s=-r+\delta\Bigr\}.
\]
In other words, $\eta_\delta(r)$ is the first instant $s$ after $S_r$
such that $p_\be(s)$
belongs to the ancestral line of $p_\be(S_r)$ and has label $-r+\delta
$. We may
also say that $\be_{S_r}-\be_{\eta_\delta(r)}$ is the minimal distance
needed when
moving from $p_\be(S_r)$ toward the root of ${\cal T}_\be$ in order
to meet a
vertex with label $-r+\delta$.

In order to state our first lemma, we need to introduce the subtrees
that branch from
the ``right side'' of the ancestral line of $p_\be(S_r)$. Formally, we consider
all (nonempty) open subintervals $(u,u')$ of $[S_r,1]$ that satisfy the property
\[
\be_u=\be_{u'}=\min_{t\in[S_r,u']}\be_t,
\]
which automatically implies that $\be_t>\be_u$ for every $t\in(u,u')$
(otherwise, this would contradict the fact
that the local minima of $\be$ are distinct).
We write $(u_{(i)},u'_{(i)})_{i\in I}$ for the collection of all these
intervals. For each $i\in I$,
we will be interested especially in the quantities
$Z_{u_{(i)}}=Z_{u'_{(i)}}$, representing the
label at the root of the subtree, and
\[
\Delta^{(i)}:= Z_{u_{(i)}} - \min_{s\in[u_{(i)},u'_{(i)}]}
Z_s,
\]
representing (minus) the minimal relative label in the subtree.

Recall the constant $\alpha_0$ introduced in Lemma~\ref{keylem}. We fix
$\alpha\in(0,1/10)$ such that $\alpha/\kappa_p <\alpha_0$. We
start by choosing
four positive constants $\alpha_1,\alpha_2,\alpha'_2,\wt\alpha$
such that
\[
\cases{ %
2\alpha_1+\wt\alpha<\alpha,
\vspace*{2pt}\cr
\displaystyle\alpha_2-\alpha'_2 >\frac{\alpha}{3},
\vspace*{2pt}\cr
\alpha_1>\alpha_2.}
\]
It is easy to verify that such a choice is possible. We then choose
$\beta_1,\beta_2 \in(2,4)$ such that
$\beta_2>\beta_1$ and $\alpha_1+\beta_1 > \alpha_2+\beta_2$. We finally
fix a constant $\lambda\in(0,1)$ such that
\[
\frac{\alpha}{3} (1+\lambda)^{1/4} < \alpha_2-
\alpha'_2\quad \mbox{and}\quad \alpha_2+
\beta_2<(1+\lambda)^{1/4} (\alpha_1+
\beta_1),
\]
and an integer $K\geq2$ such that $K\geq\alpha_1/\alpha'_2$. We
then set
\[
\ell_0:= \biggl\lfloor\frac{\log(1/\mu)}{\log K}\biggr\rfloor+3
\]
in such a way that $K^{-\ell_0+2}<\mu$.

For every integer $\ell\geq\ell_0$, we say that the event $E_\ell$
holds if $S_r<1-\kappa$, and if there exists an index $i\in I$ such that:
\begin{longlist}[(a)]
\item[(a)] $\eta_{K^{-\ell+1}}(r)<u_{(i)}<u'_{(i)}< \eta_{K^{-\ell
+2}}(r) < 1- \frac{\kappa}{2} $;
\item[(b)] $(\alpha_2+\beta_2) K^{-\ell} < \Delta^{(i)} <(\alpha_1 +
\beta_1) K^{-\ell}$;
\item[(c)] $-r + \beta_1 K^{-\ell} < Z_{u_{(i)}} < -r + \beta_2
K^{-\ell
} $;
\item[(d)] $\min\{Z_s\dvtx s\in[S_r,u_{(i)}]\cup[u'_{(i)}, \eta_{K^{-\ell
+2}}(r)]\} > -r -\alpha'_2 K^{-\ell} $;
\item[(e)] there exists a vertex $b$ of ${\cal T}_\be$ that belongs
to the
ancestral line
of $p_\be(S_r)$ in~${\cal T}_\be$, such that $\be_{\eta_{K^{-\ell
+2}}(r)}\leq
d_\be(\rho,b)\leq\be_{u_{(i)}}$,
and $Z_b < -r + \wt\alpha K^{-\ell} $;
\item[(f)] $K^{-4\ell} < u'_{(i)}-u_{(i)} < (1+\lambda) K^{-4\ell} $.
\end{longlist}
The meaning of these conditions will appear
more clearly in the forthcoming proofs. Informally, noting that the
index $i\in I$
corresponds to a subtree branching from the right side of the ancestral line
of $p_\be(S_r)$, condition (a) gives information about the level at
which this subtree branches,
and condition (f) provides bounds on its size. Condition (b) gives
bounds on the relative minimal label
of the subtree, and condition (c) is concerned with the label of its
root. Condition (d) gives (in particular)
a lower bound on the minimum of the labels ``between'' the subtree and
the vertex $p_\be(S_r)$.
Finally, condition (e), which seems mysterious at this point, will be
used together with the construction
of edges in the BDG bijection to get an upper bound for the minimal
label on a path before it enters
the subtree. Of course the choice of the various constants that appear
in (a)--(f) is made in an appropriate manner
in view of the proof of our main estimate (Lemma~\ref{mainest}
below).\looseness=-1

We note that conditions (b) and (c) imply that
\[
\min\bigl\{Z_s\dvtx s\in\bigl[u_{(i)},u'_{(i)}
\bigr]\bigr\} < -r-\alpha_2 K^{-\ell}
\]
and since $\alpha'_2<\alpha_2$, conditions (a) and (d) show that there
can be at most an
index~$i$ satisfying (a)--(f). When $E_\ell$ holds, we will write
$(u_\ell,u'_\ell)=(u_{(i)},u'_{(i)})$
and $\Delta^\ell= \Delta^{(i)}$,
where $i$ is the unique index such that properties (a)--(f) hold.

%le5.1 #&#
\begin{lemma}
\label{snakelemma}
For every given $a\in(0,1)$, we can find a constant $\ov a\in(0,1)$ and
another constant $\ov C$, which both depend on\vadjust{\goodbreak}
$\mu,A$ and $\kappa$ but not on the choice of $r\in[\mu,A]$,
such that, for every integer $\ell\geq2\ell_0$,
\[
E \bigl[ \mathbf{1}_{\{S_r<1-\kappa\}} a^{\sum_{k=\lfloor\ell/2\rfloor
}^\ell\mathbf{1}_{E_k}} \bigr] \leq\ov C \ov
a^\ell.
\]
\end{lemma}

We postpone the proof of this lemma to the \hyperref[app]{Appendix}. One can use
scaling arguments to see that the
probability of $E_k$ is bounded below by a positive constant. If the
events $E_k$, $k\geq\ell_0$ were independent
under $P(\cdot\mid S_r<1-\kappa)$, the
bound of the lemma would immediately follow. The events $E_k$ are not
independent,
in particular, because of condition (d), but in some sense there is
enough independence
to ensure that the bound of the lemma holds.

We now want to take advantage of the (almost sure) convergence (\ref
{basicconv})
to get that if the event $E_\ell$ holds, for some $\ell\geq\ell_0$, the
discrete labeled $p$-trees $\theta_n$ satisfy properties analogous to
(a)--(f) at least for all
sufficiently large values of $n$ in the sequence $(n_k)_{k\geq0}$.
From now on
until the end of this section, we consider only
values of $n$ in this sequence. We put $r_n=\lfloor r \frac
{n^{1/4}}{\kappa_p}\rfloor$ and
\[
\sigma_n:=\min\bigl\{i\geq0\dvtx \Lambda^n_i
=-r_n\bigr\},
\]
where $\min\varnothing=\infty$.
On the event $\{S_r<\infty\}$, we have
%
%e23 #&#
\begin{equation}
\label{convhitting} \lim_{n\to\infty} \frac{\sigma_n}{pn} = S_r\qquad \mbox{a.s.}
\end{equation}
This follows from the (easy) property $\min_{S_r\leq s\leq S_r+\ve} Z_s
< -r$, for every $\ve>0$, a.s.

To simplify notation, we put $\ov v_n=v^n_{\sigma_n}$. Note that we
have also $\ov v_n=\gamma_n(r_n)$,
where (in agreement with previous notation) $\gamma_n$ is the simple
geodesic from the
first corner of $\varnothing$ to $\partial$ in $M_n$.
We define $y_{n,i}$ as the (white) ancestor of $\ov v_n$ at generation
$2i$, for
every $i\in\{0,1,\ldots, \frac{1}{2}| \ov v_n|\}$. For every $\delta
\in
(0,r)$ we also set
\[
\psi_{n,r}(\delta)=\max\biggl\{i\dvtx \ell^n_{y_{n,i}} >
-r_n+\delta\frac
{n^{1/4}}{\kappa_p}\biggr\},
\]
and we let $\Psi_{n,r}(\delta)$ be the index corresponding to the
last visit
of the vertex $y_{n,\psi_{n,r}(\delta)}$ by the contour sequence of
$\tau_n^\circ$.

If $\tau$ is a subtree of $\tau_n$ branching from the right side of the
ancestral line of~$\ov v_n$,
we write $z(\tau)$ for the root of $\tau$ (this is either a vertex of
the form $y_{n,i}$ or a ``brother'' of such a vertex) and $[r(\tau
),r'(\tau)]$ for the
interval corresponding to visits of $\tau$ in the contour sequence of
$\tau_n^\circ$, and we also
let $\Delta(\tau)$ be equal to $1$ minus the minimal relative label of
white vertices of $\tau$---as previously the
relative label of a white vertex in $\tau$ is the label of this vertex
minus the label of $z(\tau)$.

Then, for every $\ell\geq\ell_0$, we say that the event $E_{n,\ell}$
holds if $\sigma_n<\infty$
and if there exists a subtree $\tau$ of $\tau_n$ branching from the
right side of the ancestral line of $\ov v_n$,
such that:
\begin{longlist}[(a$'$)]
\item[(a$'$)] $\psi_{n,r}(K^{-\ell+2}) < \frac{1}{2} |z(\tau)|< \psi_{n,r}(K^{-\ell+1})$;
\item[(b$'$)] $(\alpha_2 + \beta_2)K^{-\ell} \frac{n^{1/4}}{\kappa_p} <
\Delta(\tau)
< (\alpha_1 + \beta_1)K^{-\ell} \frac{n^{1/4}}{\kappa_p} $;
\item[(c$'$)] $-r_n+ \beta_1K^{-\ell} \frac{n^{1/4}}{\kappa_p}< \ell^n_{z(\tau)}< -r_n + \beta_2K^{-\ell} \frac{n^{1/4}}{\kappa_p} $;
\item[(d$'$)] $\min\{\ell^n_{v^n_i}\dvtx i\in[\sigma_n,r(\tau)] \cup
[r'(\tau
),\Psi_{n,r}(K^{-\ell+2})]\}
> -r_n-\alpha'_2K^{-\ell} \frac{n^{1/4}}{\kappa_p} $;
\item[(e$'$)] there exists an index $j\in\{\psi_{n,r}(K^{-\ell
+2}),\ldots,\frac{1}{2}|z(\tau)|-1\}$ such that
$\ell^n_{y_{n,j}} < -r_n+\wt\alpha K^{-\ell} \frac{n^{1/4}}{\kappa_p} $;
\item[(f$'$)] $K^{-4\ell}n < |\tau| < (1+\lambda) K^{-4\ell}n $.
\end{longlist}
Of course, (a$'$)--(f$'$) are just discrete analogs of (a)--(f). The same
argument as above shows that there can be at most one subtree $\tau$
satisfying conditions (a$'$)--(f$'$). If the event $E_{n,\ell}$ holds, we
denote this
subtree by $\tau_{n,\ell}$ and we write $z_{n,\ell}=z(\tau_{n,\ell})$
$r_{n,\ell}=r(\tau_{n,\ell})$, $r'_{n,\ell}= r'(\tau_{n,\ell})$
and $\Delta_{n,\ell}= \Delta(\tau_{n,\ell})$ to simplify notation.
Note that~(b$'$) and~(c$'$) imply
%
%e24 #&#
\begin{eqnarray}
\label{conseq1} \min_{v\in\tau_{n,\ell}} \ell^n_v &=&
\min_{r_{n,\ell}\leq i\leq
r'_{n,\ell}} \ell^n_{v^n_i}
\nonumber
\\[-8pt]
\\[-8pt]
\nonumber
&\in& \biggl(-r_n-
\alpha_1 K^{-\ell} \frac{n^{1/4}}{\kappa_p}+1,-r_n-
\alpha_2 K^{-\ell} \frac{n^{1/4}}{\kappa_p}+1 \biggr).
\end{eqnarray}

From the almost sure convergence (\ref{basicconv}) and straightforward
arguments, we get that
%
%e25 #&#
\begin{equation}
\label{estimtech1} E_\ell\subset\liminf_{n\to\infty}
E_{n,\ell}\qquad \mbox{a.s.}
\end{equation}
Hence, if $E_\ell$ holds (and discarding a set of probability zero), we
know that
$E_{n,\ell}$ also holds for all sufficiently large $n$ and,
furthermore, one has
\[
\lim_{n\to\infty} \frac{r_{n,\ell}}{pn} = u_\ell,\qquad \lim_{n\to
\infty} \frac{r'_{n,\ell}}{pn} = u'_\ell.
\]

As explained in Section~\ref{labtrees}, we may associate with $\tau_{n,\ell}$
a labeled $p$-tree $\theta_{n,\ell}$, by renaming the vertices and
subtracting the label of the
root $z_{n,\ell}$ from all labels. With this labeled $p$-tree we
associate a DMGB, which
is denoted by $\wt M_{n,\ell}$, and we write $\wt d_{\mathrm{gr}}^{n,\ell}$ for
the distance
on this DMGB.

We denote the left and right boundary geodesics in $\wt M_{n,\ell}$
by $(\gamma_{n,\ell}(j),0\leq j\leq\Delta_{n,\ell})$ and $(\gamma'_{n,\ell}(j),0\leq
j\leq\Delta_{n,\ell})$,
respectively. We may now apply the results of Section~\ref{keylemma} to
$\wt M_{n,\ell}$. To this end, we
first observe that, with the exception of (b$'$), which bounds the
minimal label in $\theta_{n,\ell}$,
and (f$'$), which bounds the size of $\tau_{n,\ell}$, the properties
(a$'$)--(f$'$) do not depend
on the labeled $p$-tree $\theta_{n,\ell}$. Hence, conditionally on
$E_{n,\ell}$
and on the size $|\tau_{n,\ell}|$, the labeled $p$-tree $\theta_{n,\ell
}$ is uniformly
distributed over labeled $p$-trees with the given size, subject to the
condition that the
minimal label satisfies condition (b$'$).

Let $(\delta_n)_{n\geq0}$ be the monotone decreasing sequence
converging to $0$ constructed in Lemma~\ref{keylem2}.
To simplify notation, we set
\[
q_{n,\ell}=\biggl\lfloor(1+\lambda)^{1/4}\frac{\alpha}{3}
K^{-\ell} \frac
{n^{1/4}}{\kappa_p}\biggr\rfloor,\qquad \delta'_{n,\ell}
= \delta_{\lfloor K^{-4\ell}n\rfloor}.
\]
We define $F_{n,\ell}$ as the subset of $E_{n,\ell}$ determined by
the following
condition: For
every integer $j$ such that $q_{n,\ell} < j\leq
\alpha K^{-\ell} \frac{n^{1/4}}{\kappa_p}$, and
for every
$j'\in\{0,1,\ldots, \lfloor\alpha K^{-\ell} \frac{n^{1/4}}{\kappa_p}
\rfloor\}$,
%
%e26 #&#
\begin{eqnarray}
\label{crucialbound}&& \wt d_{\mathrm{gr}}^{n,\ell}\bigl(\ov
\gamma_{n,\ell}(j),\ov\gamma'_{n,\ell}
\bigl(j'\bigr)\bigr)
\nonumber
\\[-8pt]
\\[-8pt]
\nonumber
&&\qquad \geq \wt d_{\mathrm{gr}}^{n,\ell}
\bigl(\ov\gamma_{n,\ell}(q_{n,\ell}), \ov\gamma'_{n,\ell}
\bigl(j'\bigr)\bigr) +j-q_{n,\ell}-\delta'_{n,\ell}
n^{1/4}-2.
\end{eqnarray}

%le5.2 #&#
\begin{lemma}
\label{crucial}
We can find a constant $\ov a_1\in(0,1)$ and another constant~$C_1$,
which both depend on $\mu,A$ and $\kappa$
but not on the choice of $r\in[\mu,A]$,
such that, for every integer $\ell\geq2\ell_0$,
\[
\limsup_{n\to\infty} P \Biggl[\bigl\{\sigma_n\leq(1-\kappa)pn
\bigr\}\cap \Biggl(\bigcap_{k=\lfloor
\ell
/2\rfloor}^\ell
F_{n,k}^c \Biggr) \Biggr] \leq C_1 (\ov
a_1)^\ell.
\]
\end{lemma}

\begin{pf} Let $\ell_1,\ldots,\ell_m$ be distinct elements of $\{
\lfloor\ell/2\rfloor,\ldots,\ell\}$ for some integer
$\ell\geq2\ell_0$. We first observe that conditionally on the event
\[
\a:=\bigl\{\sigma_n\leq(1-\kappa)pn\bigr\}\cap \biggl(\bigcap
_{j\in\{\lfloor
\ell
/2\rfloor,\ldots,\ell\}\setminus\{\ell_1,\ldots,\ell_m\}} E_{n,j}^c \biggr) \cap
\Biggl(\bigcap_{i=1}^m E_{n,\ell_i}
\Biggr)
\]
and on the variables $|\tau_{n,\ell_1}|,\ldots,|\tau_{n,\ell_m}|$, the
labeled trees $\theta_{n,\ell_1},\ldots, \theta_{n,\ell_m}$
are independent and their respective conditional distributions are as
described above.
At this point, we use the fact that $K\geq\alpha_1/\alpha'_2$: Thanks
to this fact and to (\ref{conseq1}), the property (d$'$) written at order
$\ell=k$
(assuming that $E_{n,k}$ holds) puts no additional constraint on the
labeled trees $\theta_{n,k'}$ for values
$k'\not=k$ such that $E_{n,k'}$ holds.

We now use Lemma~\ref{keylem2} with $\alpha$ replaced by
$\alpha/\kappa_p$ (recall that we assumed $\alpha/\kappa_p<\alpha_0$),
$\beta_1$ replaced by $(\alpha_2+\beta_2)/\kappa_p$
and $\beta_2$ by $(1+\lambda)^{-1/4}(\alpha_1+\beta_1)/\kappa_p$. In
applying Lemma~\ref{keylem2}, we note
that the condition
\[
\biggl(\frac{\alpha_2+\beta_2}{\kappa_p} \biggr) |\tau|^{1/4} < \Delta(\tau) < \biggl(
\frac{(1+\lambda)^{-1/4}(\alpha_1+\beta_1)}{\kappa_p} \biggr) |\tau |^{1/4}
\]
together with (f$'$) implies that (b$'$) holds.
Using formula (\ref{keylemext2}) and
supposing that~$n$ is large enough so that we
can apply the estimate of Lemma~\ref{keylem2}, we get the existence of
a constant $a_1>0$ such that,\vadjust{\goodbreak}
conditionally on the event $\a$ and on the variables $|\tau_{n,\ell
_1}|,\ldots,|\tau_{n,\ell_m}|$, the property
%
%e27 #&#
\begin{eqnarray}
\label{crucialtec1} \qquad\wt d_{\mathrm{gr}}^{n,\ell_i}\bigl(\ov
\gamma'_{n,\ell_i}\bigl(j'\bigr),\ov
\gamma_{n,\ell_i}(j)\bigr) &\geq&\wt d_{\mathrm{gr}}^{n,\ell_i}\biggl(
\ov\gamma'_{n,\ell_i}\bigl(j'\bigr), \ov
\gamma_{n,\ell
_i}\biggl(\biggl\lfloor\frac{\alpha}{3\kappa_p} |
\tau_{n,\ell
_i}|^{1/4}\biggr\rfloor\biggr)\biggr)
\nonumber
\\[-8pt]
\\[-8pt]
\nonumber
&&{} + j- \biggl\lfloor\frac{\alpha}{3\kappa_p} |\tau_{n,\ell_i}|^{1/4}
\biggr\rfloor -\delta_{|\tau_{n,\ell_i}|} |\tau_{n,\ell_i}|^{1/4}-2
\nonumber
\end{eqnarray}
holds, for every $j\in\{\lfloor\frac{ \alpha}{3\kappa_p} |\tau_{n,\ell
_i}|^{1/4}\rfloor,\ldots, \lfloor\frac{\alpha}{\kappa_p} |\tau_{n,\ell
_i}|^{1/4}\rfloor\}$ and
for every
$j'\in\{0,1,\ldots, \lfloor\frac{\alpha}{\kappa_p} |\tau_{n,\ell
_i}|^{1/4}\rfloor\}$, with probability at least $a_1$,
independently for each $i=1,\ldots,m$. By (f$'$), we have
\[
K^{-4\ell_i}n < |\tau_{n,\ell_i}| < (1+\lambda) K^{-4\ell_i}n
\]
for
$i=1,\ldots,m$,
on the event $\a$. From this observation and using also the trivial bound
\[
\wt d_{\mathrm{gr}}^{n,\ell_i}\bigl(\ov\gamma'_{n,\ell_i}
\bigl(j'\bigr),\ov\gamma_{n,\ell
_i}\bigl(q'\bigr)
\bigr)-q' \geq\wt d_{\mathrm{gr}}^{n,\ell_i}\bigl(\ov
\gamma'_{n,\ell
_i}\bigl(j'\bigr),\ov
\gamma_{n,\ell_i}(q)\bigr)-q
\]
if $0\leq q'\leq q\leq\Delta_{n,\ell_i}$,
we see that, conditionally on $\a$, the event $F_{n,\ell_i}$ holds with
probability at least $a_1$, independently for each $i=1,\ldots,m$.

It follows that, for all
sufficiently large $n$,
\begin{eqnarray*}
&&P \Biggl[\bigl\{\sigma_n\leq(1-\kappa)pn\bigr\}\cap \Biggl(\bigcap
_{k=\lfloor
\ell
/2\rfloor}^\ell F_{n,k}^c
\Biggr) \Biggr]
\\
&&\qquad=\sum_{\{\ell_1,\ldots,\ell_m\}\subset\{\lfloor\ell/2\rfloor,\ldots,\ell\}} P \Biggl[\bigl\{\sigma_n
\leq(1-\kappa)pn\bigr\}
\\
&&\hspace*{99pt}\qquad\quad{} \cap \biggl(\bigcap_{j\in\{\lfloor\ell/2\rfloor,\ldots,\ell\}\setminus\{\ell_1,\ldots,\ell_m\}} E_{n,j}^c
\biggr) \\
&&\hspace*{135pt}\qquad\quad{}\cap \Biggl(\bigcap_{i=1}^m
\bigl(E_{n,\ell_i}\cap F_{n,\ell_i}^c\bigr) \Biggr) \Biggr]
\\
&&\qquad\leq\sum_{\{\ell_1,\ldots,\ell_m\}\subset\{\lfloor\ell
/2\rfloor,\ldots,\ell\}} (1-a_1)^m
P \Biggl[\bigl\{\sigma_n\leq(1-\kappa)pn\bigr\}
\\
&&\hspace*{141pt}\qquad\quad{} \cap \biggl(\bigcap_{j\in\{\lfloor\ell/2\rfloor,\ldots,\ell\}\setminus\{\ell_1,\ldots,\ell_m\}} E_{n,j}^c
\biggr) \\
&&\hspace*{215pt}\qquad\quad{}\cap \Biggl(\bigcap_{i=1}^m
E_{n,\ell_i} \Biggr) \Biggr]
\\
&&\qquad= E \bigl[\mathbf{1}_{\{\sigma_n\leq(1-\kappa)pn\}} (1-a_1)^{\sum
_{i=\lfloor\ell/2\rfloor}^\ell\mathbf{1}_{E_{n,i}}} \bigr].
\end{eqnarray*}
Hence, using (\ref{convhitting}) and (\ref{estimtech1}),
%
%e28 #&#
\begin{eqnarray}
\label{estimtech2} &&\limsup_{n\to\infty} P \Biggl[\bigl\{\sigma_n
\leq(1-\kappa)pn\bigr\}\cap \Biggl(\bigcap_{k=\lfloor
\ell
/2\rfloor}^\ell
F_{n,k}^c \Biggr) \Biggr]
\nonumber
\\[-8pt]
\\[-8pt]
\nonumber
& &\qquad\leq E \bigl[\mathbf{1}_{\{S_r\leq1-\kappa\}}(1-a_1)^{\sum
_{i=\lfloor\ell/2\rfloor}^\ell\mathbf{1}_{E_{i}}} \bigr]
\nonumber
\end{eqnarray}
and the desired result follows from Lemma~\ref{snakelemma}.
\end{pf}

We can now state and prove our main estimate. Set
$S'_r=\sup\{s\geq0\dvtx\break
Z_s=-r\}$ and, for every $\delta\in(0,r]$,
\[
\eta_\delta'(r)=\sup\Bigl\{s<S'_r\dvtx \be_s=\min_{t\in[s,S'_r]} \be_t\mbox{ and }
Z_s=-r+\delta\Bigr\}.
\]
We let ${\cal T}_\be(\eta_\delta(r))$,
respectively ${\cal T}_\be(\eta'_\delta(r))$, be the subtree of descendants of
$p_\be(\eta_\delta(r))$, respectively of $p_\be(\eta'_\delta(r))$, in
${\cal T}_\be
$. We also let $L(s_*)$
denote the ancestral line of $p_\be(s_*)$ in ${\cal T}_{\be}$. We then
consider the event
$\h_{r,\mu,\kappa}$ where the following properties hold:
\begin{longlist}[(iii)]
\item[(i)] $S_r<\infty$ and $\kappa\vee\eta_\mu(r) < s_*
<(1-\kappa
)\wedge\eta'_\mu(r) $;
\item[(ii)] $ {\inf_{a\in{\cal T}_\be(\eta_\mu(r)),b\in L(s_*)}
D(\Pi(a),\Pi(b)) > \sup_{a\in{\cal T}_\be(\eta_\mu(r))} D(\Pi
(a),\bp(S_r))} $;
\item[(iii)] $ { \inf_{a\in{\cal T}_\be(\eta'_\mu(r)),b\in
L(s_*)} D(\Pi(a),\Pi(b)) > \sup_{a\in{\cal T}_\be(\eta'_\mu(r))}
D(\Pi(a),\bp
(S_r))}.$
\end{longlist}

We will see later that if $\mu$ and $\kappa$ are chosen sufficiently
small, the probability
of the complement of $\h_{r,\mu,\kappa}$ in $\{S_r<\infty\}$ can be
made arbitrarily small.

%le5.3 #&#
\begin{lemma}
\label{mainest}
We can find a constant $\beta\in(0,1)$ and another constant $C$, which
both depend on $A$, $\mu$ and $\kappa$,
but not on the choice of $r\in[\mu,A]$,
such that, for every $\ve\in(0,1)$,
%
%e29 #&#
\begin{eqnarray}
\label{mainestim} &&P\bigl[\h_{r,\mu,\kappa} \cap\{S_{r+\ve}<\infty\}
\cap\bigl\{\exists z\in \bp \bigl(\bigl[\eta_\mu(r),
\eta'_\mu(r)\bigr]\bigr)\dvtx
\nonumber
\\
&&\hspace*{100pt}\qquad D\bigl(z,\Gamma(r)\bigr) < D\bigl(z,\Gamma(r+\ve)\bigr)+\ve\mbox{ and}\\
&&\hspace*{120pt}\qquad D\bigl(z,\Gamma(r)\bigr) < D\bigl(z,\Gamma(r-\ve)\bigr)+\ve\bigr\}\bigr] \leq C
\ve^\beta.
\nonumber
\end{eqnarray}
\end{lemma}

\begin{pf}%[Proof of Lemma~\ref{mainest}]
Clearly, we may assume that $\ve$ is small enough so that $\sqrt{\ve
}<\mu/2$.
Consider the event
\[
\a_0^\ve=\{S_{r+\ve}<\infty\}\cap \Bigl\{
\sup_{s\in[S_{r-\ve
},S_{r+\ve
}]} Z_s< -r+\sqrt{\ve}, \sup_{s\in[S'_{r+\ve},S'_{r-\ve}]}
Z_s< -r+\sqrt{\ve} \Bigr\}.
\]
By Lemma~\ref{auxil}, the probability of
$\{S_{r+\ve}<\infty\}\setminus\a_0^\ve$ is bounded above by a
constant times $\ve^{\beta_0}$.
If $\a_0^\ve$ holds, both $p_\be(S_{r-\ve})$ and $p_\be(S_{r+\ve})$
belong to ${\cal T}_\be(\eta_{\sqrt{\ve}}(r))$ and a fortiori to
${\cal T}_\be(\eta_\mu(r))$,
and a similar statement holds for $p_\be(S'_{r-\ve})$ and $p_\be
(S'_{r+\ve})$.
This implies, in particular, that $p_\be([S_{r-\ve},S_{r+\ve
}])\subset
{\cal T}_\be(\eta_{\sqrt{\ve}}(r))$
and $p_\be([S'_{r+\ve},S'_{r-\ve}])\subset{\cal T}_\be(\eta'_{\sqrt{\ve}}(r))$.\vadjust{\goodbreak}

In the first part of the proof, we assume that both $\a_0^\ve$ and the
event considered in the lemma hold. Then
there exists a point $z$ satisfying the conditions given
in (\ref{mainestim}). Let $\omega$ be a geodesic path from
$z$ to $\Gamma(r)$ in $\mm_\infty$. If $\omega$ hits the range of
$\Gamma$ before arriving at $\Gamma(r)$, then clearly
it must stay on that range [we use the fact that $\Gamma$
is the unique geodesic from $\bp(0)$ to~$\bp(s_*)$].

Next suppose that we have $\omega(t_0)\in\bp([0,S_r])$, for some
$t_0\in[0,D(z,\Gamma(r))]$.
We claim that necessarily $\omega(t_0)\in\Pi( {\cal T}_\be(\eta_{\sqrt{\ve
}}(r)))$. To see this,
write $\omega(t_0)=\bp(s_0)$ with $s_0\in[0,S_r]$, and observe that it
is enough to verify that $s_0\geq S_{r-\ve}$.
However, if $s_0< S_{r-\ve}$, by concatenating $(\omega(t),0\leq
t\leq t_0)$
with the simple geodesic $(\Gamma_{s_0}(t),0\leq t\leq r+Z_{s_0})$,
we get a geodesic from $z$ to $\Gamma(r)$ that visits $\Gamma(r-\ve)$.
This implies that
$D(z,\Gamma(r))=D(z,\Gamma(r-\ve))+\ve$, which contradicts our
assumptions on $z$.
This contradiction proves our claim and, similarly, we get that if
$\omega(t_0)\in\bp([S'_r,1])$,
then necessarily $\omega(t_0)\in\Pi({\cal T}_\be(\eta'_{\sqrt{\ve
}}(r)))$.

Note that $\Gamma(r)\notin\bp([\eta_\mu(r),\eta'_\mu(r)])$ and set
\[
T=\inf\bigl\{t\geq0\dvtx \omega(t)\notin\bp\bigl(\bigl[\eta_\mu(r),
\eta'_\mu (r)\bigr]\bigr)\bigr\}.
\]
If $\omega(T)\in\bp([0,S_r])\cup\bp([S'_r,1])$, the preceding observations
imply that $\omega(T)\in\Pi({\cal T}_\be(\eta_{\sqrt{\ve}}(r)))
\cup\Pi({\cal T}_\be(\eta'_{\sqrt{\ve}}(r)))$. If $\omega
(T)\notin\bp
([0,S_r])\cup\bp([S'_r,1])$,
then $\omega(T)\in\bp([S_r,\eta_\mu(r)]) \cup\bp([\eta'_\mu(r),S_r])$.
Since we have $\bp([S_r,\eta_\mu(r)])\subset\Pi({\cal T}_\be(\eta_\mu(r)))$ and
$\bp([\eta'_\mu(r),S_r])\subset\Pi({\cal T}_\be(\eta'_\mu
(r)))$, we get in both
cases that
\[
\omega(T)\in\Pi\bigl({\cal T}_\be\bigl(\eta_\mu(r)\bigr)
\bigr)\cup\Pi\bigl({\cal T}_\be \bigl(\eta'_\mu(r)
\bigr)\bigr).
\]

From the
fact that $z$ satisfies the conditions in (\ref{mainestim}), we
immediately get that
\begin{eqnarray*}
D\bigl(\omega(T),\Gamma(r)\bigr)&<&D\bigl(\omega(T),\Gamma(r+\ve)\bigr)+\ve,
\\
D\bigl(\omega(T),\Gamma(r)\bigr)&<&D\bigl(\omega(T),\Gamma(r-\ve)\bigr)+\ve.
\end{eqnarray*}
Replacing $z$ by $\omega(T)$, we see that the event considered in the
lemma is a.s. contained in the union of
$\{S_{r+\ve}<\infty\}\setminus\a_0^\ve$ and of the events
$\a_0^\ve\cap\a_1^\ve$ and $\a_0^\ve\cap\a_2^\ve$,  where
\begin{eqnarray*}
\a_1^\ve&:=& \h_{r,\mu,\kappa} \cap\{S_{r+\ve}<
\infty\} \\
&&{}\cap\bigl\{ \exists z\in\bp\bigl(\bigl[\eta_\mu(r),
\eta'_\mu(r)\bigr]\bigr)\cap\Pi\bigl({\cal
T}_\be\bigl(\eta_\mu(r)\bigr)\bigr)\dvtx \\
& &\hspace*{16pt}{}D\bigl(z,\Gamma(r)\bigr)<D\bigl(z,\Gamma(r+\ve)\bigr)+\ve\mbox{ and}\\
&&\hspace*{46pt}{} D
\bigl(z,\Gamma (r)\bigr)<D\bigl(z,\Gamma(r-\ve)\bigr)+\ve\bigr\}
\end{eqnarray*}
and $\a_2^\ve$ is the analogous event with $ \Pi({\cal T}_\be(\eta_\mu(r)))$
replaced by $\Pi({\cal T}_\be(\eta'_\mu(r)))$.

Let $\ell=\ell(\ve)$ be such that $K^{-\ell-1}<\sqrt{\ve}\leq
K^{-\ell}$.
We assume that $\ve$ is small enough so that $\ell(\ve)> 2\ell_0$.
We claim that
%
%e30 #&#
\begin{equation}
\label{mainestimtech} \bigl(\a_0^\ve\cap
\a_1^\ve\bigr) \subset\liminf_{n\to\infty} \Biggl(
\bigl\{\sigma_n\leq(1-\kappa)pn\bigr\}\cap \Biggl(\bigcap
_{k=\lfloor\ell
/2\rfloor
}^\ell F_{n,k}^c \Biggr)
\Biggr)\qquad \mbox{ a.s.}
\end{equation}
By Lemma~\ref{crucial}, the probability of the set in the right-hand
side is bounded
above by $\ov C_1 \ov a_1^{\ell(\ve)}$ with a constant $\ov a_1<1$.
This gives the desired estimate
for the probability of $\a_0^\ve\cap\a_1^\ve$. An analogous argument
gives a
similar estimate for the probability of $\a_0^\ve\cap\a_2^\ve$. Since
we have already
obtained the desired bound for the probability of $\{S_{r+\ve}<\infty
\}
\setminus\a_0^\ve$, the proof
of Lemma~\ref{mainest} will be complete.

To establish (\ref{mainestimtech}), we observe that, since $\a_1^\ve
\subset\{S_{r+\ve}\leq s_*<1-\kappa\}$, we have
%
%e31 #&#
\begin{equation}
\label{main11} \a_1^\ve\subset\liminf_{n\to\infty}
\bigl\{\sigma_n\leq(1-\kappa)pn\bigr\} \qquad\mbox{a.s.}
\end{equation}
Consider the event
\[
\b^\ve= \limsup_{n\to\infty} \Biggl(\bigl\{\sigma_n
\leq(1-\kappa) pn\bigr\} \cap \Biggl(\bigcup_{k=\lfloor\ell/2\rfloor}^\ell
F_{n,k} \Biggr) \Biggr).
\]
We will prove that
%
%e32 #&#
\begin{equation}
\label{main12} \bigl(\a_0^\ve\cap\a_1^\ve
\bigr)\subset\bigl(\b^\ve\bigr)^c\qquad \mbox{a.s.}
\end{equation}
Our claim
(\ref{mainestimtech}) follows from (\ref{main11}) and (\ref{main12}).

It remains to prove (\ref{main12}). To this end, we assume that both
$\a_0^\ve\cap\a_1^\ve$
and $\b^\ve$ hold, and we will see that this leads to a contradiction (except
maybe on a set of probability zero). We first choose
$z\in\bp([\eta_\mu(r),\eta'_\mu(r)])\cap\Pi({\cal T}_\be(\eta_\mu(r)))$
such that the property
stated in the definition of $\a^\ve_1$ holds, and we let $\omega$ be a
geodesic from
$z$ to $\Gamma(r)$. Note that
$\Pi({\cal T}_\be(\eta_\mu(r)))$ is contained in $\bp([0,s_*])$ by
condition
(i) in the definition of
$\h_{r,\mu,\kappa}$.
Then, from condition (ii) in the definition of $\h_{r,\mu,\kappa}$
and the fact that $z\in\Pi({\cal T}_\be(\eta_\mu(r)))$, we get
that $\omega
$ does not visit
$\Pi(L(s_*))$. Next we observe that the boundary of $\bp([0,s_*])$ is
the union of the range of $\Gamma$
and the set $\Pi(L(s_*))$, and we already noticed that if $\omega$ hits
the range
of $\Gamma$, it then stays on this range. From these observations, we
get that
$\omega$ stays in the set
$\bp([0,s_*])\setminus\Pi(L(s_*))$. Moreover, as noted at the
beginning of the proof, we
know that $\omega$ does not visit $\bp([0,S_r])$ strictly before
entering $\Pi({\cal T}_\be(\eta_{\sqrt{\ve}}(r)))$.
We choose $s_1\in(\eta_\mu(r),s_*)$ such that $z=\bp(s_1)$. This choice
is possible
since we know that $z\in\bp([\eta_\mu(r),\eta'_\mu(r)])$ and the cases
$z=\bp(\eta_\mu(r))$ and $z\in\bp([s_*,\eta'_\mu(r)])$ are
excluded by
the preceding
discussion.

We will now argue that similar properties hold for the approximating
discrete models.
Recall the notation introduced after the statement of Lemma~\ref{snakelemma}.
In particular, $y_{n,0},y_{n,1},\ldots$ are the (white) vertices of the
ancestral line
of $\ov v_n=v^n_{\sigma_n}$, and $y_{n,\psi_{n,r}(\sqrt{\ve})}$
is the last vertex on this ancestral line with label strictly larger than
$-r_n+\sqrt{\ve}\frac{n^{1/4}}{\kappa_p}$.
We denote the subtree of
descendants of the vertex $y_{n,\psi_{n,r}(\sqrt{\ve})}$ by $\tau_{n,(\sqrt{\ve})}$.
We also set
\[
r_{n,\ve}= r_n + \biggl\lfloor\ve\frac{n^{1/4}}{\kappa_p}\biggr
\rfloor,\qquad  r'_{n,\ve}= r_n - \biggl\lfloor\ve
\frac{n^{1/4}}{\kappa_p}\biggr\rfloor.
\]
Our assumption that $\a^\ve_0$ holds
ensures that,
for $n$ large enough, $\gamma_n(r_{n,\ve})$ and $\gamma_n(r'_{n,\ve})$
both belong to the subtree
$\tau_{n,(\sqrt{\ve})}$.

We also let $i_{n,*}:= \min\{i\dvtx \ell^n_i=\min_{0\leq j\leq pn} \ell^n_j\}$. Notice that
$i_{n,*}/pn$ converges to $s_*$ as $n\to\infty$, a.s.
Recall our notation $\Psi_{n,r}(\mu)$ for the index corresponding to
the last visit
of the vertex $y_{n,{\psi_{n,r}(\mu)}}$ by the contour sequence of~$\tau_n^\circ$.
Then the convergence (\ref{basicconv}) entails that
\[
\lim_{n\to\infty} \frac{\Psi_{n,r}(\mu)}{pn}= \eta_\mu(r).
\]
We choose a sequence $(j_n)$ of integers, with $j_n\in\{0,1,\ldots,pn-1\}$, such that
$\Psi_{n,r}(\mu)<j_n<i_{n,*}$ and
$\frac{j_n}{pn}$ converges
to $s_1$ as $n\to\infty$. We set $z_n=v^n_{j_n}$.

Consider then, for every $n$, a geodesic $\omega_n$ from $z_n$ to
$v^n_{\sigma_n}$ in $M_n$,
and recall that we have $\ov v_n=\gamma_n(r_n)$, where
$\gamma_n$ is the simple geodesic from (the first corner of)
$\varnothing$ to $\partial$.
We may and will assume that if $\omega_n$ hits the range of the simple
geodesic $\gamma_n$ it
stays on that range. We first observe that, for $n$ large enough,
$\omega_n$
must stay in the set $\{v^n_i\dvtx 0 \leq i< i_{n,*}\}$. Indeed, the path
$\omega_n$ starts from
a point belonging to this set and can exit it only if it visits the
range of the simple geodesic $\gamma_n$
(but in that case $\omega_n$ will stay on this range as already
mentioned) or if it visits
the ancestral line of the minimizing vertex $v^n_{i_{n,*}}$. The latter
case is also excluded
since if it holds for infinitely many values of $n$, it follows
by an easy compactness argument that there is a point $y\in\Pi
(L(s_*))$ such that
$D(z,\Gamma(r))=D(z,y)+ D(y,\Gamma(r))$, which contradicts the fact that
geodesics from $z$ to $\Gamma(r)$ do not visit $\Pi(L(s_*))$.

Let $H_n^\ve$ denote the first hitting time of $\tau_{n,(\sqrt{\ve})}$
by the path $\omega_n$.
We then claim that, again if $n$ is large enough, $\omega_n$ does not
visit $\{v^n_i\dvtx 0\leq i \leq\sigma_n\}$ strictly before~$H_n^\ve$.
Indeed, if this occurs for infinitely many values of $n$,
a discrete version of the arguments of the beginning of the proof
(using simple geodesics) shows that
$d_{\mathrm{gr}}(z_n,\gamma_n(r_n))= d_{\mathrm{gr}}(z_n,\gamma_n(r'_{n,\ve}))+
d_{\mathrm{gr}}(\gamma_n(r'_{n,\ve}),\gamma_n(r_n))$
for these values of $n$,
and a passage to the limit $n\to\infty$ gives $D(z,\Gamma
(r))=D(z,\Gamma
(r-\ve))+\ve$, contradicting
our assumption on $z$.
It follows that, for all large enough $n$,
%
%e33 #&#
\begin{equation}
\label{inclusion} \bigl\{\omega_n(j)\dvtx 0\leq j\leq
H_n^\ve\bigr\}\subset\bigl\{v^n_i\dvtx \sigma_n\leq i < i_{*,n}\bigr\}.
\end{equation}
For $j<H_n^\ve$, this is obvious from the preceding remark, and for
$j=H_n^\ve$,
we just note that a point of $\{v^n_i\dvtx 0\leq i\leq\sigma_n\}\setminus
\{
v^n_i\dvtx \sigma_n\leq i<i_{n,*}\}$
can be connected to a point of $\{v^n_i\dvtx \sigma_n\leq i<i_{n,*}\}$ only
if the latter belongs to the
ancestral line of $v^n_{\sigma_n}$, which is again excluded by the same
argument as above.

We now choose a sufficiently large value of $n$, such that (\ref
{inclusion}) holds and the event $F_{n,k}$ holds for
some $k\in\{\lfloor\ell/2\rfloor,\ldots,\ell\}$ (recall that we assume
that $\b^\ve$ holds). Recall that the definition of
$F_{n,k}$ (or rather of $E_{n,k}\supset F_{n,k}$) involves a subtree
$\tau_{n,k}$ branching from the right side of the ancestral line
of $\ov v_n$, and that $[r_{n,k},r'_{n,k}]$ is the interval
corresponding to visits of
vertices of $\tau_{n,k}$ in the contour sequence of~$\tau^\circ_n$.
Also recall properties (a$'$)--(f$'$) listed after the statement of Lemma
\ref{snakelemma}.
Note that since $K^{-\ell_0+2}<\mu$, property (a$'$) implies
$j_n> \Psi_{n,r}(\mu)\geq\Psi_{n,r}(K^{-k+2})\geq r'_{n,k}$.
On the other hand, property (a$'$) and the fact that $\sqrt{\ve
}<K^{-\ell}$
ensure that $r_{n,k}> \Psi_{n,r}(\sqrt{\ve})$.

We then set $T'_{n,k}= 1+ \max\{j\dvtx \omega_n(j)\in\{v^n_i\dvtx i>r'_{n,k}\}
\}\leq H_n^\ve$.
We observe that $\omega_n(T'_{n,k})$ must belong to the subtree $\tau_{n,k}$. Indeed,
from properties (b$'$), (c$'$) and (d$'$) we see that the minimal label on
$\tau_{n,k}$
is strictly smaller than the minimal label in $\{\sigma_n,\ldots,r_{n,k}\}$ and, thus,
a vertex of $\{v^n_i\dvtx \sigma_n\leq i < r_{n,k}\}$ cannot be connected by
an edge to
a vertex of $\{v^n_i\dvtx r'_{n,k}< i< i_{n,*}\}$. The fact that $\omega_n(T'_{n,k})\in\tau_{n,k}$ implies
that $T'_{n,k}<H_n^\ve$.

We also set $T_{n,k}= \min\{j>T'_{n,k}\dvtx \omega_n(j)\in\{v^n_i\dvtx \sigma_n\leq i< r_{n,k}\}\}\leq H_n^\ve$.
Informally, we may say that
$\omega_n(T'_{n,k})$ is an entrance point ``from the right'' for the
tree $\tau_{n,k}$ and
$\omega_n(T_{n,k}-1)$ is an exit point ``from the left'' for this tree.
More precisely,
in the DMGB $\wt M_{n,k}$ associated with $\tau_{n,k} $, the vertex
corresponding to $\omega_n(T'_{n,k})$
is connected by an edge to a vertex in the range of the right boundary
geodesic, namely,
to the point $\ov\gamma_{n,k}'(A'_n)$, with
\[
A'_n=\ell^n_{\omega_n(T'_{n,k})} -\min\bigl\{
\ell^n_v\dvtx v\in\tau_{n,k}\bigr\},
\]
and $\omega_n(T_{n,k}-1)$ corresponds
to a point of the left boundary geodesic, namely, to the point $\ov
\gamma_{n,k}(A_n)$,
with
\[
A_n= \ell^n_{\omega_n(T_{n,k}-1)} -\min\bigl\{
\ell^n_v\dvtx v\in\tau_{n,k}\bigr\}+1.
\]
See Figure~\ref{fig5} for an illustration of the preceding definitions.

%f5 #&#
\begin{figure}

\includegraphics{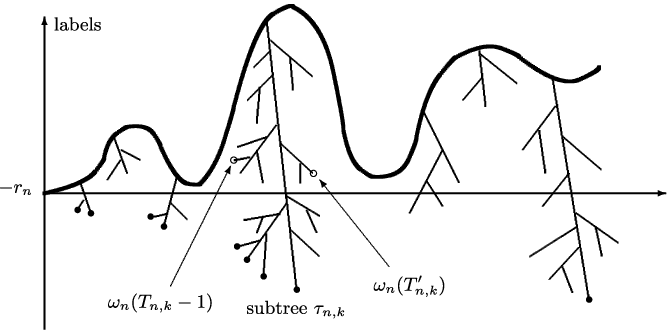}

\caption{Illustration of the proof. The thick curve represents the
evolution of labels along the
ancestral line of $\gamma_n(r_n)$ (going backward to the root). The
subtrees are those branching from
the right side of this ancestral line. The small black disks correspond
to points of the simple geodesic $\gamma_n$.
The traversal lemma makes it possible to force the geodesic $\omega_n$
to visit a point of the range of $\gamma_n$
between times $T'_{n,k}$ and $T_{n,k}-1$, without increasing too much
its length.}\label{fig5}\vspace*{-3pt}
\end{figure}

We next observe that
%
%e34 #&#
\begin{equation}
\label{estimtech3} A'_n \leq(2\alpha_1 +
\wt\alpha)K^{-k} \frac{n^{1/4}}{\kappa_p}.
\end{equation}
To see this, we use condition (e$'$) to select an index $j_0$ such that
\[
\psi_{n,r}\bigl(K^{-k+2}\bigr)\leq j_0 \leq
\tfrac{1}{2} |z_{n,k}|-1
\]
and $\ell^n_{y_{n,j_0}}
<-r_n + \wt\alpha K^{-k} \frac{n^{1/4}}{\kappa_p}$.
Notice that $\omega(T'_{n,k})$ belongs to the subtree of descendants of
$y_{n,j_0}$,
but $\omega_n(0)=z_n$ does not belong to this subtree [because
$j_n>\Psi_{n,r}(\mu)$
and $K^{-k+2}\leq K^{-\ell_0+2}<\mu$]. Then, from the construction
of edges in the BDG bijection, we get that the first vertex on the path
$\omega_n$ that belongs to the latter subtree must have a label
smaller than
or equal to the label of $y_{n,j_0}$. If $w_n$ denotes this vertex, we
have thus\vspace*{-1pt}
%
%e35 #&#
\begin{equation}
\label{estimtech33} \ell^n_{w_n}\leq\ell^n_{y_{n,j_0}}
<-r_n + \wt\alpha K^{-k} \frac{n^{1/4}}{\kappa_p}.\vspace*{-1pt}
\end{equation}
Then, using the bound (\ref{bounddgr}), and
(d$'$) and (\ref{conseq1}) to bound the minimal label between $\sigma_n$
and $\Psi_{n,r}(K^{-k+2})$,
we have\vspace*{-1pt}
\[
d_{\mathrm{gr}}(\ov v_n,w_n)\leq\ell^n_{w_n}
+ \ell^n_{\ov v_n} -2 \min_{\sigma
_n\leq i\leq\Psi_{n,r}(K^{-k+2})}\ell^n_{v^n_i}+2
\leq(2\alpha_1 + \wt\alpha)K^{-k} \frac{n^{1/4}}{\kappa_p}.\vspace*{-1pt}
\]
On the other hand, since the points $w_n$, $\omega_n(T'_{n,k})$ and
$\ov v_n$ come in that order on
the geodesic $\omega_n$, we have also
\begin{eqnarray*}
d_{\mathrm{gr}}(\ov v_n,w_n) &=& d_{\mathrm{gr}}
\bigl(\ov v_n,\omega_n\bigl(T'_{n,k}
\bigr)\bigr) + d_{\mathrm{gr}}\bigl(\omega_n\bigl(T'_{n,k}
\bigr), w_n\bigr)
\\[-1pt]
&\geq&\bigl(\ell^n_{\omega(T'_{n,k})}-\ell^n_{\ov v_n}
\bigr) + \bigl(\ell^n_{\omega
(T'_{n,k})}-\ell^n_{w_n}
\bigr)
\\[-1pt]
&\geq&2A'_n + 2\biggl(-r_n -
\alpha_1K^{-k}\frac{n^{1/4}}{\kappa_p} +1\biggr) + r_n
- \biggl(-r_n+ \wt\alpha K^{-k}\frac{n^{1/4}}{\kappa_p}\biggr)
\\[-1pt]
&=& 2A'_n - (2\alpha_1 + \wt
\alpha)K^{-k} \frac{n^{1/4}}{\kappa_p}+2.
\end{eqnarray*}
In the second inequality, we used (\ref{estimtech33}) and (\ref
{conseq1}), which gives
a lower bound on the minimal label in $\tau_{n,k}$.
Our claim (\ref{estimtech3}) follows by combining the last two
displays.

The same argument shows that (\ref{estimtech3}) still holds if we
replace $A'_n$
by $A_n-1$ (in fact, the same bound holds for the difference between
the label of any
vertex of $\tau_{n,k}$ that is visited by the path $\omega_n$ and the
minimal label
on $\tau_{n,k}$).

The crucial observation now is the lower bound
\[
T_{n,k}-T'_{n,k} \geq\wt d^{n,k}_{\mathrm{gr}}
\bigl(\ov\gamma'_{n,k}\bigl(A'_n
\bigr),\ov \gamma_{n,k}(A_n)\bigr).
\]
This is clear since between $T'_{n,k}$ and $T_{n,k}-1$ the path $\omega_n$ stays in the tree $\tau_{n,k}$ and
uses only edges that are present in the associated DMGB $M_{n,k}$.
Recalling that $2\alpha_1+\wt\alpha< \alpha$, we can then use the bound
(\ref{crucialbound}) to estimate $\wt d^{n,k}_{\mathrm{gr}} (\ov\gamma'_{n,k}(A'_n),\ov\gamma_{n,k}(A_n))$.

Suppose first that $A_n \geq(\alpha_2-\alpha'_2) K^{-k}\frac
{n^{1/4}}{\kappa_p}$ (the other case, which is
simpler, will be considered next). We use the fact that $F_{n,k}$
holds. Since $\alpha_2-\alpha'_2 > (1+\lambda)^{1/4}\frac{\alpha
}{3}$, we
can apply (\ref{crucialbound}) with $j=A_n$ and $j'=A'_n$,
and we get
\begin{eqnarray*}
&&\wt d_{\mathrm{gr}}^{n,k}\bigl(\ov\gamma_{n,k}(A_n),
\ov\gamma'_{n,k}\bigl(A'_n\bigr)
\bigr)
\\
&&\qquad \geq\wt d_{\mathrm{gr}}^{n,k}\bigl(\ov\gamma_{n,k}(q_{n,k}),
\ov\gamma'_{n,k}\bigl(A'_n\bigr)
\bigr) +A_n-q_{n,k} -\delta'_{n,k}
n^{1/4}.
\end{eqnarray*}
It follows that we can modify the part of the path $\omega_n$ between times
$T'_{n,k}$ and $T_{n,k}-1$ in such a way that the new path goes through
the vertex
$v_{(n)}=\ov\gamma_{n,k}(q_{n,k})$, and the length of $\omega_n$ is
increased by at most $\delta'_{n,k} n^{1/4}
\leq\delta'_{n,\ell} n^{1/4}$. Write $\omega'_n$
for the new path obtained after this modification. Now notice that
\begin{eqnarray*}
\ell^n_{v_{(n)}}&=& \Bigl(\ell^n_{v_{(n)}} -
\min_{v\in\tau_{n,k}} \ell^n_v \Bigr) +
\min_{v\in\tau_{n,k}} \ell^n_v
\\
&\leq& q_{n,k} -r_n-\alpha_2 K^{-k}
\frac{n^{1/4}}{\kappa_p}
\\
&\leq&-r_n - \alpha'_2K^{-k}
\frac{n^{1/4}}{\kappa_p}
\end{eqnarray*}
using (\ref{conseq1}) in the first inequality, and then the bound
$\alpha_2-\alpha'_2 > (1+\lambda)^{1/4}\frac{\alpha}{3}$.
However, by property (d$'$), the right-hand side of the last display is
strictly smaller than $\min\{\ell^n_i\dvtx i\in[\sigma_n,r_{n,k}]\}$. This
implies that the simple geodesic $\omega_{(\sigma_n)}$ starting from
$\ov v_n=v^n_{\sigma_n}$ will visit the vertex $v_{(n)}$. Hence, we can
modify the path
$\omega'_n$ without increasing its length, in such a way that it
coalesces with the (time-reversed) simple geodesic
$\gamma_n$ before entering the subtree $\tau_{n,(\sqrt{\ve})}$. In
particular, the modified path $\omega'_n$
will visit the vertex $\gamma_n(r_{n,\ve})$, which belongs to the
subtree $\tau_{n,(\sqrt{\ve})}$, and we have obtained
%
%e36 #&#
\begin{equation}
\label{lastbound} d_{\mathrm{gr}}(z_n, \ov v_n) \geq
d_{\mathrm{gr}}\bigl(z_n,\gamma_n(r_{n,\ve})
\bigr) + (r_{n,\ve}- r_n) - \delta'_{n,\ell}
n^{1/4}.
\end{equation}

If $A_n \leq(\alpha_2-\alpha'_2) K^{-k}\frac{n^{1/4}}{\kappa_p}$,\vspace*{-1pt} we
get the same bound (\ref{lastbound})
without the term $\delta'_{n,\ell} n^{1/4} $ in a much simpler way,
since the same arguments as above directly show that the simple
geodesic $\omega_{(\sigma_n)}$ starting from
$\ov v_n$ visits the point $\ov\gamma_{n,k}(A_n)=\omega_n(T_{n,k}-1)$
after visiting $\gamma_n(r_{n,\ve})$.

Finally, the lower bound (\ref{lastbound}) holds for any (sufficiently
large) $n$ such that
$F_{n,k}$ holds for some $k\in\lfloor\ell/2,\ell\rfloor$. We are
assuming that there
are infinitely many such values of $n$, and so we can pass to
the limit $n\to\infty$ in (\ref{lastbound}) after multiplying by
$\kappa_pn^{-1/4}$ to get
\[
D\bigl(z,\Gamma(r)\bigr)\geq D\bigl(z,\Gamma(r+\ve)\bigr) + \ve.
\]
This contradicts the fact that $z$ satisfies the property given in the
definition of $\a^\ve_1$.
This contradiction completes the
proof of (\ref{main12}) and of Lemma~\ref{mainest}.
\end{pf}

%s6 #&#
\section{A preliminary bound on distances}\label{sec6}

The proof of our main theorem uses a preliminary estimate, which we
state in the following
proposition.

%pr6.1 #&#
\begin{proposition}
\label{estimate-distance}
Let $\delta\in(0,1)$. There exists a (random) constant $C_\delta$
such that,
for every $x,y\in\mm_\infty$,
\[
D^*(x,y)\leq C_\delta D(x,y)^{1-\delta}.
\]
\end{proposition}

\begin{pf} We write $B_D(x,h)$, respectively $B_{D^*}(x,h)$, for the open ball
of radius~$h$
centered at $x$ in $(\mm_\infty,D)$, respectively in $(\mm_\infty,D^*)$. As
usual, the corresponding closed
balls are denoted by $\ov B_D(x,h)$ and $\ov B_{D^*}(x,h)$.
Recall from Section~\ref{convBrma} the definition of the
volume measure $\operatorname{Vol}$ on $\mm_\infty$.
From Corollary 6.2 in~\cite{AM},
there exists a (random) constant $c_\delta$ such that, for every $h\in(0,1)$,
%
%e37 #&#
\begin{equation}
\label{ballsD} \sup_{x\in\mm_\infty} \operatorname{Vol}\bigl(\ov B_D(x,h)
\bigr) \leq c_\delta h^{4-\delta}.
\end{equation}
On the other hand, it is also easy to verify that, for every $h\in(0,1)$,
%
%e38 #&#
\begin{equation}
\label{ballsDstar} \inf_{x\in\mm_\infty}\operatorname{Vol}\bigl( B_{D^*}(x,h)
\bigr) \geq c'_\delta h^{4+\delta}
\end{equation}
for some other (random) constant $c'_\delta>0$. To obtain this
estimate, just
use the bound $D^*(a,b)\leq D^\circ(a,b)$ for $a,b\in{\cal T}_\be$,
and the
fact that
the process $(Z_t)_{0\leq t\leq1}$ is H\"older continuous with
exponent $\frac{1}{4}-\ve$, for every $\ve>0$.

Let $x,y\in\mm_\infty$, and let $\omega=(\omega(t),0\leq t\leq D(x,y))$
be a geodesic
from $x$ to $y$ with respect to the metric $D$. To get the bound of the
proposition,
we may assume that $0<D(x,y)<1/2$. Put $t_0=0$ and set
\[
t_1=\sup\bigl\{t\geq0\dvtx \omega(t)\in\ov B_{D^*}
\bigl(x,D(x,y)\bigr)\bigr\}.
\]
If $t_1=D(x,y)$, we stop the construction. Otherwise we proceed by induction.
For every integer $n\geq1$ such that $t_n$ has been defined and
$t_n<D(x,y)$, we set
\[
t_{n+1}=\sup\bigl\{t\geq t_n\dvtx \omega(t)\in\ov
B_{D^*}\bigl(\omega (t_n),D(x,y)\bigr)\bigr\}.
\]
A simple argument, using the fact that the topologies induced by $D$ and
$D^*$ coincide, shows that the construction stops after a finite number
$n_{\mathrm{max}}$ of steps, such that $t_{n_{\mathrm{max}}}=D(x,y)$. The key point now
is to observe that
the balls $B_{D^*}(\omega(t_i),\frac{1}{2}D(x,y))$ and
$B_{D^*}(\omega
(t_j),\frac{1}{2}D(x,y))$
are disjoint if $0\leq i<j < n_{\mathrm{max}}$. Indeed, if this is not the case,
we get $\omega(t_j)\in B_{D^*}(\omega(t_i),D(x,y))$
and thus $t_j<t_{i+1}$, which is absurd. Using (\ref{ballsDstar})
and the bound $D\leq D^*$, it follows that
\[
n_{\mathrm{max}}\times c'_\delta\biggl(
\frac{D(x,y)}{2}\biggr)^{4+\delta} \leq\operatorname{Vol}\bigl(\ov B_D
\bigl(x,2D(x,y)\bigr)\bigr).
\]
On the other hand, (\ref{ballsD}) gives
\[
\operatorname{Vol}\bigl(\ov B_D\bigl(x,2D(x,y)\bigr)\bigr) \leq
c_\delta2^{4-\delta} D(x,y)^{4-\delta}.
\]
By combining the last two bounds, we get
$n_{\mathrm{max}}\leq\frac{c_\delta}{c'_\delta} 2^8 D(x,y)^{-2\delta}.$
Since
\[
D^*(x,y)\leq D^*\bigl(\omega(0),\omega(t_1)\bigr)+\cdots+ D^*\bigl(
\omega (t_{n_{\mathrm{max}}-1}),\omega(t_{n_{\mathrm{max}}})\bigr)\leq n_{\mathrm{max}}
D(x,y),
\]
the proof of the proposition is complete.
\end{pf}

%s7 #&#
\section{Proof of the main result}\label{sec7}

In this section we suppose that $z$ is a random point of $\mm_\infty$
distributed according to the uniform measure $\operatorname{Vol}$. We may define
$z=\mathbf{p}(U)$ where $U$ is uniformly distributed over $[0,1]$ and independent
of all other random quantities. Recall the constant $\beta$
from Lemma~\ref{mainest}.

%le7.1 #&#
\begin{lemma}
\label{badintervals}
Let $u>0$ and $A>u$, and, for every integer $k\geq1$, let $\h_k(z)$ be
the collection of
all integers $i$ with $\lfloor2^k u\rfloor< i < \lfloor2^k(\Delta
\wedge A)\rfloor$, such that
we have both
\[
D\bigl(z,\Gamma\bigl(i2^{-k}\bigr)\bigr) < D\bigl(z, \Gamma
\bigl((i+1)2^{-k}\bigr)\bigr) + 2^{-k}
\]
and
\[
D\bigl(z,\Gamma\bigl(i2^{-k}\bigr)\bigr) < D\bigl(z, \Gamma
\bigl((i-1)2^{-k}\bigr)\bigr) + 2^{-k}.
\]
Then, for every $\beta'\in(0,\beta)$,
\[
2^{-(1-\beta')k} \#\h_k(z)\build{\la}_{k\to\infty}^{\mathrm{a.s.}}
0.
\]
\end{lemma}

\begin{rema*}
Since $\Gamma$ is a geodesic, it is obvious that the weak inequality
$\leq$
holds instead of $<$ in both displayed inequalities of the lemma. The
point is that for most values
of $i$ one of these two weak inequalities can be replaced by an equality.
\end{rema*}

\begin{pf} We fix a constant $\kappa\in(0,1/4)$ and $\mu\in(0,u]$.
Recall the notation $\eta_\delta(r)$ and $\eta'_\delta(r)$ introduced
in the
previous section. We consider the subset $\h'_k(z)$ of $\h_k(z)$ that
consists of
all integers $i\in\h_k(z)$
such that $z\in\bp([\eta_\mu(i2^{-k}),\eta'_\mu(i2^{-k})])$,
\[
\kappa\vee\eta_\mu\bigl(i2^{-k}\bigr)< s_*<(1-\kappa)
\wedge\eta'_\mu\bigl(i2^{-k}\bigr)
\]
and
%
%e39 #&#
%e40 #&#
\begin{eqnarray}
\label{bad00} \qquad\inf_{a\in{\cal T}_\be(\eta_\mu(i2^{-k})),b\in L(s_*)} D\bigl(\Pi (a),\Pi(b)\bigr)&
>&
\sup_{a\in{\cal T}_\be(\eta_\mu(i2^{-k}))} D\bigl(\Pi(a),\bp (S_{i2^{-k}})\bigr),
\\
\label{bad01} \inf_{a\in{\cal T}_\be(\eta'_\mu(i2^{-k})),b\in L(s_*)} D\bigl(\Pi (a),\Pi(b)\bigr)
&>&
\sup_{a\in{\cal T}_\be(\eta'_\mu(i2^{-k}))} D\bigl(\Pi(a),\bp (S_{i2^{-k}})\bigr).
\end{eqnarray}

By Lemma~\ref{mainest}, we have for every $i$ such that $\lfloor2^k
u\rfloor< i < \lfloor2^kA\rfloor$,
\[
P\bigl[ i \in\h'_k(z)\bigr] \leq C 2^{-k\beta}.
\]
From this bound, it immediately follows that, if $0<\beta'<\beta$,
%
%e41 #&#
\begin{equation}
\label{bad1} 2^{-(1-\beta')k}\#\h'_k(z) \build{
\la}_{k\to\infty}^{\mathrm{a.s.}} 0.
\end{equation}

To complete the proof, we need to control $\#(\h_k(z)\setminus\h'_k(z))$. We first note that
the property $\kappa<s_*<1-\kappa$ holds on an event of probability
arbitrarily close
to $1$, if $\kappa$ is chosen small enough. Furthermore, on the event
$\{\kappa<s_*<1-\kappa\}$, we have
\[
\h_k(z)\setminus\h'_k(z) \subset \bigl(
\h^{(1)}_k(z)\cup\h^{(2)}_k(z)\cup
\h^{(3)}_k(z)\cup\h^{(4)}_k(z) \cup
\h^{(5)}_k(z)\cup\h^{(6)}_k(z) \bigr),
\]
where $\h^{(1)}_k(z),\ldots,\h^{(6)}_k(z)$ are the subsets of $\{
\lfloor2^k u \rfloor+1,\ldots,\lfloor2^k A\rfloor-1\}$
defined by
\begin{eqnarray*}
\h^{(1)}_k(z)&=&\bigl\{i\dvtx S_{i2^{-k}}<\infty, z\in\bp
\bigl(\bigl[S_{i2^{-k}},S'_{i2^{-k}}\bigr]\bigr)\mbox{
and } s_*\leq\eta_\mu\bigl(i2^{-k}\bigr)\bigr\},
\\
\h^{(2)}_k(z)&=&\bigl\{i\dvtx S_{i2^{-k}}<\infty, z\in\bp
\bigl(\bigl[S_{i2^{-k}},S'_{i2^{-k}}\bigr]\bigr)\mbox{
and } s_*\geq\eta'_\mu\bigl(i2^{-k}\bigr)
\bigr\},
\\
\h^{(3)}_k(z)&=&\bigl\{i\dvtx S_{i2^{-k}}<\infty \mbox{
and } z\in\bp\bigl(\bigl[S_{(i-1)2^{-k}},\eta_\mu
\bigl(i2^{-k}\bigr)\bigr]\bigr)\bigr\},
\\
\h^{(4)}_k(z)&=&\bigl\{i\dvtx S_{i2^{-k}}<\infty \mbox{
and } z\in\bp\bigl(\bigl[\eta'_\mu\bigl(i2^{-k}
\bigr),S'_{(i-1)2^{-k}}\bigr]\bigr)\bigr\},
\\
\h^{(5)}_k(z)&=&\bigl\{i\dvtx S_{i2^{-k}}<\infty, z\in\bp
\bigl(\bigl[S_{(i-1)2^{-k}},S'_{(i-1)2^{-k}}\bigr]\bigr)\mbox{
and (\ref{bad00}) fails}\bigr\},
\\
\h^{(6)}_k(z)&=&\bigl\{i\dvtx S_{i2^{-k}}<\infty, z\in\bp
\bigl(\bigl[S_{(i-1)2^{-k}},S'_{(i-1)2^{-k}}\bigr]\bigr)\mbox{
and (\ref{bad01}) fails}\bigr\}
\end{eqnarray*}
[notice that if $z\in\bp([0,S_{(i-1)2^{-k}}])$ or if $z\in\bp
([S'_{(i-1)2^{-k}},1])$, by considering a simple
geodesic starting from $z$, we get automatically $D(z,\Gamma(i2^{-k}))
= D(z, \Gamma((i-1)2^{-k})) + 2^{-k}$,
so that $i$ cannot belong to $\h_k(z)$].

Then, if $\h^{(1)}_k(z)\not=\varnothing$, we can find $r\geq u$ such
that $S_r<\infty$, $s_*\leq\eta_\mu(r)$
and $z\in\bp([S_r,S'_r])$. Hence, if we define, on the event $\{
\Delta
> u\}$,
\[
r_\mu= \inf\bigl\{r > u\dvtx S_r<\infty\mbox{ and } s_*
\leq\eta_\mu(r)\bigr\},
\]
we have
\[
\Biggl(\bigcup_{k=1}^\infty\bigl\{
\h^{(1)}_k(z)\not=\varnothing\bigr\} \Biggr) \subset \bigl\{
\Delta>u, z\in\bp\bigl(\bigl[S_{r_\mu},S'_{r_{\mu}}
\bigr]\bigr)\bigr\}.
\]
Then
%
%e42 #&#
\begin{equation}
\label{bad2} P\bigl(\Delta> u, z\in\bp\bigl(\bigl[S_{r_\mu},S'_{r_{\mu}}
\bigr]\bigr)\bigr)= E\bigl[\mathbf{1}_{\{
\Delta>
u\}} \bigl(S'_{r_\mu}-S_{r_\mu}
\bigr)\bigr].
\end{equation}

We claim that $r_\mu\la\Delta$ as $\mu\da0$, a.s. on the event
$\{\Delta> u\}$. To see this,
we observe that on the latter event we have for every $\ve\in
(0,\Delta-u)$,
%
%e43 #&#
\begin{equation}
\label{bad222} \inf_{u\leq r\leq\Delta-\ve} (r- Z_{p_\be(S_r)\wedge p_\be(s_*)}) >0.
\end{equation}
Indeed, if the infimum in (\ref{bad222}) vanishes, a compactness argument
gives $r_0\in[0,\Delta-u]$ such that either $Z_{p_\be(S_{r_0})\wedge
p_\be(s_*)}=r_0$
or $Z_{p_\be(S_{r_0+})\wedge p_\be(s_*)}=r_0$ (here $S_{r_0+}$ stands
for the
right limit of $r\to S_r$ at $r=r_0$). However, this implies that
$p_\be
(S_{r_0})=p_\be(S_{r_0})\wedge p_\be(s_*)$,
or $p_\be(S_{r_0+})=p_\be(S_{r_0+})\wedge p_\be(s_*)$, is an ancestor
of $p_\be(s_*)$ in ${\cal T}_\be$, which is impossible
since Lemma~\ref{increase} shows that all vertices of the form $p_\be
(S_r)$ or $p_\be(S_{r+})$ are leaves of ${\cal T}_\be$.

Then (\ref{bad222}) implies that $r_\mu>\Delta-\ve$ if $\mu$ is small
enough, and gives our claim. Once we
know that $r_\mu\la\Delta$ as $\mu\da0$, dominated convergence
entails that the left-hand side
of (\ref{bad2}) tends to $0$ as $\mu\to0$. So by choosing $\mu$ small
enough, we
get that all sets $\h^{(1)}_k(z)$ are empty, except on a set of
arbitrarily small probability. The same
argument applies to the sets $\h^{(2)}_k(z)$.

To deal with $\mathcal{H}^{(3)}_k(z)$, we first observe that, a.s., for each fixed $k$, there is
at most one value of $i$ such that $z\in \mathbf{p}([S_{(i-1)2^{-k}}, S_{i 2^{-k}}])$. Hence,
%
%e44 #&#
\begin{equation}
\label{bad3} \Biggl(\bigcup_{k=1}^\infty
\bigl\{\#\mathcal{H}^{(3)}_k(z)> 1\bigr\} \Biggr) \subset
\biggl\{ z \in\bp \biggl( \bigcup_{r\geq u, S_r<\infty}
\bigl[S_r,\eta_\mu (r)\bigr] \biggr) \biggr\}.
\end{equation}
Using again the fact that the vertices $p_\be(S_r)$ are leaves of
${\cal T}_\be$, one easily
verifies that the sets
\[
\bp \biggl( \bigcup_{r\geq u, S_r<\infty} \bigl[S_r,
\eta_\mu(r)\bigr] \biggr)
\]
decrease when $\mu\da0$ to the set $\{\Gamma(r)\dvtx u\leq r\leq\Delta\}$.
Since the latter
set has zero volume, we get that the probability of the event in the
right-hand side
of (\ref{bad3}) tends to $0$ as $\mu\da0$. So we can choose $\mu>0$
sufficiently
small so that all sets $\h^{(3)}_k(z)$ have cardinality at most $1$, except on a set of
arbitrarily small probability. The same
argument applies to the sets $\h^{(4)}_k(z)$.

Finally, we consider $\h^{(5)}_k(z)$. Let $\delta>0$. We observe that
\[
\bigcup_{\lfloor2^ku\rfloor<i\leq\lfloor2^k(\Delta-\delta
)\rfloor} \Pi\bigl({\cal T}_\be
\bigl(\eta_\mu\bigl(i2^{-k}\bigr)\bigr)\bigr) \subset\bigcup
_{u\leq r\leq\Delta-\delta} \Pi\bigl({\cal T}_\be\bigl(
\eta_\mu(r)\bigr)\bigr).
\]
In a way similar to the previous step of the proof, we can check that
the sets
\[
\bigcup_{u\leq r\leq\Delta-\delta} \Pi\bigl({\cal T}_\be
\bigl(\eta_\mu(r)\bigr)\bigr)
\]
are closed and decrease a.s. to $\{\Gamma(r)\dvtx u\leq r\leq\Delta
-\delta
\}$ when
$\mu\da0$. Furthermore, it is not hard to verify, again by a
compactness argument, that
\[
\sup_{u\leq r\leq\Delta} \Bigl( \sup_{x,y\in\Pi({\cal T}_\be
(\eta_\mu(r)))} D(x,y) \Bigr) \build{
\la}_{\mu\to0}^{\mathrm{a.s.}} 0.
\]
Since
\[
\inf_{u\leq r\leq\Delta-\delta, b\in L(s_*)} D\bigl(\Gamma(r),\Pi(b)\bigr) > 0\qquad \mbox{a.s.}
\]
it follows from the preceding considerations that, for any given
$\delta
>0$, we can choose $\mu>0$ sufficiently small so
that the property
\[
\sup_{u\leq r\leq\Delta} \Bigl( \sup_{x,y\in\Pi({\cal T}_\be
(\eta_\mu
(r)))} D(x,y) \Bigr) <
\inf_{u\leq r\leq\Delta-\delta} \Bigl(\inf_{a\in{\cal T}_\be
(\eta_\mu
(r)), b\in L(s_*)} D\bigl(\Pi(a),\Pi(b)\bigr) \Bigr)
\]
holds with a probability arbitrarily close to $1$. If the latter
property holds, this means that
the only indices $i$ for which (\ref{bad00}) may fail are those such
that $i2^{-k}> \Delta-\delta$.
For such indices $i$ the property $z\in\bp
([S_{(i-1)2^{-k}},S'_{(i-1)2^{-k}}])$ implies that
$z\in\bp([S_{\Delta-\delta},S'_{\Delta-\delta}])$ and if $\delta
$ has
been chosen small enough,
this also occurs with a small probability. So again
we can choose $\mu>0$ sufficiently
small so that all sets $\h^{(5)}_k(z)$ are empty, except on a set of
arbitrarily small probability.
The same
argument applies to the sets $\h^{(6)}_k(z)$.

From the preceding arguments, we can fix $\kappa$ and $\mu$
sufficiently small so
that outside a set of small probability we have $\#(\mathcal{H}_k(z)\backslash\mathcal{H}'_k(z))\leq 1$ for
every~$k$. The
conclusion of Lemma~\ref{badintervals} now follows from (\ref{bad1}).
\end{pf}

%th7.2 #&#
\begin{theorem}
\label{MainT}
We have $D(y,y')=D^*(y,y')$ for every $y,y'\in\mm_\infty$, almost surely.
\end{theorem}

As was already explained at the end of Section~\ref{convBrma},
Theorem~\ref{mainresult} (in the bipartite case) readily follows from
Theorem~\ref{MainT}.

\begin{pf} It is sufficient to verify that the identity of the theorem
holds when~$y$ and~$y'$ are independently distributed according to the volume measure
on $\mm_\infty$ [indeed, if $D(y_0,y_0')<D^*(y_0,y_0')$ for some
$y_0,y'_0\in\mm_\infty$, then
the same strict inequality holds for every $y$ and $y'$ sufficiently
close to $y_0$ and $y'_0$,
respectively, and we use the fact that the volume measure has full
support in~$\mm_\infty$].
Let $z$ be as in Lemma~\ref{badintervals}. Since the distinguished
point in
a uniformly distributed rooted and pointed $2p$-angulation is chosen
uniformly at
random among the vertices, it is easy to verify that
the random triply pointed metric spaces $(\mm_\infty, \Pi(\rho),
x_*,z)$ and
$(\mm_\infty,\Pi(\rho),y,y')$ have the same distribution. A~simple
application of Theorem~8.1
in~\cite{AM} then gives the following identity in distribution:
%
%e45 #&#
\begin{equation}
\label{triplident} \bigl(\mm_\infty, \Pi(\rho), x_*,z\bigr)
\build{=}_{}^{\mathrm{(d)}} \bigl(\mm_\infty,
y,y',x_*\bigr).
\end{equation}

Let $\wt\Gamma=(\wt\Gamma(t),0\leq t\leq D(y,y'))$ be the almost surely
unique geodesic from $y$
to $y'$ in $(\mm_\infty,D)$. The almost sure uniqueness of this geodesic
follows from Corollary 8.3(i) in~\cite{AM}. Fix $u>0$ and, for every
integer $k\geq1$, let $\h_k(y,y')$ stand for the set of
all integers $i$, with $\lfloor2^ku\rfloor\leq i < \lfloor2^k
D(y,y')\rfloor$, such that
we have both
\[
D\bigl(x_*,\wt\Gamma\bigl(i2^{-k}\bigr)\bigr) < D\bigl(x_*, \wt\Gamma
\bigl((i+1)2^{-k}\bigr)\bigr) + 2^{-k}
\]
and
\[
D\bigl(x_*,\wt\Gamma\bigl(i2^{-k}\bigr)\bigr) < D\bigl(x_*,\wt\Gamma
\bigl((i-1)2^{-k}\bigr)\bigr) + 2^{-k}.
\]
By Lemma~\ref{badintervals} and the identity in distribution (\ref
{triplident}), we have,
for every $\beta'\in(0,\beta)$,
%
%e46 #&#
\begin{equation}
\label{MainTech1} 2^{-(1-\beta')k} \#\h_k\bigl(y,y'
\bigr)\build{\la}_{k\to\infty}^{\mathrm{a.s.}} 0.
\end{equation}

Then let $\h_k^\bullet(y,y')$ stand for the set of
all integers $i$, with $\lfloor2^ku\rfloor\leq i < \lfloor2^k
D(y,y')\rfloor$, such that
\[
 \bigl|D\bigl(x_*,\wt\Gamma\bigl(i2^{-k}\bigr)\bigr) - D\bigl(x_*, \wt\Gamma
\bigl((i+1)2^{-k}\bigr)\bigr)\bigr| < 2^{-k}.
\]
If $i\notin\h_k^\bullet(y,y')$, then we have either
\[
D\bigl(x_*,\wt\Gamma\bigl(i2^{-k}\bigr)\bigr) = D\bigl(x_*,\wt\Gamma
\bigl((i+1)2^{-k}\bigr)\bigr) + 2^{-k}
\]
or
\[
D\bigl(x_*,\wt\Gamma\bigl(i2^{-k}\bigr)\bigr) = D\bigl(x_*,\wt\Gamma
\bigl((i+1)2^{-k}\bigr)\bigr) - 2^{-k}.
\]
An elementary argument shows that if $i\in\h_k^\bullet(y,y')$ and
$i'=\max\{j\leq i\dvtx j\in\h_k(y,y')\}$, with the convention $\max\varnothing=\lfloor
2^ku\rfloor$,
then, for every integer $j$ such that $i'\leq j<i$, we have
$j\notin\h^\bullet_k(y,y')$. It follows that $\#\h_k^\bullet
(y,y')\leq
\#\h_k(y,y')+1$,
and, in particular, (\ref{MainTech1}) implies also, for every $\beta'\in
(0,\beta)$,
%
%e47 #&#
\begin{equation}
\label{MainTech2} 2^{-(1-\beta')k} \#\h_k^\bullet
\bigl(y,y'\bigr)\build{\la}_{k\to\infty
}^{\mathrm{a.s.}} 0.
\end{equation}

Now suppose that $i\in\{\lfloor2^ku\rfloor,\ldots,\lfloor2^k
D(y,y')\rfloor-1\}$
is not in $\h_k^\bullet(y,y')$. Then either $\wt\Gamma(i2^{-k})$
lies on
a geodesic path from $\wt\Gamma((i+1)2^{-k})$ to $x_*$ or, conversely,
$\wt\Gamma((i+1)2^{-k})$ lies on
a geodesic path from $\wt\Gamma(i2^{-k})$ to $x_*$. By Theorem~\ref
{geodesicroot},
any of these geodesic paths is a simple geodesic,
and is also a geodesic in $(\mm_\infty,D^*)$, so that, using
(\ref{distance-root}), we have $D^*(\wt\Gamma(i2^{-k}),\wt\Gamma
((i+1)2^{-k}))=2^{-k}$.

To conclude, we write, for $u<D(z,z')$,
%
%e48 #&#
\begin{eqnarray}
\label{MainTech3} D^*\bigl(y,y'\bigr)&\leq& D^*\bigl(y,\wt\Gamma
\bigl(2^{-k}\bigl\lfloor2^ku\bigr\rfloor\bigr)\bigr)\nonumber\\
&&{} + \sum
_{i=\lfloor2^ku\rfloor}^{\lfloor2^k D(y,y')\rfloor-1} D^*\bigl(\wt\Gamma
\bigl(i2^{-k}\bigr),\wt\Gamma\bigl((i+1)2^{-k}\bigr)\bigr)
\\
&&{}+ D^*\bigl(\wt\Gamma\bigl(2^{-k}\bigl\lfloor2^k D
\bigl(y,y'\bigr)\bigr\rfloor\bigr),y'\bigr).\nonumber
\end{eqnarray}
By previous observations,
\begin{eqnarray*}
&&\sum_{i=\lfloor2^ku\rfloor}^{\lfloor2^k D(y,y')\rfloor-1} D^*\bigl(\wt\Gamma
\bigl(i2^{-k}\bigr),\wt\Gamma\bigl((i+1)2^{-k}\bigr)\bigr)
\\
&&\qquad\leq2^{-k}\bigl\lfloor2^k D\bigl(y,y'\bigr)
\bigr\rfloor
 \\
 &&\qquad\quad{}+ \bigl(\# \h_k^\bullet\bigl(y,y'\bigr)
\bigr)\sup_{0\leq i<\lfloor2^k
D(y,y')\rfloor} D^*\bigl(\wt\Gamma\bigl(i2^{-k}\bigr),
\wt\Gamma\bigl((i+1)2^{-k}\bigr)\bigr).
\end{eqnarray*}
Proposition~\ref{estimate-distance} implies that, for every $\delta
\in
(0,1)$, there
exists a (random) constant $c_\delta$ such that
\[
\sup_{0\leq i<\lfloor2^k D(y,y')\rfloor} D^*\bigl(\wt\Gamma\bigl(i2^{-k}\bigr),\wt\Gamma
\bigl((i+1)2^{-k}\bigr)\bigr) \leq c_\delta 2^{-k(1-\delta)}.
\]
Applying this bound with $\delta<\beta$ and using (\ref{MainTech2}),
we get
\[
\bigl(\# \h_k^\bullet\bigl(y,y'\bigr)\bigr)
\times\sup_{0\leq i<\lfloor2^k
D(y,y')\rfloor} D^*\bigl(\wt\Gamma\bigl(i2^{-k}\bigr),\wt
\Gamma\bigl((i+1)2^{-k}\bigr)\bigr) \build{\la}_{k\to
\infty
}^{\mathrm{a.s.}}
0.
\]

We can now pass to the limit $k\to\infty$ in (\ref{MainTech3}), using
the fact that the
topologies induced by $D$ and $D^*$ are the same, and we get
\[
D^*\bigl(y,y'\bigr)\leq D^*\bigl(y,\wt\Gamma(u)\bigr) + D
\bigl(y,y'\bigr).
\]
This holds for any $u>0$. Letting $u\to0$, we obtain $D^*(y,y')\leq
D(y,y')$. This
completes the proof since we already know that $D(y,y')\leq D^*(y,y')$.
\end{pf}

Let us state a corollary that will be useful when we deal with the case of
triangulations.

%co7.3 #&#
\begin{corollary}
\label{tritool}
Let $U$ and $V$ be two independent random variables uniformly distributed
over $[0,1]$ and such that the pair $(U,V)$ is independent of $(\be,Z)$. Then,
\[
D^*(U,V)\build{=}_{}^{\mathrm{(d)}} D^*(s_*,U)= Z_U +
\Delta\build {=}_{}^{\mathrm{(d)}} \Delta.
\]
\end{corollary}

\begin{pf}
The second equality is easy from the definition of $D^*$.
The last identity in distribution is a consequence of the invariance of
the CRT under uniform re-rooting;
see, in particular,~\cite{LGW}, Section 2.3. Let us prove the
first identity in distribution. We consider the setting of Theorem \ref
{mainIM}, and we take $p=2$,
as this simplifies the argument a little and suffices for our purposes. Let
$0=i(0)<i(1)<\cdots<i(n)$ be the indices corresponding to the first
visits of vertices of
$\tau_n^\circ$ by the white contour sequence. Then, we have
%
%e49 #&#
\begin{equation}
\label{tritooltech} \sup_{0\leq t< 1} \biggl|\frac{i(\lfloor(n+1)t\rfloor)}{2n} - t \biggr| \build{
\la}_{n\to\infty}^{\mathrm{(P)}} 0,
\end{equation}
where the notation $ \build{\la}_{}^{\mathrm{(P)}}$ indicates convergence in
probability. Noting that $2$-trees are naturally
identified with ordinary plane trees (by removing all black vertices
and putting an edge between two white vertices if
they were adjacent to the same black vertex),
the preceding convergence follows from the standard
arguments used to compare the contour function of a plane tree with its
so-called height function; see, for example, the proof of Theorem 1.17
in~\cite{trees}. From (\ref{tritooltech}) and the convergence (\ref
{basicconv}), it is a simple matter to
obtain that $D^*(U,V)=D(U,V)$ is the limit in distribution of $\kappa_p
n^{-1/4} d_{\mathrm{gr}}(X_n,Y_n)$ where
$X_n$ and $Y_n$ are independently and uniformly distributed over $\mm_n$. However,
this is also the limiting distribution of $\kappa_p n^{-1/4}
d_{\mathrm{gr}}(X_n,\partial)$, which by
(\ref{basicconv}) and (\ref{tritooltech}) again is the distribution of
$D(s_*,U)$.
 \end{pf}

%s8 #&#
\section{The case of triangulations}\label{sec8}

%s8.1 #&#
\subsection{Coding triangulations with trees}
\label{tricoding}

In this section we prove Theorem~\ref{mainresult} in the case $q=3$.
Similarly as in
the bipartite case, we will rely on certain bijections between
triangulations and trees, which we now
describe. These bijections can be found in~\cite{BDG}, and we follow
the presentation of
\cite{CLM}, to which we refer for more details.

Recall the definition of plane trees in Section~\ref{labtrees}. We
will need to
consider $4$-type plane trees. A $4$-type plane tree is just a pair
$(\tau,(\operatorname{typ}(u))_{u\in\tau})$ consisting of a plane tree $\tau
$ and
for every $u\in\tau$ of a type $\operatorname{typ}(u)\in\{1,2,3,4\}$. To simplify
notation, we systematically write $\tau$ instead of $(\tau,(\operatorname{typ}(u))_{u\in\tau})$
in what follows, as we will only be considering $4$-type plane trees. A $T$-tree
is a $4$-type plane tree that satisfies the following properties:
\begin{longlist}[(iii)]
\item[(i)] The root vertex $\varnothing$ is of type $1$ or of type $2$.
\item[(ii)] The children of a vertex of type $1$ are of type $3$.
\item[(iii)] Each vertex of type $2$ and which is not the root
$\varnothing$
has exactly one child of type $4$, and no other child. If the root
$\varnothing$
is of type $2$, it has two children, both of type $4$.
\item[(iv)] Each vertex of type $3$ has exactly one child, which is of
type~$2$.
\item[(v)] Each vertex of type $4$ has either one child of type $1$ or
two children of type~$2$,
and no other child.
\end{longlist}

If $\tau$ is a $T$-tree, we write $\tau^\circ$ for the set of all
vertices of
$\tau$ at even generation. Clearly, this is also the set of all
vertices of
type $1$ or $2$ in $\tau$. By analogy with the bipartite case, we call
the elements of $\tau^\circ$ the white vertices of $\tau$.

Let $\tau$ be a $T$-tree. An admissible labeling of $\tau$ is a
collection of labels
assigned to the white vertices of
$\tau$, such that the following properties hold:
\begin{longlist}[(a)]
\item[(a)] $\ell_\varnothing=0$ and $\ell_v\in\Z$ for each $v\in
\tau^\circ$.
\item[(b)] Let $v\in\tau\setminus\tau^\circ$, let $v_{(0)}$ be the
parent of $v$ and let
$v_{(j)}=vj$, $1\leq j\leq k$, be the children of $v$. Then for every
$j\in\{0,1,\ldots,k\}$,
$\ell_{v_{(j+1)}}\geq\ell_{v_{(j)}}-1$, where by convention
$v_{(k+1)}=v_{(0)}$.
Furthermore, if $j\in\{0,1,\ldots,k\}$ is such that $v_{(j+1)}$ is of
type $2$,
we have $\ell_{v_{(j+1)}}\geq\ell_{v_{(j)}}$.
\end{longlist}
Note the slight difference with the analogous definition in
Section~\ref{labtrees}. As a consequence of property (b), we
observe that
if a vertex $v$ of type $4$ has two children, $v1$ and $v2$ in our
formalism, and
if $u$ is the parent of $v$ (necessarily of type $2$), we have $\ell_u=\ell_{v1}=\ell_{v2}$.

A labeled $T$-tree is a pair consisting of a $T$-tree $\tau$ and an
admissible labeling
$(\ell_u)_{u\in\tau^\circ}$ of $\tau$. For every integer $n\geq
3$, we
write $\T_n$
for the set of all labeled $T$-trees with $n-1$ vertices of type $1$.
It will be
convenient to write $\T_n=\T_n^{(1)}\cup\T_n^{(2)}$, where
$\T_n^{(1)}$, respectively $\T_n^{(2)}$, corresponds to labeled $T$-trees
whose root vertex is of type~$1$, respectively of type $2$.

Let $\mathscr{T}_n$ denote the set of all rooted and pointed triangulations
with $n$ vertices [or, equivalently, $2(n-2)$ faces]. Let $M\in
\mathscr{T}_n$,
let $\partial$ be the distinguished vertex of $M$, and let $e_-$ and $e_+$
be, respectively, the origin and the target of the root edge. As previously,
write $d_{\mathrm{gr}}$ for the graph distance on the vertex set of $M$. The
triangulation
$M$ is said to be positive, respectively null, respectively negative, if
\begin{eqnarray}
d_{\mathrm{gr}}(\partial, e_+)=d_{\mathrm{gr}}(\partial,e_-)+1,\nonumber \\
\eqntext{\mbox{respectively }d_{\mathrm{gr}}(\partial, e_+)=d_{\mathrm{gr}}(\partial,e_-),
\mbox{respectively }d_{\mathrm{gr}}(\partial, e_+)=d_{\mathrm{gr}}(
\partial,e_-)-1.}
\end{eqnarray}
With an obvious notation, we can thus write $\mathscr{T}_n=\mathscr
{T}_n^+\cup\mathscr{T}_n^0\cup\mathscr{T}_n^-$. Note that
reversing the orientation of the root edge gives an obvious bijection
between $\mathscr{T}_n^+$
and $\mathscr{T}_n^-$.

\begin{figure}

\includegraphics{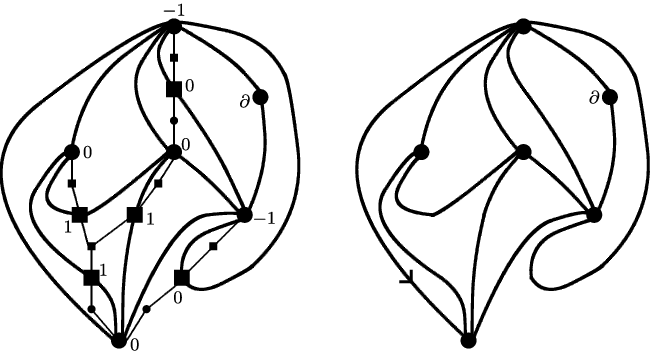}

\caption{A labeled $T$-tree in $\T_n^{(1)}$ and the associated rooted
and pointed triangulation.
Vertices of type $1$ are represented by big black disks, vertices of
type $2$ by big black squares,
vertices of type $3$ by small black disks and vertices of type $4$ by small
black squares.}\label{fig6}
\end{figure}

A special case of the results in~\cite{BDG} yields bijections between
$\T_n^{(1)}$ and
$\mathscr{T}_n^+$ on the one hand, and between $\T_n^{(2)}$ and
$\mathscr{T}_n^0$ on the other hand. Let us describe the first of these
bijections in some detail (see Figure~\ref{fig6} for an example). We start from
a labeled $T$-tree $(\tau,(\ell_u)_{u\in\tau^\circ})\in\T_n^{(1)}$,
and let $k$ be the number of edges of~$\tau$. The white contour
sequence of $\tau$ is the finite sequence
$(v_0,v_1,\ldots,v_k)$ defined exactly as in Section~\ref{labtrees},
and we set
$v_{k+i}=v_i$ for $1\leq i\leq k$. Note that
every corner of a white vertex $v$ of $\tau$ corresponds exactly to
one index $i\in\{0,1,\ldots,k-1\}$, such that $v_i=v$, and we call this
corner the
corner $v_i$ as we did previously. We assume that the tree $\tau$ is
drawn on the sphere, and
as in Section~\ref{BDGbij} we add an extra vertex $\partial$, which
is of type $1$
by convention. We then draw edges of the map according to the very same
rules as
in Section~\ref{BDGbij}: If $i\in\{0,1,\ldots,k-1\}$ is such that
$\ell_{v_i}=\min\{\ell_v,v\in\tau^\circ\}$,
we draw an edge between the corner $v_i$ and $\partial$, and, on the
other hand,
if $i$ is such that $\ell_{v_i}>\min\{\ell_v,v\in\tau^\circ\}$, we
draw an edge between
the corner $v_i$ and the corner $v_j$, where
$j$ is the successor of $i$. Note that in the latter case the vertex
$v_j$ is of
type $1$. This follows from the fact that if the vertex $v_m$ is of
type $2$, we have always
$\ell_{v_m}\geq\ell_{v_{m-1}}$ by our labeling rules.\vadjust{\goodbreak}

By property (iii) in the definition of a $T$-tree, each vertex $v$ of
type $2$ in $\tau$
has exactly two corners and, therefore, the preceding device will give
exactly two edges connecting $v$ to vertices of type $1$.
To complete the construction, we erase all vertices of type $2$ and for
each such vertex $v$, we merge the two edges incident to $v$ into a
single edge
connecting two vertices of type~$1$ (which may be the same). In this
way, we obtain
a planar map $M$ whose vertex set consists of all vertices of type~$1$
(including $\partial$), which is easily checked to be a triangulation.
This triangulation
is pointed at $\partial$ and rooted at the edge generated by the case
$i=0$ of
the construction. This edge is oriented so that its target is the
vertex $\varnothing$.
The mapping
\[
\bigl(\tau,(\ell_u)_{u\in\tau^\circ}\bigr)\la M
\]
that we have just described is a bijection from $\T_n^{(1)}$
onto $\mathscr{T}_n^+$.

%f6 #&#

A minor modification of this construction yields
a bijection from $\T_n^{(2)}$ onto $\mathscr{T}_n^0$. Edges
of the map are generated in the same way, but the root edge
of the map is now obtained
as the edge resulting of the merging of the two edges incident to
$\varnothing$ [recall that for a tree in $\T_n^{(2)}$ the root
$\varnothing$
is a vertex of type $2$ that has exactly two children, hence also two corners].
The orientation of the root edge is chosen by deciding that the
``half-edge'' coming from
the first corner of $\varnothing$ corresponds to the origin of the
root edge.

In both cases, distances in the planar map $M$ satisfy the following
analog of (\ref{formuladgr}): For every vertex $v$ of type $1$ in
$\tau
$, we have
%
%e50 #&#
\begin{equation}
\label{distance-label} d_{{\mathrm{gr}}}(\partial,v)= \ell_v -
\min\ell+1,
\end{equation}
where $\min\ell$ denotes the minimal label on the tree $\tau$.
In the left-hand side
$v$ is viewed as a vertex of the map $M$, in agreement with the
preceding construction.

An analog of (\ref{bounddgr}) also holds (again with a
proof very similar to that of~\cite{IM}, Lemma 3.1).
If $v$ and $v'$ are two vertices of type $1$ of $\tau$, such that
$v=v_i$ and $v'=v_j$ for some $i,j\in\{0,1,\ldots,k\}$ with $i< j$,
we have
%
%e51 #&#
\begin{equation}
\label{bounddgr2} d_{\mathrm{gr}}\bigl(v,v'\bigr) \leq
\ell_{v_i} + \ell_{v_j} -2 \max \Bigl(\min_{i\leq
m\leq
j}
\ell_{v_m}, \min_{j\leq m\leq i+k} \ell_{v_m} \Bigr)+2.
\end{equation}

%s8.2 #&#
\subsection{Random triangulations}

Following~\cite{MieInvar}, we now want to interpret a random labeled
$T$-tree uniformly distributed
over $\T_n^{(1)}$ as a conditioned multitype Galton--Watson tree (a
similar interpretation will hold
for a $T$-tree uniformly distributed
over $\T_n^{(2)}$). We set
\[
\beta= 1-\tfrac{\sqrt{3}}{3},\qquad  \alpha=\tfrac{1}{2} + \tfrac{\sqrt{3}}{6}
\]
and we let $\mu$ be the geometric distribution with parameter $\beta$:
\[
\mu(k)= (1-\beta) \beta^k
\]
for $k=0,1,\ldots.$ We let $\mathbf{t}$ be a random $T$-tree
satisfying the
following prescriptions:
\begin{itemize}
\item the root vertex is of type $1$;
\item each vertex of type $1$ has, independently of the
other vertices, a random number
of children distributed according to $\mu$;
\item each vertex of type $4$ has, independently of the
other vertices, either one child of type $1$
with probability
$\alpha$ or two children of type $2$ with probability $1-\alpha$.
\end{itemize}

Recalling the properties of the definition of a $T$-tree, one sees that the
preceding prescriptions completely characterize the distribution of
$\mathbf{t}$.
The random tree $\mathbf{t}$ can be viewed as a $4$-type Galton--Watson tree
in the sense of~\cite{MieGW}. Note that this Galton--Watson tree is
critical, meaning that
the spectral radius of the mean matrix of offspring distributions is
equal to~$1$. This property ensures that the $4$-type Galton--Watson tree
with these offspring distributions is finite a.s., a fact that is needed
for the existence of $\mathbf{t}$ as above (our $T$-trees are finite
by definition).
We write $\mathbf{t}^{(1)}$ for the set of all vertices of type $1$ of
$\mathbf{t}$.

Let $\bT$ be a random labeled $T$-tree obtained by assigning an
admissible labeling
to $\mathbf{t}$, uniformly at random over all possibilities.

%le8.1 #&#
\begin{lemma}
\label{GWtriangu}
Let $n\geq3$. The conditional distribution of $\bT$ knowing that
$\#\mathbf{t}^{(1)}=n-1$ is uniform over $\T_n^{(1)}$.
\end{lemma}

This lemma is a very special case of Proposition 3 in~\cite{MieInvar}, but
it is also easy to give a direct proof. Note that the values of
$\beta$ and $\alpha$ are chosen (in a unique way) so that the tree
$\mathbf{t}$ is critical and the result of the lemma holds.\looseness=1

Using Lemma~\ref{GWtriangu}, we can now apply the invariance
principles of
\cite{MieGW} to study the asymptotics of the contour and label processes
of a tree uniformly distributed over $\T_n^{(1)}$. So let
$(\tau_n,(\ell^n_v)_{v\in\tau_n^\circ})$ be a random labeled $T$-tree
uniformly distributed over $\T_n^{(1)}$ and let $k_n$
be the number of edges of $\tau_n$. If
$v^n_0,v^n_1,\ldots,v^n_{k_n}$ is the white contour sequence of $\tau_n$, the
contour process $(C^n_i)_{0\leq i\leq k_n}$ and the
label process $(\Lambda^n_i)_{0\leq i\leq k_n}$ are defined by
$C^n_i=\frac{1}{2}|v^n_i|$ and $\Lambda^n_i=\ell^n_{v^n_i}$ as in
the bipartite
case, and they are extended to the real interval $[0,k_n]$
by linear interpolation. It is also useful to define
$L^n_j$, for $0\leq j\leq k_n$, as the number of distinct vertices of
type $1$
among $v^n_0,v^n_1,\ldots, v^n_j$.

%pr8.2 #&#
\begin{proposition}
\label{contourtriangu}
Set $\lambda_{3/2}=\frac{1}{4} (3-\sqrt{3})$ and $\kappa_{3/2}=3^{1/4}$. We have
%
%e52 #&#
\begin{equation}
\label{convtricontour} \bigl( \lambda_{3/2} n^{-1/2}
C^n_{ k_nt}, \kappa_{3/2} n^{-1/4}
\Lambda^n_{k_nt} \bigr)_{0\leq t\leq1} \build{
\la}_{n\to\infty}^{\mathrm{(d)}} (\be_t, Z_t)_{0\leq t\leq1}.
\end{equation}
Furthermore,
%
%e53 #&#
\begin{equation}
\label{type1} \sup_{0\leq t\leq1} \bigl|n^{-1} L^n_{\lfloor k_nt\rfloor}
-t\bigr| \build {\la }_{n\to\infty}^{\mathrm{(P)}} 0.
\end{equation}
\end{proposition}

Using Lemma~\ref{GWtriangu}, Proposition~\ref{contourtriangu} can be
derived as
a special case of Theorems~2 and~4 in~\cite{MieGW}. In order to apply
these results, we need labels to be assigned to every vertex of $\tau_n$,
and not only to white vertices, but we can just decide that every
vertex of
type $3$ or $4$ is assigned the label of its parent. Moreover, it is assumed
in~\cite{MieGW} that the vectors of label increments (meaning the vectors
obtained by considering
the differences between the labels of the children of a vertex and the
label of this vertex) are centered,
which is not true here. However, as pointed out in~\cite{MieInvar}
in a more general setting, a very minor modification
makes the vectors of label increments centered: Just subtract $\frac
{1}{2}$ from the label
of every vertex of type~$2$ (and from the label of its unique child of
type $4$).
Obviously this modification has no effect on the validity of the convergence
in the proposition.

We also note that~\cite{MieGW}
considers the so-called height process, rather than the contour
process, and
the corresponding variant of the label process. However, as we already
mentioned in
the proof of Theorem~\ref{mainIM}, limit
theorems for the height process can be translated easily in terms of
the contour process
(see, e.g., Section~1.6 in~\cite{trees}).

At this point, it is appropriate to comment on the value of the constants
$\lambda_{3/2}$ and $\kappa_{3/2}$, since the corresponding
discussion in
\cite{MieInvar} seems to contain a miscalculation. We use the notation of\vadjust{\goodbreak}
\cite{MieGW}. The mean matrix of the offspring distributions is
\[
\pmatrix{ 0&0&\sqrt{3}-1&0
\vspace*{2pt}\cr
0&0&0& 1 \vspace*{2pt}\cr
0&1&0&0
\vspace*{2pt}\cr
\frac{1}{2} + \frac{\sqrt{3}}{6}&1-\frac{\sqrt{3}}{3}&0&0 },
\]
and the associated left and right eigenvectors are
\[
a=\tfrac{1}{6-\sqrt{3}} (1,3-\sqrt{3},\sqrt{3}-1,3-\sqrt{3}),\qquad b =
\tfrac{6-\sqrt{3}}{4} (\sqrt{3}-1,1,1,1).
\]
The quadratic forms $Q^{(i)}$, for $i=1,2,3,4$, are easily computed as
\begin{eqnarray*}
Q^{(1)}(s_1,s_2,s_3,s_4)&=&
2(\sqrt{3}-1)^2 s_3^2,\qquad Q^{(2)}=Q^{(3)}=0,
\\
Q^{(4)}(s_1,s_2,s_3,s_4)
&= &\bigl(1-\tfrac{\sqrt{3}}{3}\bigr)s_2^2.
\end{eqnarray*}
A simple calculation gives
\[
a\cdot Q(b)=\sum_{i=1}^4
a_i Q^{(i)}(b)= \frac{1}{8}(6-3\sqrt {3}) (6-
\sqrt{3}),
\]
and it follows that
\[
\lambda_{3/2}= \sqrt{a\cdot Q(b)}\times\sqrt{a_1}= \sqrt{
\tfrac
{6-3\sqrt
{3}}{8}}= \tfrac{1}{4} (3-\sqrt{3})
\]
(comparing with Theorem 2 in~\cite{MieGW}, the formula for $\lambda_{3/2}$ has an extra multiplicative
factor $2$ corresponding to the factor $\frac{1}{2}$ in the definition
of the contour process).

To compute $\kappa_{3/2}$, we then need to evaluate the quantity
$\Sigma
$ in Theorem~4 of~\cite{MieInvar}.
We note that the only label increments having nonzero variance
correspond either to
a vertex of type $3$ (having automatically one child of type $2$) or to
a vertex of type~$4$ having only one
child of type $1$, with probability $\alpha$. In both cases the
variance is $\frac{1}{4}$. It follows that
\begin{eqnarray*}
\Sigma^2&=& a_3\times b_2\times
\tfrac{1}{4} + \alpha a_4\times b_1 \times
\tfrac{1}{4}
\\
&=&\tfrac{1}{16} \bigl((\sqrt{3}-1)+ \bigl(\tfrac{1}{2} +
\tfrac{\sqrt{3}}{6}\bigr) (3-\sqrt{3}) (\sqrt{3}-1) \bigr) = \tfrac{1}{8}(
\sqrt{3}-1)
\end{eqnarray*}
and
\[
\kappa_{3/2}= \frac{1}{\Sigma} \times\sqrt{\frac{ \lambda_{3/2}}{2}}=
3^{1/4}.
\]

\begin{rema*} There are more direct ways of computing the constant
$\kappa_{3/2}$, for
instance, by considering the genealogical tree associated
with vertices of type $1$ (say that a vertex $u$ of type $1$ is a child
of another vertex $v$ of type~$1$ if $v$ is the last vertex of type~$1$
that is an ancestor of $u$
distinct from $u$). It turns out that this
tree is also a conditioned Galton--Watson tree, whose offspring
distribution can be computed
easily.
\end{rema*}

In this subsection we concentrated on the case of a labeled $T$-tree
uniformly distributed
over $\T_n^{(1)}$. However, Proposition~\ref{contourtriangu} remains
valid if we replace
$(\tau_n,(\ell^n_v)_{v\in\tau_n^\circ})$ by a random labeled $T$-tree
uniformly distributed~over~$\T_n^{(2)}$. The proof is the same up to minor modifications.

%s8.3 #&#
\subsection{Convergence of rescaled triangulations to the Brownian map}
\label{cotrimap}

We will now prove the case $q=3$ of Theorem~\ref{mainresult}. We
consider a
random triangulation $M_n$ uniformly distributed over $\mathscr
{T}_n^+$, and write
$\mm_n$ for the vertex set of $M_n$. We will prove that
%
%e54 #&#
\begin{equation}
\label{convtri} \bigl(\mm_n,\kappa_{3/2}n^{-1/4}
d_{\mathrm{gr}}\bigr) \build{\la}_{n\to\infty}^{\mathrm{(d)}} \bigl(
\mm_\infty, D^*\bigr).
\end{equation}
Obviously the same result holds if $M_n$ is uniformly distributed over
$\mathscr{T}_n^-$,
and only minor modifications would be needed to handle the case when
$M_n$ is uniformly distributed over $\mathscr{T}_n^0$. Combining all
three cases, and
using the fact that a triangulation with $n$ faces has $\frac{n}{2} +
2$ vertices, we obtain
the case $q=3$ of Theorem~\ref{mainresult}.

In order to prove (\ref{convtri}), we may assume that $M_n$ is the
image of
a random labeled $T$-tree $(\tau_n,(\ell^n_v)_{v\in\tau_n^\circ})$
uniformly distributed
over $\T_n^{(1)}$ under the bijection described in Section \ref
{tricoding}. We will
rely on Proposition~\ref{contourtriangu}, and we use the notation introduced
before this proposition. Recall from Section~\ref{tricoding} that the
bijection between triangulations and labeled $T$-trees allows us to
identify
\[
\mm_n= \tau_{n}^{(1)} \cup\{\partial\},
\]
where $\tau_{n}^{(1)}$ is the set of all vertices of type $1$ in $\tau_n$. We also set
\[
\mm'_n= \tau_{n}^{\circ} \cup\{
\partial\}
\]
and we note that the graph distance on $\mm_n$ can be extended to $\mm'_n$
in the following way. Let $u\in\mm'_n\setminus\mm_n$, then
$u$ is a vertex of type $2$ in $\tau_n$ and, as already mentioned, $u$
has two
successors $u'$ and $u''$ (possibly such that $u'=u''$), which are
vertices of type $1$.
If $v\in\mm_n$, we set
\[
d_{\mathrm{gr}}(v,u)=d_{\mathrm{gr}}(u,v)= \tfrac{1}{2} + \min
\bigl(d_{\mathrm{gr}}\bigl(u',v\bigr),d_{\mathrm{gr}}
\bigl(u'',v\bigr)\bigr).
\]
If $v\in\mm'_n\setminus\mm_n$ and $v\not= u$, we put
\[
d_{\mathrm{gr}}(u,v)= 1+ \min \bigl(d_{\mathrm{gr}}\bigl(u',v'
\bigr),d_{\mathrm{gr}}\bigl(u',v''
\bigr),d_{\mathrm{gr}}\bigl(u'',v'
\bigr),d_{\mathrm{gr}}\bigl(u'',v''
\bigr)\bigr),
\]
where $v'$ and $v''$ are the two successors of $v$.
It is straightforward to verify that $d_{\mathrm{gr}}$ thus extended is a
distance on $\mm'_n$.
Informally, we may interpret the preceding definition by saying that
in the triangulation $M_n$ we have added a new vertex at the middle of
each edge connecting two vertices at the same distance from $\partial$,
and we agree that
this new vertex is at distance $\frac{1}{2}$ from both ends of the edge
where it has
been created.

Clearly, it is enough to prove that (\ref{convtri}) holds when $\mm_n$
is replaced
by $\mm'_n$ or even by $\mm_n'':=\mm'_n\setminus\{\partial\}$. Let
$u,v\in\mm''_n$, and suppose that $u=v^n_i$ and $v=v^n_j$ for some
$i,j\in\{0,1,\ldots,k_n\}$ with $i\leq j$. Then, we get from (\ref
{bounddgr2}) that
%
%e55 #&#
\begin{equation}
\label{bounddgr3} d_{\mathrm{gr}}(u,v) \leq d_{n}^\circ(i,j),
\end{equation}
where
\[
d_n^\circ(i,j)=d_n^\circ(j,i)=
\Lambda^n_i + \Lambda^n_j -2
\max \Bigl(\min_{k\in\{i,\ldots,j\}} \Lambda^n_k,
\min_{k\in\{j,\ldots,k_n\}\cup\{0,\ldots,i\}} \Lambda^n_k \Bigr)+2.
\]
We also set $d_n(i,j)=d_{\mathrm{gr}}(v^n_i,v^n_j)$ for every $i,j\in\{
0,1,\ldots,k_n\}$, and we extend the definition
of both $d_n$ and $d^\circ_n$ to $[0,pn]\times[0,pn]$ by linear
interpolation. By (\ref{bounddgr3}),
we have $d_n\leq d^\circ_n$. On the other hand, it immediately follows from
Proposition~\ref{contourtriangu} that
%
%e56 #&#
\begin{equation}
\label{tritech1} \bigl(\kappa_{3/2}n^{-1/4}d^\circ_n(k_ns,k_nt)
\bigr)_{0\leq s\leq
1,0\leq
t\leq1} \build{\la}_{n\to\infty}^{\mathrm{(d)}}
\bigl(D^\circ(s,t) \bigr)_{0\leq
s\leq
1,0\leq t\leq1},
\end{equation}
where $D^\circ$ is as in Section~\ref{CRTlabels}, and this
convergence holds jointly with (\ref{convtricontour}). From the convergence
(\ref{tritech1}) and the bound $d_n\leq d^\circ_n$, the same argument as
in the proof of Proposition 3.2 in~\cite{IM} shows that the sequence
of the laws of the processes
\[
\bigl(n^{-1/4}d_n(k_ns,k_nt)
\bigr)_{0\leq s\leq1,0\leq t\leq1}
\]
is tight in the space of probability measures on $C([0,1]^2,\R)$.
Hence, from any
monotone increasing sequence of positive integers, we can extract a subsequence
$(n_j)_{j\geq1}$ along which we have the convergence in distribution
%
%e57 #&#
\begin{eqnarray}
\label{basicconv3} && \bigl(\lambda_{3/2} n^{-{1}/{2}}
C^n_{k_nt}, \kappa_{3/2} n^{-{1}/{4}}
\Lambda^n_{k_nt}, \kappa_{3/2} n^{-{1}/{4}}
d^\circ_n(k_ns,k_nt),
\kappa_{3/2} n^{-{1}/{4}} d_n(k_ns,
k_nt) \bigr)
\nonumber\hspace*{-30pt}
\\[-4pt]
\\[-12pt]
\nonumber
 &&\qquad \build{\la}_{n\to\infty}^{\mathrm{(d)}} \bigl(
\be_t,Z_t,D^\circ(s,t),
D'(s,t) \bigr),
\nonumber\hspace*{-30pt}
\end{eqnarray}
where $(D'(s,t))_{0\leq s\leq1,0\leq t\leq1}$ is a random process
such that
$D'\leq D^\circ$. Using the Skorokhod representation theorem, we may
and will
assume that the convergence~(\ref{basicconv3}) holds a.s., uniformly on
$[0,1]^2$,
along the sequence $(n_j)_{j\geq1}$. Since $d_n$ is symmetric and
satisfies the
triangle inequality, one immediately obtains that $D'$ is a (random)
pseudometric
on $[0,1]$.

We claim that
%
%e58 #&#
\begin{equation}
\label{D'claim} D'(s,t)=D^*(s,t) \qquad\mbox{for every
}s,t\in[0,1]\mbox{ a.s.}
\end{equation}
To prove this, we start by observing that we have $D'(s,t)=0$ for every
$s,t\in[0,1]$ such that
$p_\be(s)=p_\be(t)$, a.s. This follows by exactly the same argument as
in the proof of~\cite{IM}, Proposition 3.3(iii).
Recalling the definition of $D^*$ in Section~\ref{CRTlabels}, and
using the fact that
$D'$ satisfies the triangle inequality, we see that the property
$D'\leq D^\circ$
implies $D'(s,t)\leq D^*(s,t)$ for every $s,t\in[0,1]$.

Then let $U$ and $V$ be two independent random variables uniformly distributed
over $[0,1]$ and such that the pair $(U,V)$ is independent of all other
random quantities. By a continuity argument, our claim (\ref{D'claim})
will follow if we can verify that $D'(U,V)=D^*(U,V)$ a.s. Since
$D'(U,V)\leq D^*(U,V)$, it will be sufficient
to verify that $D'(U,V)$ and $D^*(U,V)$ have the same distribution. To
see this, let $0=i(1)<i(2)<\cdots<i(n-1)$
be the first visits by the white contour sequence of $\tau_n$ of the
vertices of type $1$
in $\tau_n$. Also set $U_n=\lceil(n-1)U\rceil$ and $V_n=\lceil
(n-1)V\rceil$, which are both
uniformly distributed over $\{1,2,\ldots,n-1\}$. It follows from (\ref
{type1}) that
\[
\frac{i(U_n)}{k_n} \build{\la}_{n\to\infty}^{\mathrm{(P)}} U,\qquad
\frac{i(V_n)}{k_n} \build{\la}_{n\to\infty}^{\mathrm{(P)}} V.
\]
Together with (\ref{basicconv3}), this now implies that
\[
\kappa_{3/2} n^{-1/4} d_{\mathrm{gr}}\bigl(v^n_{i(U_n)},v^n_{i(V_n)}
\bigr) \build{\la }_{n\to\infty}^{\mathrm{(P)}} D'(U,V)
\]
as $n\to\infty$ along the sequence $(n_j)_{j\geq1}$.
Hence, the distribution of $D'(U,V)$ is the limiting distribution
[along $(n_j)_{j\geq1}$]
of $\kappa_{3/2} n^{-1/4} d_{\mathrm{gr}}(X_n,Y_n)$,
where $X_n$ and $Y_n$ are independently uniformly distributed over $\mm_n$.
Obviously, this is also the limiting distribution of
\[
\kappa_{3/2} n^{-1/4} d_{\mathrm{gr}}\bigl(\partial,
v^n_{i(U_n)}\bigr)=\kappa_{3/2} n^{-1/4}
\bigl( \Lambda^n_{i(U_n)} - \min\bigl\{\Lambda^n_i\dvtx 0
\leq i\leq k_n\bigr\} +1 \bigr)
\]
using (\ref{distance-label}) in the last equality. From (\ref
{basicconv3}), we now get that $D'(U,V)$ has the same distribution as
$Z_U + \Delta$. Corollary~\ref{tritool} shows that this is the same as
the distribution of $D^*(U,V)$, thus completing the
proof of (\ref{D'claim}).

Recall from Theorem~\ref{mainIM} that $s\approx t$ if and only if
$D^*(s,t)=0$, and that
$\mm_\infty$ is the quotient space $[0,1]/\approx$.
Once we know that $D'=D^*$, it is an easy matter to deduce from the
(almost sure) convergence (\ref{basicconv3})
that we have, along the sequence $(n_j)_{j\geq1}$,
\[
\bigl(\mm''_{n},\kappa_{3/2}n^{-1/4}d_{\mathrm{gr}}
\bigr)\build{\la}_{n\to\infty
}^{\mathrm{a.s.}} \bigl(\mm_\infty, D^*
\bigr)
\]
in the Gromov--Hausdorff sense. To see this, define a correspondence
between the metric spaces $(\mm''_{n},\kappa_{3/2}n^{-1/4}d_{\mathrm{gr}})$ and
$(\mm_\infty, D^*)$ by saying that
a vertex $v\in\mm''_n$ is in correspondence with the equivalence class
of $s\in[0,1]$
if and only if $v=v^n_i$, where $i=\lfloor k_ns\rfloor$. From
(\ref{basicconv3}), the distortion of this correspondence tends to $0$
along the sequence $(n_j)_{j\geq1}$, a.s., and this gives the desired
Gromov--Hausdorff
convergence.

Consequently,
we have obtained that from any monotone increasing sequence of positive
integers one
can extract a subsequence along which (\ref{convtri}) holds. This
suffices for the
desired result.

\begin{rema*} The argument of the preceding proof could also be used to
deduce the
convergence of Theorem~\ref{mainresult} for all even values of $q\geq
4$ from
the special case $q=4$ (note that the statement of Corollary~\ref{tritool}
already follows
from this special case).\vadjust{\goodbreak} We have chosen not to do so because restricting
ourselves to quadrangulations would not simplify much the proof in the
bipartite case
and because some of the intermediate results that we derive for
$2p$-angulations are
of independent interest.
\end{rema*}

%s9 #&#
\section{Extensions and problems}\label{sec9}

%s9.1 #&#
\subsection{Boltzmann weights on bipartite planar maps}

The argument we have used to handle triangulations can be applied to
other classes of random planar maps. In this paragraph we briefly discuss
Boltzmann distributions on bipartite planar maps, which have been studied
by Marckert and Miermont~\cite{MaMi}. We consider a sequence $\bw
=(\bw_i)_{i\geq1}$
of nonnegative real numbers, such that there exists at least one
integer $i\geq2$
such that $\bw_i>0$. We assume that the sequence $\bw$ is regular
critical in
the sense of~\cite{MaMi}.

For every integer $n\geq2$, we let $\b_n$ stand for the set of all
rooted bipartite planar maps with $n$
vertices. Recall that a planar map is bipartite if and only if all its
faces have even degree.
If $M$ is a planar map, we denote the set of all its faces by $F(M)$,
and for every face
$f$ of $M$ we write $\operatorname{deg}(f)$ for the degree of the face $f$.

%th9.1 #&#
\begin{theorem}
\label{Boltzmann}
For every large enough integer $n$, let $P^\bw_n$ denote the unique probability
measure on $\b_n$ such that, for every $M\in\b_n$,
\[
P^\bw_n(M)= c_{\bw,n} \prod
_{f\in F(M)} \bw_{\operatorname{deg}(f)/2},
\]
where $c_{\bw,n}$ is a constant depending only on $\bw$ and $n$.
Let $M_n$ be a random planar map distributed according to $P^\bw_n$. Then,
if $V(M_n)$ stands for the vertex set of $M_n$ and $d_{\mathrm{gr}}$ is the
graph distance on $V(M_n)$, there exists a constant $a_\bw>0$ such that
\[
\bigl(V(M_n), a_\bw n^{-1/4}d_{\mathrm{gr}}
\bigr)\build{\la}_{n\to\infty}^{\mathrm{(d)}} \bigl(\mm_\infty, D^*
\bigr)
\]
in the Gromov--Hausdorff sense.
\end{theorem}

This theorem can be proved along the lines of Section~\ref{cotrimap}.
The main technical tool
is the BDG bijection for bipartite planar maps, as described in
Section~2.3 of~\cite{MaMi} (this is very close to the BDG bijection
for $2p$-angulations described above).
Analogously to Lemma~\ref{GWtriangu},
\cite{MaMi}, Proposition 7, allows us to interpret
the tree associated with $M_n$ as a conditioned
$2$-type Galton--Watson tree. Then~\cite{MaMi}, Theorem 8,
yields an analog of Proposition~\ref{contourtriangu}.
The remaining part of the proof is similar to Section~\ref{cotrimap},
and we will leave the details to the reader. In contrast with
Theorem~\ref{mainresult}, there is in general no explicit formula
for the constant $a_\bw$; see, however, the discussion in
Section~3.2 of~\cite{MaMi}.

\begin{rema*} In Theorem~\ref{Boltzmann} we consider random planar maps
having a fixed large number of vertices. In the case of
$q$-angulations treated in\vadjust{\goodbreak} Theorem~\ref{mainresult}, Euler's relation
shows that
conditioning on the number of vertices is equivalent to conditioning on the
number of faces. This is no longer true for general Boltzmann weights.
The results of
\cite{MaMi} are stated for both kinds of conditionings, but they are
concerned with
rooted and pointed planar maps. In proving Theorem~\ref{Boltzmann} one
implicitly
uses the (obvious) fact that for a planar map in $\b_n$ there are
exactly $n$
possibilities of choosing a distinguished vertex in order to get a
rooted and pointed
planar map.
\end{rema*}

%s9.2 #&#
\subsection{Brownian maps with geodesic boundaries}

In Proposition~\ref{convergenceDMGB} we saw that, along a suitable
sequence $(n_k)_{k\geq1}$,
the rescaled DMGBs associated with uniformly distributed
$2p$-angulations with $n$ edges converge to
a limiting random metric space, which was identified
in Proposition~\ref{identBMGB}. We may now remove the restriction to a
subsequence.

We keep the setting of Section~\ref{scalim}. In particular, $M_n$ is
a rooted
$2p$-angulation uniformly distributed over $\m^p_n$, $\wt M_n$ is the
associated
DMGB as defined in Section~\ref{Dimageobo}, $\wt\mm_n$
is the vertex set of $\wt M_n$ and $\wt d_{\mathrm{gr}}$ is the graph distance
on $\wt\mm_n$.
We also let $(\mm_\infty^\bullet,D^\bullet)$ be the random metric space
obtained via the construction of the end of Section~\ref{scalim} with
$D=D^*$.

%pr9.2 #&#
\begin{proposition}
\label{convBMGB}
We have
\[
\bigl(\wt\mm_n,\kappa_pn^{-1/4}\wt
d_{\mathrm{gr}}\bigr)\build{\la}_{n\to\infty
}^{\mathrm{(d)}} \bigl(
\mm_\infty^\bullet,D^\bullet\bigr)
\]
in the Gromov--Hausdorff sense.
\end{proposition}

\begin{pf} This readily follows from Propositions~\ref{convergenceDMGB}
and~\ref{identBMGB} once we know
that $D=D^*$ in these statements. Indeed, Proposition \ref
{convergenceDMGB} shows that from any
monotone increasing sequence of integers, we can extract a subsequence
along which the
convergence of the proposition holds, and Proposition~\ref{identBMGB}
shows that the limiting law
is uniquely determined as the law of $(\mm_\infty^\bullet,D^\bullet)$.
\end{pf}

The limiting random metric space in Proposition~\ref{convBMGB} may be
called the
Brownian map with geodesic boundaries. As a motivation for the
preceding statement, we expect that this random metric space
will play a significant role in the study of further properties of the
Brownian map.

\begin{rema*} A result analogous to Proposition~\ref{convBMGB} holds for
uniformly distributed triangulations. Since we did not introduce DMGBs
in the
setting of triangulations, we will leave this statement to the reader.
\end{rema*}

%s9.3 #&#
\subsection{Questions}

It is very plausible that Theorem~\ref{mainresult} holds for uniformly
distributed $q$-angulations
for any choice of the integer $q$ (and even for nonbipartite planar
maps distributed according
to Boltzmann weights satisfying suitable conditions). Extending the\vadjust{\goodbreak}
proof we gave in the case
of triangulations would require an analog of Proposition \ref
{contourtriangu} for the random
trees associated with uniformly distributed $q$-angulations. As
observed by Miermont
\cite{MieInvar}, such an analog holds, but only for a ``shuffled''
version of the trees, and this is
not sufficient for our purposes. The reason why a shuffling operation
is needed
is the fact that the vectors of label increments in the trees
associated with $q$-angulations
are no longer centered when $q$ is odd and $q\geq5$. Nonetheless, it
is likely
that one can avoid the shuffling operation and get a full analog of
Proposition~\ref{contourtriangu}.

Another interesting question is to extend Theorem~\ref{mainresult} to
triangulations
satisfying additional connectedness properties (and, in particular, to
type-II or
type-III triangulations in the terminology of~\cite{ADJ}). Via the BDG
bijections,
this would lead to analyzing labeled $T$-trees with extra constraints,
for which
it is again plausible but not obvious that an analog of Proposition
\ref
{contourtriangu}
holds.

Finally, our results raise a number of interesting questions about
Brownian motion
indexed by the CRT. Note that the functions $D^\circ$ and $D^*$
are defined in terms of the pair $(\be,Z)$. Two crucial properties
of these random functions are
%
%e59 #&#
\begin{equation}
\label{cru1}\quad D^*(a,b)=0 \quad\mbox{if and only if}\quad D^\circ(a,b)=0\qquad
\mbox{for every }a,b\in{\cal T}_\be\mbox{ a.s.}
\end{equation}
and
%
%e60 #&#
\begin{equation}
\label{cru2} D^*(U,V)\build{=}_{}^{\mathrm{(d)}} D^*(s_*,U),
\end{equation}
where $U$ and $V$ are independent and uniformly distributed over $[0,1]$
and independent of the pair $(\be,U)$. The equivalence (\ref{cru1}) is
proved in~\cite{IM}, Theorem 3.4,
and~(\ref{cru2}) appears in Corollary~\ref{tritool} above. In both
cases, the proof
relies on the use of approximating labeled trees and the associated
random planar maps.
Since~(\ref{cru1}) and~(\ref{cru2}) are properties of the pair $(\be,Z)$, it would seem
desirable to have a more direct argument for these statements. A direct
proof of (\ref{cru2}),
in particular, would yield a simpler approach to our main result
Theorem~\ref{mainresult}
along the lines of Section~\ref{cotrimap}.

%sA #&#
\begin{appendix}\label{app}
\section*{Appendix}
In this appendix we prove Lemmas~\ref{auxil} and~\ref{snakelemma},
which are both concerned
with properties of the Brownian snake. It will be convenient to argue under
the excursion measure $\N_0$ of the one-dimensional Brownian snake (see
\cite{Zu}).
Recall that the Brownian snake $(W_s)_{s\geq0}$ is a strong Markov
process taking
values in the space $\W$ of all stopped paths. Here a stopped path is simply
a continuous map $w\dvtx [0,\zeta]\la\R$, where $\zeta=\zeta_{(w)}$ is
called the lifetime
of $w$. We write $\wh w=w(\zeta_{(w)})$ for the endpoint of the path $w$.
We may and will assume that $(W_s)_{s\geq0}$ is defined on the
canonical space
$C(\R_+,\W)$ of continuous functions from $\R_+$ into $\W$, and we write
$\zeta_s:=\zeta_{(W_s)}$ for the lifetime of~$W_s$.
Under $\N_0$, the
lifetime process $(\zeta_s)_{s\geq0}$ is distributed
according to It\^o's measure of positive excursions of linear Brownian motion
[normalized so that $\N_0(\sup\{\zeta_s\dvtx s\geq0\}>\ve)=(2\ve)^{-1}$,
for every $\ve>0$]. We use the notation $\sigma=\sup\{s\geq0\dvtx \zeta_s>0\}$
for the duration of the excursion and, for every
$r> 0$, we set $S_r=\inf\{s\geq0\dvtx \wh W_s=-r\}$. This is consistent
with our previous notation
since under the conditional measure $\N_0(\cdot\mid\sigma=1)$ the pair
$(\zeta_s,\wh W_s)_{0\leq s\leq1}$ has the
same distribution as the process $(\be_s,Z_s)_{0\leq s\leq1}$ of the
preceding
sections. For every $t\geq0$,
$\g_t$ denotes the $\sigma$-field generated by $(W_s,0\leq s\leq t)$.
We will use
the explicit form of the distribution of $W_{S_r}$ under $\N_0$, which follows
from the results of~\cite{DLG}, Section 4.6. We first recall~\cite{Zu}, page 91,
that\looseness=1
%
%eA.1 #&#
\begin{equation}
\label{lawmin} \N_0(S_r<\infty) =\N_0
\Bigl(\inf_{s\geq0} \wh W_s \leq-r \Bigr) =
\frac{3}{2r^2}.
\end{equation}\looseness=0
If $(B_t)_{t\geq0}$ denotes a standard linear Brownian motion, the
random path
$(W_{S_r}(t),0\leq t\leq\zeta_{S_r})$ is distributed under $\N_0(\cdot
\mid S_r<\infty)$
as the solution of the stochastic differential equation
\[
\cases{ %
 dX_t= dB_t -
\displaystyle\frac{2}{r+X_t}\,dt,
\vspace*{2pt}\cr
X_0=0, }
\]
stopped when it first hits $-r$. Equivalently, $(r+W_{S_r}(t),0\leq
t\leq\zeta_{S_r})$
is a Bessel process with index $-\frac{5}{2}$ started from $r$ and
stopped when
it hits $0$. By a classical reversal theorem of Williams~\cite{Wi}, Theorem 2.5, the reversed path
$(r+W_{S_r}(\zeta_{S_r}-t),0\leq t\leq\zeta_{S_r})$ is distributed as
a Bessel
process with index~$\frac{5}{2}$, or equivalently with dimension $7$,
started from $0$ and stopped at its last passage
time at level $r$. To simplify notation, we will set
\[
Y^{(r)}_t:= r+W_{S_r}(\zeta_{S_r}-t),\qquad 0\leq t\leq\zeta_{S_r}.
\]
The definition of $Y^{(r)}$ makes sense under $\N_0(\cdot\mid
S_r<\infty)$.

Applying the strong Markov property at $S_r$ will lead us to
consider the Brownian snake ``subexcursions'' branching from
$W_{S_r}$ after time $S_r$ [this really corresponds to the subtrees branching
from the right side of the ancestral line of $p_\be(S_r)$ that were discussed
at the beginning of Section~\ref{mainEstim}, with the difference
that we now argue under the excursion measure]. We consider all nontrivial
subintervals $(v,v')$ of $[S_r,\sigma]$ such that
\[
\zeta_{v}=\zeta_{v'}= \min_{s\in[S_r,v']}
\zeta_s.
\]
We let $(v_i,v'_i)_{i\in I}$ be the collection of all these intervals. For
every $i\in I$ we define a path-valued process $(W^i_s)_{s\geq0}$ by
setting
\[
W^i_s(t) = W_{(v_i+s)\wedge v'_i}(\zeta_{v_i}
+t) - W_{v_i}(\zeta_{v_i}),\qquad 0\leq t\leq\zeta^i_s:=
\zeta_{(v_i+s)\wedge v'_i} -\zeta_{v_i}.
\]
Then, under the probability measure $\N_0(\cdot\mid S_r<\infty)$,
conditionally on $\g_{S_r}$, the point measure
\[
{\cal N}= \sum_{i\in I} \delta_{(\zeta_{v_i},W^i)}(dt\, d
\omega)
\]
is Poisson with intensity $2 \mathbf{1}_{[0,\zeta_{S_r}]}(t)\,dt \N_0(d\omega)$. This follows from
Lemma V.5 in~\cite{Zu} after applying the strong Markov property
of the Brownian snake~\cite{Zu}, Theorem IV.6, at time $S_r$.

\begin{pf*}{Proof of Lemma~\ref{auxil}} We start by explaining how the
bound of the lemma can be reduced to an estimate under the
excursion measure. We write $P$ for the probability measure
$\N_0(\cdot\mid\sigma=1)$. For
every $t<1$, the restriction of $P$ to $\g_t$ is absolutely continuous with
respect to the restriction of $\N_0$ to the same $\sigma$-field. This readily
follows from the analogous property for the law of the normalized Brownian
excursion and the It\^o measure. Moreover, the Radon--Nikodym derivative
of $P_{|\g_t}$ with respect to ${\N_0}_{|\g_t}$ is bounded above by
a constant depending only on $t$.

We then observe that the event
\[
\{ S_{r}\leq1-\kappa\}\cap \Bigl\{ \sup_{s\in[S_{r-\ve},S_r]} \wh
W_s\geq-r+\sqrt{\ve} \Bigr\}
\]
is measurable with respect to $\g_{1-\kappa}$. If we are able to bound
the $\N_0$-measure
of this event, we will immediately get the same bound for its
$P$-measure, up to
a multiplicative constant depending on $\kappa$. So, using the
above-mentioned fact that the
law of $(\zeta_s,\wh W_s)_{0\leq s\leq1}$ under $P$ is the same as
the law
of the process $(\be_s,Z_s)_{0\leq s\leq1}$ of the
preceding
sections, it suffices to verify that the $\N_0$-measure
of the latter event satisfies the bound of Lemma~\ref{auxil}.

To this end, it is enough to prove that
for $r\in[\mu,A]$ and $\ve\in(0,\mu/2)$,
%
%eA.2 #&#
\begin{equation}
\label{auxiltech1} \N_0 \Bigl(S_{r+\ve}<\infty,
\sup_{s\in[S_{r},S_{r+\ve}]} \wh W_s\geq -r+\sqrt{\ve} \Bigr)\leq
C_{A,\mu} \ve^\beta
\end{equation}
with $\beta\in(0,1)$ and a constant $C_{A,\mu}$ depending only on $A$
and $\mu$.
We use the notation introduced at the beginning of this \hyperref[app]{Appendix}, and
we also set
(in this proof only)
\[
T^{(r)}_\ell=\inf\bigl\{t\geq0\dvtx Y^{(r)}_t=2^{-\ell}
\bigr\}
\]
for every integer $\ell\geq0$ such that $2^{-\ell}\leq r$.

For every $\ve\in(0,1)$,
let $\ell_0(\ve)$ and $\ell_1(\ve)$ be the nonnegative integers
such that
\[
2^{-\ell_0(\ve)-1}<\ve\leq2^{-\ell_0(\ve)}, 2^{-\ell_1(\ve
)}<\ve^{3/4}
\leq2^{-\ell_1(\ve)+1}.
\]
Define two events $A_{(\ve)}$ and $B_{(\ve)}$ by
\begin{eqnarray*}
A_{(\ve)}&=&\{S_r<\infty\}\cap \Biggl(\bigcup
_{\ell=\ell_1(\ve
)}^{\ell_0(\ve
)} \Bigl\{ \inf_{\{j\in I\dvtx T^{(r)}_{\ell+1}<\zeta_{S_r}-\zeta_{v_j} <
T^{(r)}_{\ell
}\}} \Bigl(
\inf_{s\geq0} \wh W^j_s \Bigr) < -2
\cdot2^{-\ell} \Bigr\} \Biggr),
\\
B_{(\ve)}&=& \{S_r<\infty\}\cap \biggl\{
\sup_{\{j\in I\dvtx \zeta
_{S_r}-\zeta
_{v_j} < T^{(r)}_{\ell_1(\ve)}\}} \Bigl( \sup_{s\geq0} \wh W^j_s
\Bigr) <\frac{1}{2} \sqrt{\ve} \biggr\}.
\end{eqnarray*}
We may assume that $\ve$ is small enough so that $\ve^{3/4} \leq
\frac
{1}{2}\sqrt{\ve}$. Then if
$B_{(\ve)}$ holds, one immediately checks that $\wh W_s <-r+\sqrt{\ve}$
for $S_r\leq s\leq\inf\{t\geq S_r\dvtx \zeta_t=\zeta_{S_r}-T^{(r)}_{\ell
_1(\ve)}\}$.
On the other hand, if $A_{(\ve)}$ holds, there is a value of $s$ in the
same interval such that $\wh W_s<-r-\ve$. By
combining these observations, we get
\[
(A_{(\ve)}\cap B_{(\ve)}) \subset \Bigl\{S_{r+\ve}<
\infty, \sup_{s\in
[S_{r},S_{r+\ve}]} \wh W_s< -r+\sqrt{\ve} \Bigr\}
\]
and, therefore,
\[
\Bigl\{S_{r+\ve}<\infty, \sup_{s\in[S_{r},S_{r+\ve}]} \wh W_s\geq
-r+\sqrt{\ve} \Bigr\} \subset \bigl( \{S_r<\infty\}\setminus
A_{(\ve)} \bigr) \cup \bigl( \{ S_r<\infty\}\setminus
B_{(\ve)} \bigr).
\]
In view of proving (\ref{auxiltech1}), we bound separately $\N_0(\{
S_r<\infty\}\setminus A_{(\ve)})$
and $\N_0(\{S_r<\infty\}\setminus B_{(\ve)})$.

From the exponential formula for Poisson measures,
and formula (\ref{lawmin}),
we have first
\begin{eqnarray*}
&&\N_0\bigl(\{S_r<\infty\}\setminus B_{(\ve)}
\bigr)
\\
&&\qquad=\N_0(S_r<\infty) \N_0 \bigl(1-\exp
\bigl(-3(\sqrt{\ve }/2)^{-2}T^{(r)}_{\ell_1(\ve)} \bigr) |
S_r<\infty \bigr)
\\
&&\qquad=\N_0(S_r<\infty) E \bigl[1-\exp \bigl(-12
\ve^{-1}2^{-2\ell_1(\ve)} T_{(1)} \bigr) \bigr]
\\
&&\qquad\leq\N_0(S_r<\infty) E \bigl[1-\exp \bigl(-12
\ve^{1/2}T_{(1)} \bigr) \bigr],
\end{eqnarray*}
where, for every $u>0$, $T_{(u)}$ stands
for the hitting time of $u$ by
a seven-dimensional Bessel process started from $0$, and we used the
scaling property $T_{(u)}\build{=}_{}^{\mathrm{(d)}} u^2T_{(1)}$.
Clearly, the right-hand side is bounded above by a constant times $\ve^{1/2}$
(we use the fact that $T_{(1)}$ has moments of
any order).

Then,
\begin{eqnarray*}
&&\N_0\bigl(\{S_r<\infty\}\setminus A_{(\ve)}
\bigr)
\\
&&\qquad=\N_0(S_r<\infty) \N_0 \Biggl(\prod
_{\ell=\ell_1(\ve)}^{\ell
_0(\ve)} \exp \bigl(-3
\bigl(T^{(r)}_\ell- T^{(r)}_{\ell+1}\bigr)
\bigl(2\cdot2^{-\ell
}\bigr)^{-2} \bigr) \Big| S_r<
\infty \Biggr)
\\
&&\qquad=\N_0(S_r<\infty) \prod
_{\ell=\ell_1(\ve)}^{\ell_0(\ve)} E \biggl[ \exp \biggl(-\frac{3}{4}
2^{2\ell}(T_{(2^{-\ell})}-T_{(2^{-\ell-1})}) \biggr) \biggr]
\\
&&\qquad=\N_0(S_r<\infty) E \biggl[ \exp\biggl(-
\frac{3}{4} (T_{(1)}-T_{(1/2)})\biggr)
\biggr]^{\ell_0(\ve)-\ell_1(\ve)+1}
\end{eqnarray*}
using the strong Markov property and the scaling property of the Bessel
process. Since $E[ \exp(-\frac{3}{4} (T_{(1)}-T_{(1/2)})]<1$ and
since $\ell_0(\ve)-\ell_1(\ve)$ behaves like a constant times $\log
(1/\ve)$ when $\ve$ is small, we arrive at the
desired bound. This completes the proof.
\end{pf*}

\begin{pf*}{Proof of Lemma~\ref{snakelemma}}
In this proof $C$ will denote a constant (which may depend on $\mu, A$
and $\kappa$, but not on $r$)
that may vary from line to line.
As explained previously, we can assume that $\be_s=\zeta_s$ and
$Z_s=\wh W_s$,
under the probability measure $P=\N_0(\cdot\mid\sigma=1)$.
We then observe that, for every integer
$\ell\geq\ell_0$, the event $E_\ell$
is measurable with respect to $\g_{1-\kappa/2}$. Thanks to
this observation, it will suffice to prove the bound of Lemma~\ref{snakelemma}
when the expectation is replaced by an integral under $\N_0$. We use
the notation introduced
at the beginning of the \hyperref[app]{Appendix}, and we now set
\[
T^{(r)}_\ell=\inf\bigl\{t\geq0\dvtx Y^{(r)}_t=K^{-\ell}
\bigr\}
\]
for every integer $\ell\geq0$ such that $K^{-\ell}\leq r$. By convention,
we also put \mbox{$T^{(r)}_\infty=0$}. Finally, for every
choice of the integers $k\leq k'\leq\infty$, such that $K^{-k}\leq r$,
we put
\[
I\bigl(k,k'\bigr):=\bigl\{i\in I\dvtx T^{(r)}_{k'}<
\zeta_{S_r}-\zeta_{v_i}< T^{(r)}_k\bigr
\}.
\]
If we view $(\zeta_s)_{0\leq s\leq\sigma}$ as coding a real tree, the
indices $i\in I(k,k')$
correspond to the ``subtrees'' that branch from the ancestral line of
the vertex
corresponding to $S_r$ at a distance between $T^{(r)}_{k'}$ and
$T^{(r)}_{k}$ from this vertex.

Let $E'_\ell$ be the event defined by the same properties as $E_\ell$,
except that we remove the bound $\eta_{K^{-\ell+2}}(r)<1-\frac
{\kappa}{2}$
in (a). Then of course $E_\ell\subset E'_\ell$, and
\[
\Biggl\{ \sum_{k=\lfloor\ell/2\rfloor}^\ell{
\bf1}_{E'_k} \not= \sum_{k=\lfloor\ell/2\rfloor}^\ell{
\bf1}_{E_k} \Biggr\} \subset \biggl( \{S_r<1-\kappa\} \cap
\biggl\{\eta_{K^{-\lfloor\ell
/2\rfloor
+2}}(r)-S_r >\frac{\kappa}{2}\biggr\}
\biggr).
\]
From the
strong Markov property at time $S_r$, we get that the distribution of
$\eta_{K^{-\lfloor\ell/2\rfloor+2}}(r)-S_r$
under $\N_0(\cdot\mid S_r<\infty)$ coincides with the distribution of
the hitting time
of $T^{(r)}_{\lfloor\ell/2\rfloor-2}$ by an independent linear Brownian
motion starting from $0$.
Straightforward estimates now give the bound
\[
\N_0\biggl(\{S_r<1-\kappa\} \cap\biggl\{
\eta_{K^{-\lfloor\ell/2\rfloor+2}}(r)-S_r >\frac{\kappa}{2}\biggr\}\biggr)\leq C
b^\ell,
\]
where the constant $b\in(0,1)$ does not depend on $\ell$ or on $r$.
Hence, the proof
of the lemma reduces to checking the existence of $a'\in(0,1)$ such
that, for $\ell\geq2\ell_0$,
\[
\N_0 \bigl( \mathbf{1}_{\{S_r<1-\kappa\}} a^{\sum_{k=\lfloor\ell
/2\rfloor
}^\ell\mathbf{1}_{E'_k}} \bigr)\leq C
{a'}^\ell.
\]
Finally, thanks to the presence of the indicator function $\mathbf{1}_{\{
S_r<1-\kappa\}}$, we may also
replace $E'_k$ by the event $G_k$, which is defined by the same
properties (a)--(f)
[without the bound $\eta_{K^{-\ell+2}}(r)<1-\frac{\kappa}{2}$ in (a)]
but without imposing that
$S_r\leq1-\kappa$. Then it will be enough to get the bound
%
%eA.3 #&#
\begin{equation}
\label{apptech1} \N_0 \bigl( \mathbf{1}_{\{S_r<\infty\}}
a^{\sum_{k=\lfloor\ell
/2\rfloor
}^\ell\mathbf{1}_{G_k}} \bigr)\leq C {a'}^\ell.
\end{equation}

For every integer $\ell\geq\ell_0$, let $A_\ell$ be the event where
the following holds: There exists
an excursion interval $(v_i,v'_i)$ such that
\begin{longlist}[$(\alpha)$]
\item[$(\alpha)$] $i\in I(\ell-2,\ell-1)$ [or equivalently
$T^{(r)}_{\ell-1} < \zeta_{S_r}-\zeta_{v_i}< T^{(r)}_{\ell-2} $];
\item[$(\beta)$] $-(\alpha_1+\beta_1)K^{-\ell} < \inf\{\wh W^i_s\dvtx s\geq0\} < -(\alpha_2+\beta_2)K^{-\ell} $;
\item[$(\gamma)$] $-r+ \beta_1 K^{-\ell}< \wh W_{v_i} < -r + \beta_2
K^{-\ell} $;
\item[$(\delta)$] $ {\inf_{j\in I(\ell-2,\ell-1)\setminus\{
i\}} (
\inf_{s\geq0}\wh W^j_s )> -\alpha'_2 K^{-\ell}} $;
\item[$(\varepsilon)$] there exists $t\in[\zeta_{S_r}-\zeta_{v_i},
T^{(r)}_{\ell-2}]$ such that $Y^{(r)}_t < \wt\alpha K^{-\ell} $;
\item[$(\varphi)$] $K^{-4\ell} < v'_i-v_i < (1+\lambda) K^{-4\ell} $.
\end{longlist}
Then the events $A_\ell$, $\ell\geq\ell_0$ are independent under
$\N_0(\cdot\mid S_r<\infty)$.
If we condition on $W_{S_r}$ or, equivalently, on the random path
$Y^{(r)}$, this independence
property follows from the independence properties of Poisson measures,
and we can then
use the fact that the processes $(Y^{(r)}_{(T^{(r)}_{\ell-1}+t)\wedge
T^{(r)}_{\ell}})_{t\geq0}$, $\ell\geq\ell_0$
are independent, by the strong Markov property of the Bessel process.
Furthermore, a
scaling argument shows that $\N_0(A_\ell\mid S_r<\infty)=c$, where
$c>0$ is a constant that
does not depend on $\ell$ [notice that the property $\beta_1<\beta_2<4\leq K^2$ ensures that
$(\alpha)$ and $(\gamma)$ are not incompatible].

We also set
\[
B_\ell= \Bigl\{ \inf_{j\in I(\ell-1,\infty)} \Bigl( \inf_{s\geq0} \wh
W^j_s \Bigr) > -\alpha'_2
K^{-\ell} \Bigr\}
\]
and we observe that $A_\ell\cap B_\ell\subset G_\ell$ by
construction. So the bound (\ref{apptech1})
will follow if we can verify that
%
%eA.4 #&#
\begin{equation}
\label{apptech2} \N_0 \bigl( \mathbf{1}_{\{S_r<\infty\}}
a^{\sum_{k=\lfloor\ell
/2\rfloor
}^\ell\mathbf{1}_{A_k\cap B_k}} \bigr)\leq C {a'}^\ell.
\end{equation}
If we had $\mathbf{1}_{A_k}$ instead of $\mathbf{1}_{A_k\cap B_k}$ in (\ref
{apptech2}),
this bound would immediately follow from the independence properties
mentioned above.
The events $A_k\cap B_k$ are not independent, because $B_\ell$ involves
all ``subtrees'' branching above level $\zeta_{S_r}- T^{(r)}_{\ell-1}$,
but still we will prove
that there is enough independence to give a bound of the form (\ref{apptech2}).

To this end, it will be convenient to modify slightly the definition of
$A_\ell$ and $B_\ell$.
We fix an integer $q\geq1$, whose choice will be made precise later,
and we restrict our attention
to integers that are multiples of $q$. Precisely, for every $k\geq
\lfloor\frac{\ell_0}{q}\rfloor+1$, we
let $\wt A_k$ be defined by the same properties as $A_{qk}$,
with the difference that in $(\delta)$ we require
\[
\inf_{j\in I(qk-2,q(k+1)-2)\setminus\{i\}} \Bigl( \inf_{s\geq0}\wh W^j_s
\Bigr)> -\alpha'_2 K^{-qk}.
\]
For the same values of $k$, we put
\[
\wt B_k = \Bigl\{ \inf_{j\in I(q(k+1)-2,\infty)} \Bigl( \inf_{s\geq0}
\wh W^j_s \Bigr) > -\alpha'_2
K^{-qk} \Bigr\}.
\]
It is then immediate that $\wt A_k\cap\wt B_k = A_{qk} \cap B_{qk}$
and so if we can prove that, for
a suitable choice of $q$ and
for every $\ell\geq2(\lfloor\frac{\ell_0}{q}\rfloor+1)$,
%
%eA.5 #&#
\begin{equation}
\label{apptech3} \N_0 \bigl( \mathbf{1}_{\{S_r<\infty\}}
a^{\sum_{k=\lfloor\ell
/2\rfloor
}^\ell\mathbf{1}_{\wt A_k\cap\wt B_k}} \bigr)\leq C {a'}^\ell,
\end{equation}
the bound (\ref{apptech2}) will follow (with a different value of $a'$).

The events $\wt A_k$ are again independent, and (by scaling) they have
the same probability
$c_{(q)}>0$ under $\N_0(\cdot\mid S_r<\infty)$. From a standard large
deviation estimate, we get
\[
\N_0 \Biggl(\sum_{k=\lfloor\ell/2\rfloor}^\ell{
\bf1}_{\wt A_k} <\frac
{c_{(q)}}{4}\ell \Big| S_r<\infty \Biggr)
\leq C \theta_{(q)}^\ell
\]
with some constant $\theta_{(q)}\in(0,1)$. On the other hand, write
$H_\ell$
for the event where we have both
\[
\sum_{k=\lfloor\ell/2\rfloor}^\ell\mathbf{1}_{\wt A_k}
\geq\frac
{c_{(q)}}{4}\ell
\]
and
\[
\sum_{k=\lfloor\ell/2\rfloor}^\ell\mathbf{1}_{\wt B_k}
\geq\ell- \lfloor\ell/2\rfloor-\frac{c_{(q)}}{8}\ell.
\]
On the event $H_\ell$, we have obviously
\[
{\sum_{k=\lfloor\ell/2\rfloor}^\ell\mathbf{1}_{\wt A_k\cap\wt B_k}}
\geq \frac{c_{(q)}}{8}\ell
\]
and, therefore,
\[
\N_0 \bigl( \mathbf{1}_{H_\ell} a^{\sum_{k=\lfloor\ell/2\rfloor}^\ell
\mathbf{
1}_{\wt A_k\cap\wt B_k}} |
S_r<\infty \bigr)\leq a^{c_{(q)}\ell/8}.\vadjust{\goodbreak}
\]
Therefore, the proof of (\ref{apptech3}) will be complete if we can
verify that, for
a suitable choice of $q$,
%
%eA.6 #&#
\begin{equation}
\label{apptech4} \N_0 \Biggl( \{S_r<\infty\} \cap
\Biggl\{ \sum_{k=\lfloor\ell
/2\rfloor
}^\ell{
\bf1}_{\wt B_k^c} \geq\frac{c_{(q)}}{8}\ell \Biggr\} \Biggr) \leq C
{a''}^\ell
\end{equation}
with some $a''\in(0,1)$.

To simplify notation, we write $P_{(r)}$ instead of $\N_0(\cdot\mid
S_r<\infty)$
and $E_{(r)}$ for the associated expectation in what follows.
We start by observing that $c_{(q)}$ is bounded below by a positive
constant $c'$ that does not depend on $q$. Indeed, it follows from
independence properties of Poisson measures that, for every $k\geq
\lfloor\frac{\ell_0}{q}\rfloor+1$,
\begin{eqnarray*}
c_{(q)}&=&P_{(r)}(\wt A_k)
\\
&=&P_{(r)}(A_{qk})\times E_{(r)} \biggl[\exp-2\int
_{T^{(r)}_{q(k+1)-2}}^{T^{(r)}_{qk-1}} \,dt \N_0 \Bigl(
\inf_{s\geq0} \wh W_s > -\alpha'_2K^{-qk}
\Bigr) \biggr]
\\
&=&c E_{(r)}\biggl[\exp\biggl(-2\bigl(T^{(r)}_{qk-1}-T^{(r)}_{q(k+1)-2}
\bigr)\frac
{3}{2{\alpha'_2}^2}K^{2qk}\biggr)\biggr]
\\
&\geq& c E_{(r)}\biggl[\exp\biggl(-\frac{3}{{\alpha'_2}^2}K^{2qk}T^{(r)}_{qk-1}
\biggr)\biggr]
\\
&=& c E\biggl[\exp\biggl(-\frac{3K^2}{{\alpha'_2}^2}T_{(1)}\biggr)\biggr],
\end{eqnarray*}
where as above $T_{(1)}$ stands for the hitting time of $1$ by
a seven-dimensional Bessel process started from $0$, and we
used (\ref{lawmin}) in the second equality. We conclude that
$c_{(q)}\geq c':= cE[\exp(-3K^2({\alpha'}_2)^{-2} T_{(1)})]$. Obviously,
it is enough to prove (\ref{apptech4}) with $c_{(q)}$ replaced by $c'$.

We now specify the choice of $q$. To this end, we first note that, for
any choice of
$\ell_0\leq k < k'\leq\infty$ and $x>0$, we have, with the
convention $T^{(r)}_\infty=0$,
%
%eA.7 #&#
\begin{eqnarray}
\label{apptech5} P_{(r)} \Bigl(\inf_{j\in I(k,k')} \Bigl(
\inf_{s\geq0} \wh W^j_s \Bigr)\leq -x \Bigr)
&=&1- E_{(r)} \biggl[\exp-2\int_{T^{(r)}_{k'}}^{T^{(r)}_{k}}
\frac
{3}{2x^2} \,dt \biggr]
\nonumber\\
&\leq&1 - E_{(r)} \biggl[ \exp\biggl(-\frac{3 T^{(r)}_k}{x^2}\biggr) \biggr]
\nonumber
\\[-8pt]
\\[-8pt]
\nonumber
&=& 1- E \biggl[\exp\biggl(-\frac{3\cdot K^{-2k}}{x^2} T_{(1)}\biggr) \biggr]
\\
&\leq& M \frac{K^{-2k}}{x^2},
\nonumber
\end{eqnarray}
where $M$ is a constant. We choose the integer $q\geq2$ sufficiently
large so that
\[
M \frac{K^{-q+4}}{(\alpha'_2)^2} \leq\frac{1}{2} \quad\mbox{and}\quad q\frac{c'}{8}>1.
\]

Let $\ell\geq\lfloor\frac{2\ell_0}{q}\rfloor+1$. We have
\[
P_{(r)} \Biggl( \sum_{k=\lfloor\ell/2\rfloor}^\ell{
\bf1}_{\wt B_k^c} \geq\frac{c'}{8}\ell \Biggr) \leq\sum
_{k_1,\ldots,k_m} P_{(r)}\bigl( \wt B_{k_1}^c
\cap\cdots\cap\wt B_{k_m}^c\bigr),
\]
where $m=\lfloor\frac{c'}{8} \ell\rfloor$, and the sum in the
right-hand side is over
all choices of $k_1,\ldots,k_m$ such that $\lfloor\ell/2\rfloor\leq
k_1<k_2<\cdots<k_m\leq\ell$.
Obviously, the number of such choices is bounded above by $2^\ell$, and
so the proof of (\ref{apptech4}) will be complete if we can verify
that, for any
choice of $k_1,\ldots,k_m$ as above, we have
%
%eA.8 #&#
\begin{equation}
\label{apptech6} P_{(r)}\bigl( \wt B_{k_1}^c\cap
\cdots\cap\wt B_{k_m}^c\bigr) \leq K^{-qm}
\end{equation}
(recall that $K\geq2$ and $q\frac{c'}{8}>1$). We prove by induction
that the
bound (\ref{apptech6}) holds for any $m\geq1$. If $m=1$, we use the bound
(\ref{apptech5}) with $k'=\infty$, $k$ replaced by $q(k+1)-2$
and $x=\alpha'_2 K^{-qk}$ to get
\[
P_{(r)}\bigl(\wt B_k^c\bigr)\leq M
\frac{K^{-2q+4}}{(\alpha'_2)^2} \leq K^{-q}.
\]
Then, if $m\geq2$,
%
%eA.9 #&#
\begin{eqnarray}
\label{apptech7} &&P_{(r)}\bigl( \wt B_{k_1}^c
\cap\cdots\cap\wt B_{k_m}^c\bigr)
\nonumber\\
&&\qquad\leq P_{(r)} \Bigl(\inf_{j\in I(q(k_1+1)-2,q(k_2+1)-2)} \Bigl( \inf_{s\geq0}
\wh W^j_s \Bigr) \leq -\alpha'_2
K^{-qk_1} \Bigr)
\nonumber
\\[-8pt]
\\[-8pt]
\nonumber
&&\qquad\quad{}\times P_{(r)}\bigl( \wt B_{k_2}^c
\cap\cdots\cap\wt B_{k_m}^c\bigr)
\\
&&\qquad\quad{} +P_{(r)}\bigl( B_{k_2}^{(k_1)}\cap\wt
B_{k_3}^c\cap\cdots\cap\wt B_{k_m}^c
\bigr),
\nonumber
\end{eqnarray}
where, for $\lfloor\ell/2\rfloor\leq k<k'$, we use the notation
\[
B_{k'}^{(k)}= \Bigl\{\inf_{j\in I(q(k'+1)-2,\infty)} \Bigl(
\inf_{s\geq0} \wh W^j_s \Bigr)\leq -
\alpha'_2 K^{-qk} \Bigr\}.
\]
The first term in the right-hand side of (\ref{apptech7}) is bounded by
the quantity
$P_{(r)}(B^{(k_1)}_{k_1})P_{(r)}( \wt B_{k_2}^c\cap\cdots\cap\wt B_{k_m}^c)$.
By iterating the argument, we obtain
\begin{eqnarray*}
P_{(r)}\bigl( \wt B_{k_1}^c\cap\cdots\cap\wt
B_{k_m}^c\bigr) &\leq&\sum_{j=1}^m
P_{(r)}\bigl(B^{(k_1)}_{k_j}\bigr) P_{(r)}
\bigl( \wt B_{k_{j+1}}^c\cap\cdots\cap\wt B_{k_m}^c
\bigr)
\\
&\leq&\sum_{j=1}^m K^{-q(m-j)}
P_{(r)}\bigl(B^{(k_1)}_{k_j}\bigr)
\end{eqnarray*}
using the induction hypothesis. On the other hand, the bound (\ref{apptech5})
gives for $k\leq k'$
\[
P_{(r)}\bigl(B^{(k)}_{k'}\bigr) \leq M
\frac{K^{-2q(k'+1)+4}}{(\alpha'_2)^2K^{-2qk}} \leq\frac{1}{2} K^{-2q(k'-k)-q}
\]
by our choice of $q$. We thus obtain
\begin{eqnarray*}
P_{(r)}\bigl( \wt B_{k_1}^c\cap\cdots\cap\wt
B_{k_m}^c\bigr)&\leq&\frac{1}{2} \sum
_{j=1}^m K^{-q(m-j)} K^{-2q(k_j-k_1)-q}
\\
&\leq&\frac{1}{2} K^{-qm}\sum_{j=1}^m
K^{-q(j-1)}\leq K^{-qm}.
\end{eqnarray*}
This completes the proof of (\ref{apptech6}) and of Lemma \ref
{snakelemma}.
\end{pf*}
\end{appendix}

\section*{Acknowledgments}
I would like to thank Gr\'egory Miermont for
many useful discussions.
I am also grateful to two anonymous referees for their careful reading
and several
helpful observations.

% imsref loaded by akundreckaite, 2013-02-05 12:21:50
%

%

%suskaldyti doi

\printaddresses

\end{document}